\documentclass[11pt, a4paper]{article}

\usepackage[english]{babel}
\usepackage{dutchcal}
\usepackage[T1]{fontenc}

\usepackage{amsmath}
\usepackage{amssymb}
\usepackage{amsthm}
\usepackage{tikz}
\usepackage{tikz-cd}
\usepackage{todonotes}
\usepackage{hyperref}
\usepackage[capitalize]{cleveref}
\usepackage{footmisc}
\usepackage{bbm}
\usepackage{mathrsfs}
\usepackage{extarrows}
\usepackage{enumitem}
\usepackage{stmaryrd}

\usepackage{mathtools}

\newcommand{\llpar}{\llparenthesis}
\newcommand{\rrpar}{\rrparenthesis}

\newcommand{\aff}{\mathrm{aff}}

\newcommand{\fpqc}{{\mathrm{fpqc}}}
\newcommand{\cub}{{\mathrm{cub}}}

\newcommand{\qua}{{\mathrm{quad}}}

\newcommand{\cp}{{\mathbb{C}\mathrm{p}}}

\newcommand{\cn}{\mathrm{cn}}

\newcommand{\Ell}{\mathrm{Ell}}

\newcommand{\Tori}{\mathrm{Tori}}

\newcommand{\et}{\mathrm{\acute{e}t}}

\newcommand{\ori}{\mathrm{or}}

\newcommand{\st}{\mathrm{st}}

\renewcommand{\top}{\mathrm{top}}

\newcommand{\1}{\mathbf{1}}

\DeclareMathOperator{\KU}{KU}

\DeclareMathOperator{\MOP}{MOP}

\DeclareMathOperator{\KO}{KO}

\DeclareMathOperator{\ko}{ko}

\newcommand{\MU}{\mathrm{MU}}
\newcommand{\MO}{\mathrm{MO}}
\newcommand{\MUP}{\mathrm{MP}}
\newcommand{\MP}{\mathrm{MP}}

\newcommand{\Sph}{\mathbf{S}}

\DeclareMathOperator{\TMF}{TMF}
\DeclareMathOperator{\tmf}{tmf}
\DeclareMathOperator{\Tmf}{Tmf}

\DeclareMathOperator{\Aff}{Aff}

\DeclareMathOperator{\CAlg}{CAlg}

\DeclareMathOperator{\Cat}{Cat}

\DeclareMathOperator{\Fun}{Fun}

\newcommand{\FG}{\mathrm{FG}}
\newcommand{\FGh}{{\smash{{\mathrm{FG}{\leq h}}}}}
\newcommand{\FGhp}{{\smash{{\mathrm{FG}{\leq \vec{h}(p)}}}}}
\newcommand{\FGvech}{{\smash{{\mathrm{FG}{\leq \vec{h}}}}}}
\newcommand{\FGvechgeq}{{\smash{{\mathrm{FG}{\geq \vec{h}}}}}}

\DeclareMathOperator{\Mod}{Mod}
\DeclareMathOperator{\LMod}{LMod}

\renewcommand{\Pr}{\mathrm{Pr}}
\newcommand{\PrLst}{\mathrm{Pr}^L_{\mathrm{st}}}

\DeclareMathOperator{\QCoh}{QCoh}

\DeclareMathOperator{\Spec}{Spec}
\DeclareMathOperator{\Spet}{Sp\acute et}

\newcommand{\relSpec}{\underline{\mathrm{Spec}}}
\DeclareMathOperator{\Shv}{Shv}
\DeclareMathOperator{\Stk}{Stk}

\DeclareMathOperator{\GeoStk}{GeoStk}
\DeclareMathOperator{\Sp}{Sp}
\DeclareMathOperator{\SpDM}{SpDM}
\DeclareMathOperator{\SpDMnc}{SpDM^{nc}}
\DeclareMathOperator{\SpDMncet}{SpDM^{nc}_{\et}}

\newcommand{\Spc}{\mathcal{S}}
\newcommand{\spaces}{\mathcal{S}}

\DeclareMathOperator{\Tot}{Tot}

\DeclareMathOperator{\colim}{colim}
\DeclareMathOperator{\fib}{fib}
\DeclareMathOperator{\cof}{cof}

\newcommand{\id}{\mathrm{id}}
\DeclareMathOperator{\Map}{Map}
\DeclareMathOperator{\map}{map}

\DeclareMathOperator{\Hom}{Map}

\newcommand{\op}{\mathrm{op}}

\newcommand{\A}{\mathsf{A}}
\newcommand{\calA}{\mathcal{A}}
\newcommand{\B}{\mathsf{B}}

\newcommand{\C}{\mathcal{C}}
\newcommand{\csite}{\mathcal{C}}
\newcommand{\CP}{\mathbf{CP}}
\newcommand{\D}{\mathsf{D}}
\newcommand{\calD}{\mathcal{D}}
\newcommand{\E}{\mathbf{E}}

\newcommand{\F}{\mathbf{F}}

\newcommand{\G}{\mathbf{G}}

\newcommand{\sfG}{\mathsf{G}}

\newcommand{\M}{\mathsf{M}}

\newcommand{\calM}{\mathcal{M}}
\newcommand{\calN}{\mathcal{N}}
\newcommand{\N}{\mathbf{N}}

\renewcommand{\O}{\mathcal{O}}

\renewcommand{\P}{\mathbf{P}}

\newcommand{\QQ}{\mathcal{Q}}
\newcommand{\R}{\mathsf{R}}
\renewcommand{\S}{\mathsf{X}}

\newcommand{\T}{\mathrm{T}}

\newcommand{\sfU}{\mathsf{U}}

\newcommand{\sfV}{\mathsf{V}}
\newcommand{\X}{\mathsf{X}}
\newcommand{\Y}{\mathsf{Y}}

\newcommand{\Z}{\mathbf{Z}}

\newcommand{\sfZ}{\mathsf{Z}}
\newcommand{\sfX}{\mathsf{X}}
\newcommand{\sfY}{\mathsf{Y}}
\newcommand{\calO}{\mathcal{O}}
\newcommand{\sfC}{\mathsf{C}}
\newcommand{\calC}{\mathcal{C}}
\newcommand{\sfD}{\mathsf{D}}

\newcommand{\sfA}{\mathsf{A}}
\newcommand{\sfF}{\mathsf{F}}
\newcommand{\sfQ}{\mathsf{Q}}
\newcommand{\sfW}{\mathsf{W}}
\newcommand{\sfB}{\mathsf{B}}

\newcommand{\nc}{\mathrm{nc}}

\newcommand{\dcat}{\mathsf{D}}
\newcommand{\ccat}{\mathsf{C}}

\newcommand{\jdiagram}{\mathsf{J}}

\newcommand{\spherespectrum}{\Sph}
\newcommand{\monoidalunit}{\1}
\newcommand{\cstablecat}{\mathsf{C}}
\newcommand{\cmonoidalunit}{\1_\cstablecat}
\newcommand{\dmonoidalunit}{\1_\dstablecat}
\newcommand{\dstablecat}{\mathsf{D}}
\newcommand{\Thick}{\mathrm{Thick}}

\newcommand{\ol}[1]{\overline{#1}}

\newcommand{\bs}{{-}}

\newcommand{\presheaves}{\mathcal{P}}

\newcommand{\al}{\alpha}
\newcommand{\be}{\beta}
\newcommand{\ga}{\gamma}
\newcommand{\Ga}{\Gamma}

\newcommand\noloc{%
	\nobreak
	\mspace{6mu plus 1mu}
	{:}
	\nonscript\mkern-\thinmuskip
	\mathpunct{}
	\mspace{2mu}
}

\theoremstyle{theorem}\numberwithin{equation}{subsubsection}
\newtheorem{theorem}[equation]{Theorem}
\Crefname{theorem}{{Th}.\!\!}{{Ths}.\!\!}
\newtheorem{theoremalph}{Theorem}

\Crefname{theoremalph}{{Th}.\!\!}{{Ths}.\!\!}

\Crefname{coralph}{{Cor}.\!\!}{{Cors}.\!\!}

\Crefname{conjalph}{{Conj}.\!\!}{{Conjs}.\!\!}

\Crefname{problem}{{Prb}.\!\!}{{Prbs}.\!\!}
\newtheorem{prop}[equation]{Proposition}
\Crefname{prop}{{Pr}.\!\!}{{Prs}.\!\!}
\newtheorem{lemma}[equation]{Lemma}
\Crefname{lemma}{{Lm}.\!\!}{{Lms}.\!\!}
\newtheorem{cor}[equation]{Corollary}
\Crefname{cor}{{Cor}.\!\!}{{Cors}.\!\!}

\Crefname{conjecture}{{Conj}.\!\!}{{Conjs}.\!\!}

\theoremstyle{definition}\numberwithin{equation}{subsubsection}
\newtheorem{mydef}[equation]{Definition}
\Crefname{mydef}{{Df}.\!\!}{{Dfs}.\!\!}

\Crefname{recall}{{Rcl}.\!\!}{{Rcls}.\!\!}

\Crefname{construction}{{Con}.\!\!}{{Cons}.\!\!}

\Crefname{ass}{{As}.\!\!}{{As}.\!\!}

\Crefname{notation}{{Nt}.\!\!}{{Nts}.\!\!}

\Crefname{situation}{{St}.\!\!}{{Sts}.\!\!}

\theoremstyle{remark}\numberwithin{equation}{subsubsection}
\newtheorem{example}[equation]{Example}
\Crefname{example}{{Ex}.\!\!}{{Exs}.\!\!}
\newtheorem{nonexample}[equation]{Non-example}
\Crefname{nonexample}{{NonEx}.\!\!}{{NonEx}.\!\!}

\Crefname{claim}{{Clm}.\!\!}{{Clms}.\!\!}
\newtheorem{remark}[equation]{Remark}
\Crefname{remark}{{Rmk}.\!\!}{{Rmks}.\!\!}

\Crefname{idea}{{Id}.\!\!}{{Ids}.\!\!}
\newtheorem{warn}[equation]{Warning}
\Crefname{warn}{{Warn}.\!\!}{{Warns}.\!\!}

\Crefname{question}{{Qn}.\!\!}{{Qns}.\!\!}

\Crefname{figure}{{Fig.}\!\!}{{Figs.}\!\!}
\Crefname{footnote}{{Fn.}\!\!}{{Fn.}\!\!}

\Crefname{part}{{\textsection}\!\!}{{\textsection}\!\!}
\Crefname{chapter}{{\textsection}\!\!}{{\textsection}\!\!}
\Crefname{section}{{\textsection}\!\!}{{\textsection}\!\!}
\Crefname{subsection}{{\textsection}\!\!}{{\textsection}\!\!}
\Crefname{appendix}{{\textsection}\!\!}{{\textsection}\!\!}

\usepackage{parskip}

\begingroup
\makeatletter
\@for\theoremstyle:=definition,remark,plain,theorem\do{%
	\expandafter\g@addto@macro\csname th@\theoremstyle\endcsname{%
		\addtolength\thm@preskip\parskip
	}%
}
\endgroup

\begin{document}
	\title{Affineness and reconstruction\\ in complex-periodic geometry}
	\author{William Balderrama\footnote{\url{williamb@math.uni-bonn.de}}, Jack Morgan Davies\footnote{\url{davies@uni-wuppertal.de}}, and Sil Linskens\footnote{\url{sil.linskens@mathematik.uni-regensburg.de}}}
    \date{November 27, 2025}
	\maketitle

	\begin{abstract}
		Working in a generic derived algebro-geometric context, we lay the foundations for the general study of affineness and local descendability. When applied to $\E_\infty$ rings equipped with the fpqc topology, these foundations give an $\infty$-category of \emph{spectral stacks}, a viable functor-of-points alternative to Lurie's approach to nonconnective spectral algebraic geometry in \cite{sag}. Specializing further to spectral stacks over the moduli stack of oriented formal groups, we use chromatic homotopy theory to obtain a large class of $0$-affine stacks, generalizing Mathew--Meier's famous $0$-affineness result. We introduce a spectral refinement of Hopkins' stack construction of an $\E_\infty$ ring, and study when it provides an inverse to the global sections of a spectral stack. We use this to show that a large class of stacks, which we call \emph{reconstructible}, are naturally determined by their global sections, including moduli stacks of oriented formal groups of bounded height and the moduli stack of oriented elliptic curves.
	\end{abstract}
	\setcounter{tocdepth}{3}
	
	{\small
		\tableofcontents}
        
	\newpage
	\section{Introduction}
	Let $\X$ be a derived stack: roughly, a pair $(X,\O_\X)$ of a topological space $X$ and a sheaf $\O_\X$ on $X$ taking values in commutative algebras in a stable symmetric monoidal $\infty$-category $\calC$.\footnote{Our actual generalized stacks are defined in \Cref{ssec:stackssection} and take a functor-of-points perspective, rather than working with sheaves on spaces or $\infty$-topoi.} Associated with $\X$ are algebraic and categorical invariants, such as its derived ring of global sections $\Ga(\O_\X)$ and its $\infty$-category of quasi-coherent sheaves $\QCoh(\X)$. The titular concepts of \emph{reconstruction} and ($0$-)\emph{affineness} are concerned with when $\X$ can be recovered from $\Ga(\O_\X)$, and when $\QCoh(\X)$ can be recovered from $\Mod_{\Ga(\O_\X)}(\calC)$, respectively.\footnote{The vertical relationships in (\ref{eq:diagramintro}), between $\X$ and $\QCoh(\X)$, and $\Ga(\O_\X)$ and $\Mod_{\Ga(\O_\X)}(\csite)$ are more formal. The first is the fact that the assignment $\X$ to the symmetric monoidal $\infty$-category $\QCoh(\X)$ is often fully faithful, often called \emph{Tannaka duality}, discussed in \Cref{sec:appendix}. The second is the affine version of this fact, as shown by Lurie in \cite[Cor.4.8.5.21]{ha}.}
	
	\begin{equation}\label{eq:diagramintro}\begin{tikzcd}
			{\X} &&&& {\Gamma(\O_\X)} \\
			{\QCoh(\X)} &&&& {\Mod_{\Gamma(\O_\X)}(\calC)}
			\arrow["{\text{Reconstruction, \Cref{maintheorem:reconstruction,maintheorem:reconstruction2}}}", dashed, from=1-1, to=1-5]
			\arrow["{\text{Tannaka duality, \Cref{sec:appendix}}}"', dashed, from=1-1, to=2-1]
			\arrow["{\text{\cite[\textsection4.8]{ha}}}"', 
			dashed, from=1-5, to=2-5]
			\arrow["{\text{0-affineness, \Cref{maintheorem:boundedaffineness,maintheorem:moaoriaffineness}}}"', dashed, from=2-1, to=2-5]
	\end{tikzcd}\end{equation}
	
	It is perhaps surprising that these are interesting notions to study. For example, the existence of \emph{$0$-affine} derived stacks, so those $\X$ with an equivalence $\QCoh(\X) \simeq \Mod_{\Ga(\O_\X)}(\calC)$, is a purely \emph{derived phenonemon}: the classical version of this definition, where $\X$ is replaced with a classical quasi-compact scheme $X$ and stable $\infty$-categories with abelian $1$-categories, is uninteresting, as Serre's famous affineness criterion of \cite{serreaffineness} immediately shows that such a scheme $X$ is affine. However, in the derived setting, there are many interesting $0$-affine stacks which are not themselves affine.
	
	For example, consider the world of \emph{spectral algebraic geometry}, the case where $\calC$ is the universal stable symmetric monoidal $\infty$-category $\Sp$ of \emph{spectra}.
	Variants of affineness in this context have been explored by To\"{e}n \cite{affinestackstoen}, Lurie \cite[\textsection4]{dagviii} and \cite{sag}, and Mathew--Meier \cite{akhilandlennart}. For example, Lurie shows \cite[Pr.2.4.1.4]{sag} that if $\X$ is a \emph{quasi-affine} nonconnective spectral Deligne--Mumford stack, then $\X$ is $0$-affine. Mathew--Meier prove a \emph{relative} version of this statement over $\M_\FG^\ori$, the \emph{moduli stack of oriented formal groups}: if $\X$ is a Noetherian nonconnective spectral Deligne--Mumford stack which admits a flat quasi-affine map to $\M_\FG^\ori$, then $\X$ is $0$-affine. This result established a deep connection between the notion of $0$-affineness and the study of chromatic homotopy theory. The $0$-affineness of the \emph{moduli stack of oriented elliptic curves} $\M_\Ell^\ori$, for instance, is key to our understanding of the $\E_\infty$ ring $\TMF$ of \emph{topological modular forms}, with applications to its Galois theory \cite[\textsection10]{akhilgalois} and associated descent spectral sequence \cite{osyn}.

\subsection{Main results}
In this article, we strengthen and clarify the connection between $0$-affineness and chromatic homotopy theory. Along the way, we develop a convenient foundation for non-connective spectral algebraic geometry particularly suited to studying properties of quasi-coherent sheaves. Our first goal is to generalize the affineness result of Mathew--Meier to include stacks which are as close to $\M_\FG^\ori$ as possible. To do this, we have to work in a larger $\infty$-category than that of nonconnective spectral Deligne--Mumford stacks, as $\M_\FG^\ori$ is not itself Deligne--Mumford---it only admits a \emph{flat} cover by an affine stack, as opposed to an \'etale cover. Unfortunately, there are delicate size issues with the $\infty$-category of flat sheaves; already classically, flat sheafification need not exist within a given universe, see \cite{waterhousefpqc}. To sidestep this issue, we set up general foundations for studying derived algebraic geometry built from sheaves on a possibly large Grothendieck site; see \Cref{sec:sheaves}. Applying this setup to $\CAlg$, the $\infty$-category of commutative algebras in $\Sp$, or \emph{$\E_\infty$ rings}, equipped with the flat topology, we obtain the $\infty$-category $\Stk$ of \emph{spectral stacks}; see \Cref{sec:nonconsag}. Finally, a \emph{complex-periodic} stack is a spectral stack $\X$ equipped with a map $\X \to \M_\FG^\ori$. As we explain in \Cref{ssec:MFGOR}, a spectral stack $\X$ admits a map to $\M_\FG^\ori$ if and only if $\X$ admits a flat cover by \emph{complex-periodic} affines $\Spec \A_i$, in which case the resulting map to $\M_\FG^\ori$ is unique, explaining our terminology.
	
	\begin{theoremalph}[$0$-affineness]
		\label{maintheorem:boundedaffineness}
		Let $\X$ be a complex-periodic stack of bounded height. If the unique morphism $\X \to \M_\FG^\ori$ is quasi-affine, then $\X$ is \emph{$0$-affine}.
	\end{theoremalph}

In fact, we further generalize this statement by replacing the geometric quasi-affineness condition with the more categorical condition we call \emph{universally $0$-affine}; we will come back to this point in \Cref{introdef:u0a}.

Intuitively, the ``bounded height'' condition means that we restrict to complex-periodic stacks whose associated formal group has finite height. The theorem above is general enough to capture the universal examples of stacks of bounded height, given by \smash{$\X = \M_{\smash{\FGvech}}^\ori$}, the moduli stack of oriented formal groups of \emph{height $\leq \vec{h}$}, where $\vec{h}$ is some bounded height function on the set of prime numbers.
	
	\Cref{maintheorem:boundedaffineness} generalizes the main theorem of \cite{akhilandlennart}. Achieving this generalization requires a different proof, which is considerably more conceptual and categorical. To emphasize this, we present our new simple proof, essentially in its entirety, in \cref{ssec:proofofMM}. The length of this paper reflects our attempt to distill this categorical argument into its most basic form, with potential applications to other variants of derived algebraic geometry. As an additional benefit, the proof we give also naturally leads one to interesting variations of the result above.
	
	\begin{theoremalph}[$0$-affineness of complex-periodifications]
		\label{maintheorem:moaoriaffineness}
		Let $\A$ be an $\MU$-nilpotent $\E_\infty$-ring. Then the \emph{complex-periodification of $\A$}, given by the stack $\M_\A^\ori = \Spec \A \times \M_\FG^\ori$, is $0$-affine.
	\end{theoremalph}
	
	The stack $\M_{\A}^\ori$ is a complex-periodic lift of Hopkins' classical stack associated with the ring spectrum $\A$ found in \cite{coctalos,handbooktmf}. By design, the descent spectral sequence for $\M_\A^\ori$ is the Adams--Novikov spectral sequence for $\A$. \Cref{maintheorem:moaoriaffineness} applies to $\A=\ko$ and $\tmf$, for example; this produces stacks which fall outside the more classical purview of \Cref{maintheorem:boundedaffineness}, as both $\ko$ and $\tmf$ have infinite chromatic height. The spectral stack $\M_{\ko}^\ori$ is a complex-periodic lift of the \emph{moduli stack of quadratic equations}, and $\M_{\tmf}^\ori$ is a complex-periodic lift of the \emph{moduli stack of Weierstraß curves} or \emph{cubic equations}. These complex-periodic stacks are not flat over $\M_\FG^\ori$, and therefore cannot be constructed using any spectral lift of the Landweber exact functor theory, such as the Goerss--Hopkins--Miller theorem.
	
	Just as $0$-affineness compares quasi-coherent sheaves over a stack with modules over its global sections, the functor $\M_{(-)}^\ori$ from $\E_\infty$ rings to complex-periodic stacks can be used to directly compare a stack $\X$ with the complex-periodification of its global sections.
	
	\begin{theoremalph}[Reconstruction]\label{maintheorem:reconstruction}
		Let $\X$ be a $0$-affine complex-periodic stack for which the unique map $\X \to \M_\FG^\ori$ is affine. Then $\X$ is \emph{reconstructible}: the canonical map $\X \to \M_{\Ga(\O_\X)}^\ori$ is an equivalence.
	\end{theoremalph}
	
	This general reconstructibility result captures all of the examples of complex-periodification explored by Gregoric in \cite[Th.I.2]{rokeven}. For example, we obtain equivalences of complex-periodic stacks
	\begin{equation}\label{eq:examplesofreconstructiblestackies}\M_{\FGvech}^\ori \simeq \M_{\Sph_{\vec{h}}}^\ori, \qquad \M_\Ell^\ori \simeq \M_{\TMF}^\ori, \qquad \M_{\Tori}^\ori \simeq \M_{\KO}^\ori,\end{equation}
	where $\vec{h}$ is a bounded height function and $\Sph_{\vec{h}}$ is the $L_{\vec{h}}$-local sphere, and $\M_{\Tori}^\ori = \Spec \KU / C_2$ is the \emph{moduli stack of oriented tori} with $\KO$ its $\E_\infty$ ring of global sections. See \Cref{sssec:exampleofreconstruction} for more details and examples.

     As an application of reconstruction, we show in \Cref{cor:orientedeillipeiticcurves} that for a complex-periodic stack $\X$, the $\infty$-groupoid of \emph{oriented elliptic curves} over $\X$ is naturally equivalent to maps of $\E_\infty$ rings from $\TMF$ to $\Ga(\O_\X)$:
     \[\Ell^\ori(\X) \simeq \Map_{\Stk}(\X, \M_\Ell^\ori) \xrightarrow{\simeq} \Map_{\CAlg}(\TMF, \Ga(\O_\X)\]
     For instance, this implies that an oriented elliptic curve is determined by its rationalisation, all of its $p$-completions, together with some adelic gluing data. The second-named authors used this fact in \cite{globaltate} to construct a derived Tate curve.
	
	The definition of a reconstuctible stack immediately implies that the global sections functor $\Ga\colon \Stk_{/\M_\FG^\ori} \to \CAlg^\op$ is fully faithful when restricted to reconstructible stacks; see \Cref{prop:ffforreconstruction}. If we further restrict to those stacks of a particular bounded height, we can also identify the essential image of $\Ga$.
	
	\begin{theoremalph}[Bounded reconstruction]\label{maintheorem:reconstruction2}
		Let $\vec{h}$ be a bounded height function. Then the functor of global sections
		\[\Ga \colon \Stk_{/\M_\FGvech^\ori}^\aff \xrightarrow{\simeq} \CAlg(\Sp_{\vec{h}})^\op,\]
		from complex-periodic stacks $\X$ of height $\leq \vec{h}$ for which $\X \to \M_\FGvech^\ori$ is affine to $L_{\vec{h}}$-local $\E_\infty$ rings, is an equivalence.\footnote{Given a height function $\vec{h}$, the category $\Sp_{\vec{h}}$ is defined as the Bousfield localization of $\Sp$ at $\prod_p E_{\vec{h}(p)}$, a product of Morava $E$-theories over all primes of height dictated by $\vec{h}$. We use this notation to highlight that we need not focus on a single prime $p$, allowing for integral examples such as $\KO$ and $\TMF$.}
	\end{theoremalph}
	
	Reconstructible stacks and $0$-affine stacks form different subcategories of spectral stacks, and neither subcategory contains the other. For example, $\M_\FG^\ori = \M_{\Sph}^\ori$ is clearly reconstructible, but it is not $0$-affine as $\Ga\colon \QCoh(\M_\FG^\ori) \to \Sp$ is not an equivalence. On the other hand, the \emph{compactification of the moduli stack of oriented elliptic curves} $\overline{\M}_\Ell^\ori$, as defined in \cite{globaltate}, is $0$-affine, but not reconstructible. Indeed, $0$-affineness is already contained in \cite{akhilandlennart}. For reconstructibility, first notice that all reconstructible stacks are \emph{affine} over $\M_\FG^\ori$ by construction, and $\overline{\M}_\Ell^\ori$ is merely quasi-affine over $\M_\FG^\ori$. One can also read this discrepancy off of the fact that the Adams--Novikov spectral sequence for $\Tmf = \Ga(\O_{\overline{\M}_\Ell^\ori})$ does not agree with the descent spectral sequence for $\overline{\M}_\Ell^\ori$; see \cite[\textsection1.2.1]{smfcomputation} for more on this difference. The complex-periodic stacks of (\ref{eq:examplesofreconstructiblestackies}) happen to all be both $0$-affine and reconstructible.
	
	\subsection{A proof of \Cref{maintheorem:boundedaffineness}}\label{ssec:proofofMM}
	As promised, we now sketch the proof of \Cref{maintheorem:boundedaffineness}. In particular, this will highlight the roles played by notions such as ``universally $0$-affine'' and ``bounded height''.
		
	\begin{proof}[Proof of \Cref{maintheorem:boundedaffineness}]
		The $\infty$-category $\Mod_{\Ga(\O_\X)}$ is determined as a symmetric monoidal stable $\infty$-category by the fact that its unit $\1$ is a compact generator and that the $\E_\infty$-ring of endomorphisms of $\1$ is precisely $\Ga(\O_\X)$. The latter is true for $\QCoh(\X)$ by design, so it suffices to show that $\O_\X$ is a compact generator of $\QCoh(\X)$. The global sections functor $\Ga \colon \QCoh(\X) \to \Sp$ is corepresented by $\O_\X$, so it further suffices to show this functor is conservative and cocontinuous. By assumption, $\X$ is of bounded height, so $f$ factors through the open substack $\M_\FGvech^\ori \subseteq \M_\FG^\ori$ of oriented formal groups of height $\leq \vec{h}$ for some bounded height function $\vec{h}$. The global sections functor $\Ga$ can be written as the composite of pushforward functors
		\[\Ga\colon \QCoh(\X) \xrightarrow{f_\ast} \QCoh(\M_\FGvech^\ori) \xrightarrow{\pi_\ast} \QCoh(\Spec \Sph) \simeq \Sp,\]
		hence it suffices to show that $f_\ast$ and $\pi_\ast$ are both conservative and cocontinuous.
		
		\paragraph{Conservativity of $f_\ast$:} Let $\calM \in \QCoh(\X)$ be such that $f_\ast \calM =0$. We begin with the observation that $\M_\FGvech^\ori$ is covered by $\Spec E$, where $E$ is any Landweber exact $\E_\infty$-ring with the correct height. We can then consider the pullback of stacks
		\[\begin{tikzcd}
			{\Y}\ar[r, "{f'}"]\ar[d, "{p'}"]    &   {\Spec E}\ar[d, "p"]    \\
			{\X}\ar[r, "f"] &   {\M_\FGvech^\ori.}
		\end{tikzcd}\]
		As $f$ is quasi-affine, we have a Beck--Chevalley equivalence of $E$-modules $p^\ast f_\ast \calM \xrightarrow{\simeq} f'_\ast p'^\ast \calM$, so in particular, $f'_\ast p'^\ast \calM$ vanishes. The definition of $f$ being quasi-affine means that $f'$ can be written as the composite of an open immersion $f''$ followed by a map between affines $f'''$:
		\[f' \colon \Y \xrightarrow{f''} \Spec \Ga(\O_\Y) \xrightarrow{f'''} \Spec E.\]
		Pushforwards of open immersions are conservative, being inclusions of local subcategories
		, and pushforwards between affines are clearly conservative; they are just restriction of scalars. In particular, $f'_\ast = f'''_\ast f''_\ast$ is conservative, so $p'^\ast \calM$ vanishes. However, $p'$ is again a cover by base change, so $p'^\ast$ is also conservative, hence $\calM = 0$, proving that $f_\ast$ is conservative.
		
		\paragraph{Cocontinuity of $f_\ast$:} Given a small diagram $\calM_i$ in $\QCoh(\X)$, it suffices to show that the natural map $\colim f_\ast \calM_i \to f_\ast \colim \calM_i$ is an equivalence after applying $p^\ast$, as this pullback is conservative. Using similar arguments from the above paragraph, and the facts that $f''_\ast$, $f'''_\ast$, and pullbacks are all cocontinuous, we obtain the desired equivalences:
		\[p^\ast \colim f_\ast \calM_i \simeq \colim p^\ast f_\ast \calM_i \simeq \colim f'_\ast p'^\ast \calM_i \simeq f'_\ast p^\ast \colim \calM_i \simeq p^\ast f_\ast \colim \calM_i\]
		
		\paragraph{The functor $\pi_\ast$:} The conservativity and cocontinuity of $\pi_\ast$ follows from a different collection of facts. Combining the cover of $\M_{\FGvech}^{\ori}$ by $\Spec E$ with a refinement of the \emph{Hopkins--Ravenel smash product theorem} we obtain an equivalence
		\[\QCoh(\M_\FGvech^\ori) \simeq \lim \Mod_{E^{\otimes(\bullet+1)}} \simeq L_E \Sp\]
		between the category of quasi-coherent sheaves on $\M_\FGvech^\ori$ and the category of $E$-local spectra. In particular, $\pi_\ast$ can be identified with the inclusion of $E$-local spectra into $\Sp$. This inclusion is clearly conservative, and another application of the smash product theorem shows that it also preserves colimits. This finishes the proof.
	\end{proof}

	Note that this proof is straightforward and made completely internally to the language of spectral algebraic geometry and higher algebra; no sheafifications of homotopy groups or explicit invocations of spectral sequences are necessary. Instead, the computational ingredients underlying the proof are entirely contained in our use of the smash product theorem.
	
    Also note that the proof above is divided into proving the cocontinuity and conservativity of $f_\ast$ and of $\pi_\ast$ separately. Moreover, to prove these facts about $f_\ast$, we only needed them for its base change $f'_\ast$. This leads to the following definition.
	
	\begin{mydef}[{\Cref{def:univ0affine}}]\label{introdef:u0a}
		A map of spectral stacks $f\colon \Y \to \X$ is \emph{universally $0$-affine} if for each $\Spec \A \to \X$, the pullback $f'\colon \Spec \A\times_\X \Y = \Y_\A \to \Spec \A$ induces a cocontinuous and conservative pushforward functor $f'_\ast$, and for each $\A$-algebra $\A'$, the natural map $\Ga(\O_{\Y_\A}) \otimes_\A \A' \to \Ga(\Y_{\A'})$ is an equivalence.
	\end{mydef}
	
	The second condition is akin to an \emph{adjointability} or \emph{Beck--Chevalley condition}, and is necessary for $f_\ast$ to commute with pullbacks. It roughly follows from the arguments made above that quasi-affine morphisms are universally $0$-affine; see \Cref{ex:quasiaffineareU0A} for details. Moreover, replacing ``quasi-affine'' with ``universally $0$-affine'' in the proof of \Cref{maintheorem:boundedaffineness} above then yields the following generalization.
	
	\begin{cor}[{\Cref{thm:affinenessbounded}}]\label{cor:toboundedaffinenessintro}
		Let $\X$ be a complex-periodic stack of bounded height. If the unique morphism $\X \to \M_\FG^\ori$ is universally $0$-affine, then $\X$ is $0$-affine.
	\end{cor}
	
	The behaviour of $\pi_\ast$ and $\M_\FGvech^\ori$ also generalizes. The expression for $\M_\FGvech^\ori$ as the colimit of affines $\colim \Spec E^{\otimes(\bullet+1)}$ exhibits $\M_\FGvech^\ori$ as the \emph{$(-1)$-truncation} of the unique map $\Spec E \to \Spec \Sph$. This leads to the following definition.
	
	\begin{mydef}[{\Cref{def:descentstacks}}]\label{introdef:descentstack}
		For a map of stacks $f\colon \Y \to \X$, we define the \emph{descent stack} $\sfD_f$ of $f$ as its $(-1)$-truncation, so that $\sfD_f$ is equivalent to the geometric realization of the \v{C}ech nerve of $f$:
		\[
		\sfD_f = \colim \left(\begin{tikzcd} \Y&\Y\times_\X \Y\ar[l,shift left=1mm]\ar[l,shift right=1mm]&\cdots\ar[l,shift left=1mm]\ar[l]\ar[l,shift right=1mm]\end{tikzcd}\right).
		\]
		If $\Y=\Spec \Sph$ is the terminal stack, we write $\D_f = \D_\X$, and if moreover $\X=\Spec \A$ is affine then we write $\D_f = \D_\X = \D_\A$.
	\end{mydef}
	
	One can then show that $\D_E \simeq \M_\FGvech^\ori$ for the $E$ described in the proof of \Cref{maintheorem:boundedaffineness} (\Cref{pointsofleqvech}), and that $\M_\FG^\ori \simeq \D_{\MUP}$, where $\MUP$ is any complex-periodic form of $\MU$ (\Cref{identicationofMFGOR}).
	
	For a map of stacks $f\colon \Y \to \X$, there is always a natural inclusion $j\colon \D_f \to \X$ inducing a pushforward $j_\ast \colon\QCoh(\D_f) \to \QCoh(\X)$. For this functor to behave as $\D_E \simeq \M_\FGvech^\ori$ does in the proof of \Cref{maintheorem:boundedaffineness}, we need one further assumption on $f$, originally studied by Balmer \cite{balmer2016separable}, Mathew \cite{akhilgalois}, and the first-named author \cite[\textsection A]{balderrama2024total}.
	
	\begin{mydef}[{\Cref{def:ldgeneral,def:ld}}]\label{introdef:ld}
		An object $X$ in a stable symmetric monoidal $\infty$-category $\sfC$ is \emph{locally descendable} if the Bousfield localization $L_X \1$ lies in the thick tensor ideal generated by $X$. 
A map $f\colon \Y \to \X$ of stacks is \emph{locally descendable} if $f$ is universally $0$-affine and $f_\ast \O_\Y \in \QCoh(\X)$ is locally descendable.
	\end{mydef}
	
	We show in \Cref{thm:locallydescendable} that if $f$ is locally descendable then $j_\ast$ induces an equivalence
	\[\QCoh(\D_f) \simeq L_{f_\ast \O_\Y} (\QCoh(\X)),\]
	where we take $f_\ast \O_\Y$-local objects in $\QCoh(\X)$ on the right-hand side. One form of the \emph{smash product theorem} shows that $\Spec E \to \Spec \Sph$ is a locally descendable map of stacks, which then yields the symmetric monoidal identification
	\begin{equation}\label{eq:introdescentableetheorey}\QCoh(\M_\FGvech^\ori) \simeq L_E(\QCoh(\Spec \Sph)) \simeq L_E(\Sp) = \Sp_{\vec{h}};\end{equation}
	this has already been observed by Mathew \cite[\textsection10.2]{akhilgalois} and Gregoric \cite[Th.2]{rokfiltration} in the $p$-local case. 
    
    The identification (\ref{eq:introdescentableetheorey}) also points towards the proof of \Cref{maintheorem:reconstruction2}. Indeed, the equivalence of $\infty$-categories of \Cref{maintheorem:reconstruction2} is the composition
    \[\Stk^{\aff}_{/\M_\FGvech^\ori} \xrightarrow{\Phi} \CAlg(\QCoh(\M_\FGvech^\ori))^\op \xrightarrow{\simeq} \CAlg(\Sp_{\vec{h}}),\]
    the first functor being the general equivalence of \Cref{relativespectrumtheorem} sending affine maps of stacks $f\colon \Y \to \X$ to $f_\ast \O_\Y$, and the second being induced by (\ref{eq:introdescentableetheorey}).
	
	\subsection{Outline}
	The goal of this paper is to study two broad concepts, universal $0$-affineness and locally descendability of morphisms of stacks, at three different levels of generality. We start with maximal generality, working with stacks on a general Grothendieck site (\Cref{sec:sheaves}); next, we specialize this to nonconnective spectral algebraic geometry (\Cref{sec:nonconsag}), a variant on Lurie's set-up in \cite{sag}; and finally, we specialize further to stacks over $\M_\FG^\ori$, the titular \emph{complex-periodic geometry} (\Cref{sec:complexperiodicgeometry}). In more detail:
	
	\paragraph{(\Cref{sec:sheaves})}
	In \Cref{ssec:stackssection}, we define an $\infty$-category of \emph{stacks} $\Stk_\tau(\csite)$ associated with a site $(\csite,\tau)$. We do not assume $\csite$ to be small, however, if $\csite$ is small, then $\Stk_\tau(\csite)$ is equivalent to the category of sheaves on $(\csite,\tau)$. This applies to $\CAlg = \CAlg(\Sp)$ equipped with the $\fpqc$ topology, which has famously caused many headaches, even classically; see \cite{waterhousefpqc}. We then show that $\Stk_\tau(\csite)$ avoids these technical difficulties by demonstrating some of its topos-like behaviour and the functoriality of this construction. We will use this general notion of stacks in nonconnective spectral algebraic geometry in the next section, but this notion also fits into contexts such as derived algebraic geometry based on other variants of $\E_\infty$ rings. We also axiomatise the notion of a \emph{quasi-coherent sheaf context} on $\Stk_\tau(\csite)$, which is the minimal data needed to set up a workable theory of quasi-coherent sheaves in this generality (\Cref{sssec:qcohcontexts}).
	
	With these foundations, we then study the general notions of \emph{affine}, \emph{$0$-affine}, \emph{$0$-semiaffine}, and \emph{universally $0$-affine} morphisms in \Cref{ssec:genericaffineness}. The arguments in this section mainly depend on the algebra of $\CAlg(\PrLst)$, and so the work is largely categorical. Perhaps most importantly, we show in \Cref{prop:univ0affinefacts} that universally $0$-affine morphisms, an affine-local condition, are actually $0$-affine, a global condition. We also prove a version of Lurie's \emph{weak Tannaka duality} in \Cref{sec:appendix} in this general context. Finally, in \Cref{ssec:descentstacks} we define the descent stack $\D_f$ associated to a map of stacks, and explain its relationship with local descendability. In particular, we prove \Cref{thm:locallydescendable}, which allows us to deduce that $j \colon \D_f \to \X$ is locally descendable if the original map $f\colon \Y \to \X$ was.
	
	\paragraph{(\Cref{sec:nonconsag})}
	The goal of this section is to apply the foundations of the previous section to spectral algebraic geometry. In \Cref{sssec:basicdefinitions}, we define $\Stk = \Stk_\fpqc(\CAlg^\op)$ as the $\infty$-category of stacks associated with the fpqc topology on $\E_\infty$ rings. We also make similar definitions for connective stacks $\Stk^\cn$ and classical stacks $\Stk^\heartsuit$. The functoriality of our generalized stack construction immediately induces colimit preserving connective cover and truncation functors
	\[\begin{tikzcd}
		{\Stk^\heartsuit} & {\Stk^\cn} & {\Stk,}
		\arrow[shift left=2, from=1-2, to=1-1, "{\tau_{\leq 0}}"]
		\arrow[shift left=2, from=1-1, to=1-2]
		\arrow[shift right=2, from=1-2, to=1-3]
		\arrow[shift right=2, from=1-3, to=1-2, "{\tau_{\geq 0}}", swap]
	\end{tikzcd}\]
    with left adjoints written on top, whose adjoints in both cases are fully-faithful; see \Cref{th:coversandtruncations}. We equip $\Stk$ with the canonical quasi-coherent sheaf context determined by $\QCoh(\Spec\A) = \Mod_\A$ for an $\E_\infty$ ring $\A$ (\Cref{def:thecanonicalqcohcontext}). In \Cref{ssec:spdm}, we connect our definition of $\Stk$ with Lurie's notion of a \emph{nonconnective spectral Deligne--Mumford stack}. In particular, we show that there is a fully faithful embedding $\SpDMnc \to \Stk_\et(\CAlg^\op)$ into \emph{\'{e}tale} stacks (\Cref{thm:dmstacksintoetalestacks}), and that the composite with fpqc sheafification induces a functor $\SpDMnc \to \Stk$ which is fully faithful on those objects satisfying flat descent (\Cref{prop:dmstacksintoflatstacks}). We then show that objects of $\SpDMnc$ with quasi-affine diagonal satisfy flat descent (\Cref{prop:qaffdiagonalmeansdescent}). Independent to this discussion of flat descent, we show that our notion of quasi-coherent sheaves on $\Stk$ agrees with that of Lurie (\Cref{ex:qcohofspDM}) for \emph{all} nonconnective spectral Deligne--Mumford stacks. In particular, $\Stk$ is a topos-like cocomplete extension of $\SpDMnc$ built to study questions regarding quasi-coherent sheaves.
	
	In \Cref{ssec:quasiaffine}, we define quasi-affine morphisms of stacks and show that such morphisms are universally $0$-affine (\Cref{ex:quasiaffineareU0A}). In particular, this shows how \Cref{cor:toboundedaffinenessintro} is a generalization of \Cref{maintheorem:boundedaffineness}. In \Cref{ssec:flatnessandgeometricstacks}, we study stacks that admit a flat atlas, the so-called \emph{geometric stacks}, and study flat maps between such stacks. These concepts allow one to pass easily between $\Stk$ and either $\Stk^\cn$ or $\Stk^\heartsuit$. For example, both the connective cover and underlying stack constructions preserve pullbacks of geometric stacks along flat maps (\Cref{flatnessislovely}). Another useful result roughly states that a flat map between geometric stacks is (quasi-) affine  if and only if the underlying map of classical stacks is (quasi-) affine (\Cref{cor:quasiaffinenessfromtheheart}).
	
	\paragraph{(\Cref{sec:complexperiodicgeometry})}
	Now we inject some chromatic homotopy theory into our spectral algebraic geometry. In \Cref{ssec:MFGOR}, we identify complex-periodic stacks as precisely those stacks living over $\M_\FG^\ori$, the \emph{moduli stack of oriented formal groups}, and show that mapping into $\M_\FG^\ori$ is a \emph{property} rather than \emph{structure}. This allows us to compare $\M_\FG^\ori$ with $\D_{\MUP}$, the descent stack associated to $\MUP$, any $\E_\infty$ form of periodic complex cobordism, immediately leading to many fundamental properties of $\M_\FG^\ori$. In \Cref{ssec:classicalmoduliofformalgroups}, we discuss the connection between stacks living over $\M_\FG^\ori$ and their underlying classical stacks living over the classical moduli stack $\M_\FG^\heartsuit$ of formal groups. For example, we show that a complex-periodic stack is Landweber exact if and only if the underlying classical stack is Landweber exact (\Cref{eqdefofLEFT}). This allows us to view Mathew--Meier's approach to chromatic homotopy theory using spectral algebraic geometry as a special case of complex-periodic stacks (\Cref{directcomparisontoMATHEWMEIER}), as well as to produce many examples of complex-periodic stacks which are flat or (quasi-) affine over $\M_\FG^\ori$ (\Cref{ssec:examplesofLEFTandQAFF}). In \Cref{ssec:complextheoriesofboundedheight}, we introduce substacks of $\M_\FG^\ori$ parametrizing the \emph{height} of the underlying formal group, written as \smash{$\M_{\smash{\FGvech}}^\ori$}. We then show that \smash{$\M_{\smash{\FGvech}}^\ori$} is also the descent stack $\sfD_E$ for a particular choice of Landweber exact $\E_\infty$ ring $E$, and as a result of \Cref{ssec:descentstacks} and a refinement of the smash product theorem, we compute $\QCoh(\M_\FGvech^\ori)$ (\Cref{thm:qcohofDvech}). Finally, in \Cref{ssec:affinenesssection}, we prove \Cref{maintheorem:boundedaffineness} (in its improved form as in \Cref{cor:toboundedaffinenessintro}). As advertised, having established the correct formalism, the proof amounts to two sentences.
	
	In \Cref{sssec:maori}, we study the \emph{complex-periodification} of an $\E_\infty$-ring $\A$, written as $\M_\A^\ori$, a spectral refinement of Hopkins' stack associated to an $\E_\infty$ ring. We prove basic properties of these stacks $\M_\A^\ori$ (\Cref{sssec:complexperiodicifcofeinftyrings}), give the same salient examples (\Cref{sssec:maoriexamples}), and then show that many of these stacks are $0$-affine (\Cref{sssec:0affinenessforMAORI}), proving \Cref{maintheorem:moaoriaffineness}. In the last section \Cref{ssec:reconstructionsection}, we prove both \Cref{maintheorem:reconstruction,maintheorem:reconstruction2}, our two main results concerning \emph{reconstruction}. We begin with the definition of reconstructible stacks and their basic properties (\Cref{sssec:basicreconstructionproeprtioes}), and follow that by proving a large collection of examples of such stacks (\Cref{sssec:exampleofreconstruction}) including all of those from (\ref{eq:examplesofreconstructiblestackies}).
	
	\subsection*{Acknowledgements}
    The authors would like to thank Emma Brink, Christian Carrick, Rok Gregoric, Lennart Meier, Marius Verner Nielsen, Sven van Nigtevecht, and Lucas Piessevaux, for enlightening conversations surrounding these topics. 
	
	The second author was supported as an associate member of the Hausdorff Center for Mathematics at the University of Bonn (\texttt{DFG GZ 2047/1}, project ID \texttt{390685813}), the Max-Planck Institute for Mathematics, and the DFG-funded research 
training group GRK 2240: Algebro-Geometric Methods in Algebra, 
Arithmetic and Topology. The third author is an associate member of the SFB: Higher invariants. The second and third author would like to thank the Isaac Newton Institute for Mathematical Sciences, Cambridge, for support and hospitality during the programme \emph{Equivariant homotopy theory in context} where work on this paper was undertaken. This work was supported by EPSRC grant no EP/K032208/1.

\section{General derived geometries}\label{sec:sheaves}
The generic behaviour of $0$-affine and locally descendable morphisms of stacks is largely captured by the algebra of $\CAlg(\PrLst)$ and general sheaf contexts. For this reason, we begin by axiomatising our general set-up for stacks on a site $(\csite, \tau)$, defining universally $0$-affine morphisms and locally descendable morphisms between such stacks, and proving their basic properties.

\subsection{Generalized stacks and quasi-coherent sheaves}\label{ssec:stackssection}
Our main example of a stack is $\M_\FG^\ori$, which can, at least informally, be written as the colimit of $\Spec \MUP^{\otimes(\bullet+1)}$ in the category of fpqc sheaves on $\CAlg$. In general, if a category of sheaves admits a sheafification functor, then colimits of sheaves are computed by sheafifying the associated colimit in presheaves. Unfortunately, fpqc sheaves on $\CAlg$ do \textbf{not} generally admit sheafifications for size reasons, hence such colimits do not generally exist. This problem already occurs classically for fpqc sheaves on discrete rings; see \cite{waterhousefpqc}. In this section, we define a category of stacks $\Stk(\csite, \tau)$ for a possibly large site $(\csite,\tau)$ that sidesteps this issue, and prove some of its basic properties. We also axiomatize the behaviour of quasi-coherent sheaves to \emph{quasi-coherent sheaf contexts}, see \Cref{def:qcohcontext}. This generalized set-up is then used in all further subsections to discuss generic affineness phenomen{\ae}.

\subsubsection{Stacks on large categories}
We begin by fixing some terminology and recalling some basic properties of sheaves on sites in the setting of $\infty$-category theory. A convenient reference for this material is \cite[\textsection A.1]{syntheticspectra}.

\begin{mydef}\label{def:coaccessiblesite}
A \emph{site} is a pair $(\csite,\tau)$ consisting of a category $\csite$ which admits all finite limits, together with a Grothendieck pretopology $\tau$ which we further assume to be subcanonical.\footnote{Recall that a Grothendieck topology $\tau$ on $\csite$ is \emph{subcanonical} if the Yoneda embedding $\csite \to \Fun(\csite^\op, \Spc)$ factors through $\tau$-sheaves.} Given a small site $(\csite,\tau)$, meaning that $\csite$ is a small category, write
\[
\Shv_\tau(\csite)\subseteq\presheaves(\csite)
\]
for the category of $\tau$-sheaves on $\csite$ and
\[
L^\tau\colon \presheaves(\csite) \to \Shv_\tau(\csite)
\]
for the functor of sheafification. This is by definition left adjoint to the inclusion, which exists by \cite[Pr.6.2.2.7]{htt}.
\end{mydef}

\begin{mydef}
A \emph{morphism} $f\colon (\csite,\tau) \to (\csite',\tau')$ of sites is a functor $f\colon \csite \to \csite'$ which preserves cover families and preserves pullbacks of coverings. An \emph{exact morphism} $f\colon (\csite,\tau)\to (\csite',\tau')$ of sites is a morphism that preserves all finite limits.
\end{mydef}

\begin{prop}\label{prop:morphismsmallsites}
Let $f\colon (\csite,\tau) \to (\csite',\tau')$ be a morphism of small sites. Then restriction
\[
f^\ast\colon \presheaves(\csite') \to \presheaves(\csite)
\]
preserves sheaves, and therefore its restriction to sheaves is the right adjoint in an adjunction
\[
f_! : \Shv_\tau(\csite) \rightleftarrows \Shv_{\tau'}(\csite') : f^\ast,
\]
where $f_!$ is isomorphic to the composite $L_{\tau'}\mathrm{Lan}_f$ of the unique cocontinuous extension of $f$ followed by $\tau'$-sheafification. Moreover, if $f$ is an exact morphism, then $f_!$ preserves finite limits.
\end{prop}

\begin{proof}
All but the last claim are proved in \cite[Pr.A.1.3]{syntheticspectra}. It therefore remains to verify that if $f$ is exact, then $f_!$ preserves finite limits. As $f_!$ factors as the composite
\[
\Shv_\tau(\csite) \subseteq \presheaves(\csite) \xrightarrow{f_\ast}\presheaves(\csite') \xrightarrow{L^{\tau'}} \Shv_{\tau'}(\csite'),
\]
it suffices to verify just that $f_\ast$ preserves finite limits, which is \cite[Pr.6.1.5.2]{htt}.
\end{proof}

\begin{prop}\label{prop:univpropertysheaves}
Let $(\csite,\tau)$ be a small site. Then for any category $\calD$ with small limits, restriction along the Yoneda embedding $\csite\to\Shv_\tau(\csite)$ induces a fully faithful embedding
\[
\Fun^{\mathrm{R}}(\Shv_\tau(\csite)^\op,\calD)\to\Fun(\csite^\op,\calD),
\]
with essential image spanned by the $\tau$-sheaves.
\end{prop}
\begin{proof}
This is \cite[Pr.1.3.1.7]{sag}.
\end{proof}

Care must be taken when considering the notion of sheaves on a large site. We take the perspective that \cref{prop:univpropertysheaves} is the defining property of $\Shv_\tau(\csite)$, and give a definition that extends this to the situation where $\csite$ is not small.

The rest of this subsection works with a locally small site $(\csite,\tau)$, meaning the category $\csite$ is locally small. Given any small subcategory $\csite_0\subseteq\csite$, we may form the closure of $\csite_0$ under finite limits $(\csite_0',\tau)\subseteq(\csite,\tau)$ containing $\csite_0$. 
In this way, we see that the class of small subsites of $(\csite,\tau)$ forms a cofinal subsystem of the class of small subcategories of $\csite$.

\begin{mydef}\label{def:stacks}
Let $(\csite,\tau)$ be a locally small site. The category of \emph{stacks} on $(\csite,\tau)$ is defined as the colimit
\[
\Stk_\tau(\csite)\simeq\colim \Shv_\tau(\csite_0)
\]
indexed over all small subsites $(\csite_0,\tau) \subseteq (\csite,\tau)$.
\end{mydef}

If $(\csite,\tau)$ is a small site, then $\Stk_\tau(\csite)\simeq\Shv_\tau(\csite)$, as the above diagram has a terminal object. There is another key recurring example.

\begin{mydef}\label{def:small_presheaf}
Let $\csite$ be a category. A presheaf $X\in \Fun(\csite^{\op},\Spc)$ is \emph{small} if it is a small colimit of representables. We denote the category of small presheaves by $\presheaves(\csite)$.
\end{mydef}

\begin{example}\label{ex:trivialtopology}
If $\tau$ generates the trivial Grothendieck topology on $\csite$, where the only covers are isomorphisms, then by definition we have the identification
\[
\Stk_\tau(\csite) \simeq \presheaves(\csite).
\]
\end{example}

As we assumed that $\tau$ was subcanonical, for each inclusion of small subsites $i\colon \csite_0 \subseteq \csite_1$, the induced diagram of Yoneda embeddings commutes:
\[\begin{tikzcd}
    {\csite_0}\ar[r]\ar[d, "i"] &   {\Shv_\tau(\csite_0)}\ar[d, "{i_!}"]   \\
    {\csite_1}\ar[r] &   {\Shv_\tau(\csite_1)}
\end{tikzcd}\]
In particular, in this case, there is a well-defined Yoneda embedding
\begin{equation}\label{eq:yonedaforstacks}
\Spec \colon \csite \to \Stk_\tau(\C)    .
\end{equation}

\begin{mydef}\label{def:affinestacks}
    Let $(\csite, \tau)$ be a locally small site. An \emph{affine stack} is a stack in the essential image of the Yoneda embedding of (\ref{eq:yonedaforstacks}).
\end{mydef}

\subsubsection{Basic properties of stacks}
We now establish some properties of the category $\Stk_\tau(\csite)$ of stacks on a potentially large (locally small) site. First, we explain why $\Stk_\tau(\csite)$ has colimits.

\begin{prop}
For a locally small site $(\csite, \tau)$, the category $\Stk_{\tau}(\calC)$ admits small colimits and finite limits.
\end{prop}

\begin{proof}
Let $\kappa > \lambda$ be two regular cardinals, and recall that $\kappa$-filtered colimits of categories with $\lambda$-small (co)limits along $\lambda$-(co)continuous functors again admit $\lambda$-small (co)limits; the same proof as \cite[Lm.7.3.5.10]{ha}, which considers the special case of colimits when $\kappa = \omega$, applies. Since $\Stk_{\tau}(\calC)$ is a colimit over a diagram of functors which preserve all small colimits and finite limits, and which is filtered with respect to all $\kappa$ in our universe, the claim follows.
\end{proof}

Next, we see how functors out of $\csite$ satisfying $\tau$-descent extend uniquely over $\Stk_\tau(\csite)$, giving our desired universal property for $\Stk_\tau(\csite)$.

\begin{prop}\label{prop:limitextension}
For any category $\calD$ with small limits, restriction along $\Spec\colon \csite\to\Stk_\tau(\csite)$ induces a fully faithful embedding
\[
\Fun^{\mathrm{R}}(\Stk_\tau(\csite)^\op,\calD) \to \Fun(\csite^\op,\calD),
\]
with essential image spanned by the $\tau$-sheaves $\csite^\op\to\calD$.
\end{prop}

\begin{proof}
By construction, we have
\[
\Fun^{\mathrm{R}}(\Stk_\tau(\csite)^\op,\calD)\simeq\lim\Fun^{\mathrm{R}}(\Shv_\tau(\csite_0)^\op,\calD)\qquad \text{and} \qquad \Fun(\csite^\op,\calD)\simeq\lim\Fun(\csite_0^\op,\calD).
\]
As a functor $f\colon \csite\to\calD$ is a $\tau$-sheaf if and only if its restriction to any small subsite $\csite_0\subseteq\csite$ is a $\tau$-sheaf, the proposition follows from \cref{prop:univpropertysheaves}.
\end{proof}

\begin{remark}
The proposition above shows that our definition of stacks always admits the expected universal property, and so is the correct definition in the generality in which we operate. As an alternative to $\Stk_\tau(\csite)$, one can define the category of \emph{small sheaves} on $(\csite,\tau)$ to be the full subcategory
\[
\Shv_\tau(\csite)\subseteq\presheaves(\csite)
\]
of small presheaves on $\csite$ which have the property of being $\tau$-sheaves. By \cref{prop:pseudosheafify} below, there is a natural functor
\[
\Shv_\tau(\csite) \to \Stk_\tau(\csite).
\]
We expect that this is an equivalence for many of the large sites that arise in practice. This is the case, for example, if $(\csite,\tau)$ is an \emph{accessible site} in the sense of \cite[\textsection1.4]{pyknotic}.\footnote{This is revisited in work-in-progress of Emma Brink, where the notions of accessible sheaves and sheaves on large categories are explored in more depth.} 
Outside of this generality, we consider $\Stk_{\tau}(\csite)$ to be better behaved. Our main example of interest is the fpqc topology on affine schemes (based on $\E_\infty$ rings) and unfortunately, we have been unable to verify that this site is accessible.
\end{remark}

The above proposition is incredibly useful to formally extend concepts for $\E_\infty$ rings which satisfy flat descent to concepts defined on stacks, such as notions of elliptic curves, formal groups, $\P$-divisible groups, as well as their (pre)oriented variants; this will be heavily used in \cite{temperedglobal}.

Moreover, the above universal property yields the following useful corollary.

\begin{cor}\label{cor:bigsheaves}
Restriction along the functor $\Spec$ defines a fully faithful embedding
\[
\Stk_\tau(\csite) \to \Fun(\csite^\op,\widehat{\spaces}),
\]
with essential image contained in the category $\widehat{\Shv}_\tau(\csite)$ of large $\tau$-sheaves on $\csite$ with values in large spaces.
\end{cor}
\begin{proof}
Both functors
\[
\Stk_\tau(\csite) \to \Fun^{\mathrm{R}}(\Stk_\tau(\csite)^\op,\widehat{\spaces})\to\Fun(\csite^\op,\widehat{\spaces})
\]
are fully faithful by the Yoneda lemma and \Cref{prop:limitextension}.
\end{proof}

Next, we show that, in a precise sense, $\Stk_\tau(\csite)$ can be thought of as a category of large sheaves on $\csite$ which are determined by a small amount of data. 

\begin{cor}\label{cor:accessiblesites}
Suppose that $\csite^\op$ is accessible as a category, and suppose that the full subcategory $\csite_\kappa\subseteq\csite$ of $\kappa$-cocompact objects forms a subsite for all sufficiently large cardinals $\kappa$. Then the natural map
\[
\colim_\kappa\Shv_\tau(\csite_\kappa) \xrightarrow{\simeq} \Stk_\tau(\csite)
\]
is an equivalence, where the colimit is indexed over sufficiently large regular cardinals $\kappa$.
\end{cor}
\begin{proof}
If $\tau$ generates the trivial topology, then this follows from \Cref{ex:trivialtopology} and \cite[Pr.A.2(3)]{diracii}. From here, the same proof as \cref{prop:limitextension} applies equally well to show that $\colim_\kappa\Shv_\tau(\csite_\kappa)$ has the same universal property as $\Stk_\tau(\csite)$, so the two categories must be equivalent.
\end{proof}

We next prove some functoriality properties of the construction of $\Stk_\tau(\csite)$.

\begin{prop}\label{prop:accessiblemorphismofsites}
Any morphism $f\colon (\csite,\tau)\to (\csite',\tau')$ of sites induces a small colimit-preserving functor
\[
f_!\colon \Stk_{\tau}(\csite) \to \Stk_{\tau'}(\csite').
\]
If $f$ is exact, then $f_!$ preserves finite limits.
\end{prop}
\begin{proof}
For each inclusion of full subsites $\csite_0\subseteq\csite_1\subseteq\csite$, one may find a corresponding inclusion of full subsites $\csite'_0\subseteq\csite'_1\subseteq \csite'$ for which $f$ induces a commutative diagram
\begin{center}\begin{tikzcd}
(\csite_0,\tau)\ar[r]\ar[d,"f"]&(\csite_1,\tau)\ar[d,"f"]\\
(\csite'_0,\tau')\ar[r]&(\csite'_1,\tau')
\end{tikzcd}\end{center}
of morphisms of sites. This, in turn, induces a commutative diagram
\begin{center}\begin{tikzcd}
\Shv_\tau(\csite_0)\ar[r]\ar[d, "{f_!}"] & \Shv_\tau(\csite_1)\ar[d, "{f_!}"]\\
\Shv_{\tau'}(\csite_0')\ar[r]&\Shv_{\tau'}(\csite_1')
\end{tikzcd}\end{center}
of left adjoints, left exact if $f$ is exact. Passing to colimits yields the functor
\[
\Stk_{\tau}(\csite) \to \Stk_{\tau'}(\csite').
\]
We claim that this functor preserves small colimits and preserves finite limits if $f$ is left exact. Since these are computed in $\Shv_{\tau}(\csite_0)$ for some small $\csite_0\subset \csite$, this follows from \cref{prop:morphismsmallsites}.
\end{proof}

As a corollary, we deduce that small presheaves admit \emph{stackifications}.

\begin{cor}\label{prop:pseudosheafify}
The functors of sheafification $L^\tau\colon \presheaves(\csite_0)\to\Shv_\tau(\csite_0)$ assemble into a small colimit and finite limit preserving functor
\[
L^\tau\colon \presheaves(\csite) \to \Stk_\tau(\csite).
\]
\end{cor}
\begin{proof}
Apply the previous proposition to the identity of $\csite$, viewed as a map of sites
\[
\id\colon (\csite,\mathrm{isos})\to (\csite,\tau),
\]
from $\csite$ with the trivial topology, in which the only covers are the isomorphisms, to $\csite$ with the topology $\tau$; see \Cref{ex:trivialtopology}.
\end{proof}

\begin{prop}\label{cor:adjunctionsandmapsofsites}
Let $f\colon (\csite, \tau) \to (\csite', \tau')$ and $g\colon (\csite', \tau') \to (\csite, \tau)$ be morphisms of sites where $f$ is left adjoint to $g$. Then $f_!\colon \Stk_\tau(\csite) \to \Stk_{\tau'}(\csite')$ is left adjoint to $g_!\colon \Stk_{\tau'}(\csite') \to \Stk_{\tau}(\csite)$. Moreover, if $f$ (resp.\ $g$) is fully faithful, then $f_!$ (resp.\ $g_!$) is fully faithful.
\end{prop}

\begin{proof}
    We can immediately reduce the case where $\csite$ and $\csite'$ are small. By the 2-functoriality of left Kan extension, we obtain an adjunction 
    \[
    \mathrm{Lan}_f\colon \presheaves(\csite)\rightleftarrows \presheaves(\csite') \noloc \mathrm{Lan}_g.
    \]
    Now we claim that this adjunction restricts to an adjunction on categories of sheaves. Recall that $f_! \simeq L^{\tau} \mathrm{Lan}_f$, and similarly for $g_!$. Therefore, we may define the natural transformation
    \[
    L^{\tau'} \mathrm{Lan}_f L^{\tau} \mathrm{Lan}_g \simeq L^{\tau'} \mathrm{Lan}_f \mathrm{Lan}_g \xrightarrow{L^{\tau'} (\epsilon)} L^{\tau'}.
    \]
    It is simple to show that this defines the necessary counit when restricted to $\Shv_{\tau'}(\csite')$.
\end{proof}

The following is a summary of the useful categorical properties of $\Stk_\tau(\csite)$, reflecting how it behaves similarly to a topos; compare Giraud's axioms from \cite[Th.6.1.0.6(3)]{htt}.

\begin{prop}\label{pr:stackisbasicallyatopos}
Let $(\csite,\tau)$ be a locally small site.
\begin{enumerate}
\item Coproducts in $\Stk_{\tau}(\csite)$ are disjoint.
\item Colimits in $\Stk_{\tau}(\csite)$ are universal.
\item Groupoids in $\Stk_{\tau}(\csite)$ are effective.
\item The functor $\Stk_\tau(\csite)^\op\to\widehat{\Cat}$ sending $X$ to $\Stk_\tau(\csite)_{/X}$ preserves small limits. 
\end{enumerate}
\end{prop}
\begin{proof}
(1,2,3)~~Each of these statements concerns finite limits and colimits in $\Stk_\tau(\csite)$, and their interaction. Since these are inherited from $\presheaves(\csite)$, the statements follow.

(4)~~Note that taking slices is a pullback, which commutes with filtered colimits of categories. Therefore, if $X = \colim_i X_i$ is a small colimit in $\Stk_{\tau}(\csite)$, we may factor the canonical map to the limit as the following composite of equivalences:
\[\begin{tikzcd}
	{\Stk_{\tau}(\csite)_{/X} } && {\lim \Stk_{\tau}(\csite)_{/X_i}} \\
	{\colim \Shv_\tau(\csite_0)_{/X} } & {\colim\lim\Shv_{\tau}(\csite_0)_{/X_i}} & {\lim\colim\Shv_{\tau}(\csite_0)_{/X_i}}
	\arrow[from=1-1, to=1-3]
	\arrow["\sim"', from=1-1, to=2-1]
	\arrow["\sim", from=2-1, to=2-2]
	\arrow["\sim", from=2-2, to=2-3]
	\arrow["\sim", from=2-3, to=1-3]
\end{tikzcd}\]
\end{proof}

Finally, let us also consider the relationship between categories of stacks and Grothendieck constructions/unstraightenings. For the statement, we require the following notation: Given a Cartesian fibration $p\colon \int_\calC F\to \calC$ we say a diagram $I\to \int_\calC F$ is \emph{Cartesian} if it factors through the subcategory spanned by the ($p$-)Cartesian edges. 

\begin{prop}\label{prop:stacksonunstraightenings}
    Let $f\colon \C^\op \to \Cat$ be a $\tau$-sheaf, with unique limit preserving extension $F\colon \Stk_\tau(\csite)^\op \to \Cat$. Then:
    \begin{enumerate}
        \item The collection of Cartesian lifts of covers in $\C$ define a Grothendieck topology on $\int_\C f$, the Cartesian unstraightening of $f\colon\C^{\op}\to \Cat$.
        \item $\int_{\Stk_\tau(\csite)} F$ admits colimits of Cartesian diagrams.
        \item The inclusion $\int_\csite f \to \int_{\Stk_\tau(\csite)} F$ induces a fully faithful embedding 
        \[
        \Fun^{\textup{R}_{\mathrm{ct}}}(\Bigl(\int_{\Stk_\tau(\csite)} F\Bigr)^{\op},\calD) \to \Fun(\Bigl(\int_\csite f\Bigr)^{\op} ,\calD)
        \]
        from functors which preserve colimits of Cartesian diagrams to functors on $\int_{\calC} f$, with essential image given by the $\tau$-sheaves $\int_{\calC}f\to \Cat$.
    \end{enumerate}
\end{prop}

\begin{proof}
    Part 1 follows as $\int_\C f$ admits pullback of Cartesian morphisms, preserved by the projection to $\csite$. For part 2, consider a Cartesian diagram 
    \[G\colon I\to \int_{\Stk_\tau(\csite)} F, \quad i\mapsto (X_i,G_i).
    \]
    Let us write $X$ for the colimit of the composite of $G$ with the projection to $\csite$. Recalling that Cartesian sections of the unstraightening is one model of the limit of categories, $G$ is equivalent to an object of $\lim F(X_i) \simeq F(X)$ and so canonically an object of $\int_{\Stk(\csite)}F$ living in the fiber over $X$. We will show that $G$, regarded as an object of $\int_{\Stk(\csite)}F$, is the colimit of $G$, regarded as a diagram $I\to \int_{\Stk_\tau(\csite)} F$. To do this, we consider the following diagram
    \[\begin{tikzcd}
	{\Hom_{F(X)}(G,f^*G')} & {\lim\Hom_{F(X_i)}(G_i,f_i^*G'_i)} \\
	{\Hom_{\int F}((X,G),(Y,G'))} & {\lim\Hom_{\int F}((X_i,G_i),(Y,G'))} \\
	{\Hom_{\Stk_\tau(\csite)}(X,Y)} & {\lim\Hom_{\Stk_\tau(\csite)}(X_i,Y)} \\
	{\{f\}} & {\{(f_i)\}}
	\arrow["\sim", from=1-1, to=1-2]
	\arrow[from=1-1, to=2-1]
	\arrow[from=1-2, to=2-2]
	\arrow[from=2-1, to=2-2]
	\arrow[from=2-1, to=3-1]
	\arrow[from=2-2, to=3-2]
	\arrow[""{name=0, anchor=center, inner sep=0}, "\sim", from=3-1, to=3-2]
	\arrow[hook, from=4-1, to=3-1]
	\arrow[equals, from=4-1, to=4-2]
	\arrow[hook, from=4-2, to=3-2]
\end{tikzcd}\]
By definition, the bottom map is an equivalence, since $X = \colim X_i$. On the other hand, we may consider the fiber of this map over any point $\{f\}$, which is drawn at the top of this diagram. This is an equivalence since $F(X) = \lim F(X_i)$. Therefore, the middle map is an equivalence, thus proving the statement.

For part 3, we can immediately reduce to the case where $\csite$ is small. Now we first consider the case where the topology is trivial. In this case, we need to prove that $\int_{\presheaves(\csite)} F$ is freely generated by $\int_\csite f$ under colimits of Cartesian diagrams. For this, it suffices to prove that $\int_{\Stk_{\tau}(\csite)} F$ is generated by $\int_{\csite}f$ under colimits of Cartesian diagrams and that mapping out of an object in $\int_{\csite}f$ preserves colimits of Cartesian diagrams. The first is immediately from the description of colimits of Cartesian diagrams above. For the second, consider $(C,H)\in \int_\csite f$. Since $C$ is completely compact in $\presheaves(\csite)$, we obtain an equivalence 
\[
\Hom_{\presheaves(\csite)}(C,X) \simeq \colim \Hom_{\presheaves(\csite)}(C,X_i),
\]
and so it suffices to show that the natural map
\[
\colim \Hom((C,H),(X_i,G_i)) \to \Hom((C,H),(X,G))
\]
is an equivalence after pulling back along the maps $(g_i)_*\colon\Hom_{\presheaves(\csite)}(C,X_i)\to \Hom_{\presheaves(\csite)}(C,X)$ for all $i\in I$. However, this follows from the fact that
\[\begin{tikzcd}
	{\Hom_{\int F}((C,H),(X_i,g_i^*G'))} & {\Hom_{\int F}((C,H),(X,G'))} \\
	{\Hom_{\presheaves(\csite)}(C,X_i)} & {\Hom_{\presheaves(\csite)}(C,X)}
	\arrow[from=1-1, to=1-2]
	\arrow[from=1-1, to=2-1]
	\arrow[from=1-2, to=2-2]
	\arrow[""{name=0, anchor=center, inner sep=0}, "{(g_i)_*}", from=2-1, to=2-2]
	\arrow["\lrcorner"{anchor=center, pos=0.125}, draw=none, from=1-1, to=0]
\end{tikzcd}\]
is a pullback, where the top map is given by postcomposing by a Cartesian lift of $g_i$, see \cite[Pr.2.4.4.3]{htt}. Finally, we consider the case of a non-trivial topology $\tau$ on $\csite$. To begin with, we consider the inclusion
\[
\int_{\Shv_{\tau}(\csite)} F \hookrightarrow \int_{\presheaves(\csite)} F.
\]
We claim that this admits a left adjoint, which is then a Bousfield localization. So consider an object $(X,G)$ in $\int_{\presheaves(\csite)} F$. Since $f$ is a $\tau$-sheaf, we obtain an isomorphism $G(L^\tau X)\xrightarrow{\sim} G(X)$ and so an object $(L^\tau X,G)$ in $\int_{\Shv_{\tau}(\csite)} F$ under $(X,G)$. A simple calculation shows that this is the localization of $(X,G)$ into $\int_{\Shv_{\tau}(\csite)} F$. Just as in \cite[Pr.5.5.4.20]{htt}, as a consequence, we find that for every category $\calD$ with small limits, we obtain a fully faithful inclusion 
\[
\Fun^{\mathrm{R}_{\mathrm{ct}}}(\int_{\Shv_{\tau}(\csite)} F,\calD)\hookrightarrow \Fun^{\mathrm{R}_{\mathrm{ct}}}(\int_{\presheaves(\csite)} F,\calD) \simeq \Fun(\int_{\csite} f,\calD),
\]
which another simple calculation shows has essential image precisely the $\tau$-sheaves.
\end{proof}

\subsubsection{Quasi-coherent sheaf contexts and global sections}\label{sssec:qcohcontexts}
Let $\PrLst$ denote the category of stable presentable categories and left adjoints.

\begin{mydef}\label{def:affinefunctor}
A map $f^\ast \colon \ccat\to\dcat$ in $\CAlg(\PrLst)$ is \emph{affine} if its right adjoint $f_\ast\colon \dcat \to \ccat$ induces a (necessarily symmetric monoidal) equivalence
\[
\dcat \xrightarrow{\simeq} \Mod_{f_\ast\dmonoidalunit}(\ccat).
\]
\end{mydef}

This categorical form of affineness has been studied by several people, including \cite{akhilgalois,nilpotenceanddescent, behrensshah,barthelcarmelischlankyanovski2024chromatic}. We will summarize some of its main properties below in \cref{sssec:0affine}.

\begin{mydef}\label{def:qcohcontext}
A \emph{quasi-coherent context} $(\csite,\tau,\sfQ)$ consists of a locally small site $(\C,\tau)$ together with a $\tau$-sheaf $\sfQ\colon \C^\op\to\CAlg(\PrLst)$ of stable presentably symmetric monoidal categories satisfying the following conditions:
\begin{enumerate}
\item $\sfQ$ sends morphisms in $\csite$ to affine morphisms in $\CAlg(\PrLst)$;
\item $\sfQ$ sends Cartesian squares in $\C$ to coCartesian squares in $\CAlg(\PrLst)$.
\end{enumerate}
Given a quasi-coherent context $(\csite,\tau,\sfQ)$, we write
\[
\QCoh\colon \Stk_\tau(\csite)^\op\to\CAlg(\PrLst)
\]
for the unique limit-preserving extension of $\sfQ$ guaranteed by \cref{prop:limitextension}.
\end{mydef}

\begin{example}\label{ex:basicqcohcontext}
Let $\cstablecat \in \CAlg(\PrLst)$. Then the functor
\begin{equation}\label{eq:modulecategoryfunctor}
\CAlg(\cstablecat) \to \CAlg(\Pr^L_\st),\qquad A \mapsto \Mod_A(\cstablecat)
\end{equation}
satisfies the (covariant analogues of the) conditions of \cref{def:qcohcontext}; see \cite[Th.4.8.4.6 \& Cor.7.1.3.4]{ha}. It follows that 
\[
(\CAlg(\cstablecat)^\op,\tau,\Mod_{(\bs)}(\cstablecat))
\]
defines a quasi-coherent context for any pretopology $\tau$ on $\CAlg(\cstablecat)^\op$ with respect to which (\ref{eq:modulecategoryfunctor}) satisfies descent. When this holds, we call it the \emph{tautological quasi-coherent sheaf context} of $(\csite,\tau)$.
\end{example}

If $f\colon \Y \to \X$ is a morphism of stacks, then we obtain a map $f^\ast \colon \QCoh(\Y) \to \QCoh(\X)$ in $\CAlg(\Pr^L_{\st})$, which we call \emph{pullback} along $f$. We denote by $f_\ast\colon \QCoh(\X)\to \QCoh(\Y)$ the right adjoint of $f^*$.
    
Associated to a quasi-coherent sheaf context is a \emph{global sections functor}.

\begin{prop}\label{pr:endomorphismeinftyring}
Given a quasi-coherent sheaf context $(\csite, \tau, \sfQ)$, where $\csite$ has terminal object $\ast$, the functor
\[\CAlg(\sfQ(\ast)) \to \CAlg(\PrLst)_{\sfQ(\ast)/},\qquad \A \mapsto \Mod_\A(\sfQ(\ast))\]
admits a right adjoint given sending a map $f\colon \sfQ(\ast) \to \sfD$ of symmetric monoidal categories to the pushforward $f_\ast \1_\dcat \in \sfQ(\ast)$ of the unit of $\sfD$.
\end{prop}

\begin{proof}
By the proof of \cite[Th.4.8.5.11]{ha} and the following paragraph, the functor
\[\Mod_{(-)}\colon \mathrm{Alg}(\sfQ(\ast)) \to \LMod_{\sfQ(\ast)} (\PrLst)_{\sfQ(\ast)/}\]
admits a right adjoint. By \cite[Th.4.8.5.16]{ha}, this functor is strong monoidal. The desired adjunction then follows by applying $\CAlg(-)$.
\end{proof}

\begin{mydef}\label{def:globalsectionsonstacks}
    Given a quasi-coherent sheaf content $(\csite, \tau, \sfQ)$, we define the \emph{global sections functor} as the composite
    \[\Ga\colon \Stk_\tau(\csite)^\op \xrightarrow{\QCoh} \CAlg(\PrLst)_{\sfQ(\ast)/} \to \CAlg(\sfQ(\ast)),\]
    where the latter functor is the right adjoint of \Cref{pr:endomorphismeinftyring}.
\end{mydef}

By construction, we have $\Ga(\X) = \lim \map_{\sfQ(\A)}(\1,\1) = \lim \Ga(\A)$ taken over a presentation of $\X$ as a colimit of representables $\Spec \A$. Given a fixed stack $\X\in \Stk_\tau(\csite)$, there is also a global sections functor relative to the given quasi-coherent sheaf context.

\begin{mydef}\label{def:globalsectionsonqcoh}
    Let $(\csite, \tau, \sfQ)$ be a quasi-coherent sheaf content and $\X\in\Stk_\tau(\csite)$ be a stack. Then the counit of the adjunction of \Cref{pr:endomorphismeinftyring} evaluated on $\X$ induces a colimit preserving symmetric monoidal functor $\Mod_{\Ga(\X)}(\sfQ(\ast)) \to \QCoh(\X)$ as a morphism in $\CAlg(\PrLst)$. The right adjoint to this functor
    \[\Ga\colon \QCoh(\X) \to \Mod_{\Ga(\X)}(\sfQ(\ast))\]
    is called the \emph{global sections functor} on $\QCoh(\X)$. It is canonically a lax symmetric monoidal functor by \cite[Pr.2.5.5.1]{sag}.
\end{mydef}

To avoid overcrowded notation, we write $\Ga$ for both of the functors appearing in \Cref{def:globalsectionsonstacks,def:globalsectionsonqcoh}. The different (co)domains of these functors should avoid any confusion.

\subsection{Affineness properties of morphisms between stacks}\label{ssec:genericaffineness}
The goal of this subsection is to study various forms of ``affineness'' for morphisms $\Y \to \X$ between general stacks. To this end, we fix a quasi-coherent sheaf context $(\csite, \tau,\sfQ)$, and we abbreviate $\Stk_\tau(\csite) = \Stk$.

\subsubsection{0-affine morphisms}\label{sssec:0affine}
We start with a geometric reformulation of \cref{def:affinefunctor}.

\begin{mydef}
A map $f\colon \Y \to \X$ in $\Stk$ is \emph{$0$-affine} if $f_\ast\colon \QCoh(\Y)\to\QCoh(\X)$ induces a (necessarily symmetric monoidal) equivalence
\[
\QCoh(\Y) \simeq \Mod_{f_\ast\calO_\Y}(\QCoh(\X)).
\]
A stack $\X$ itself is \emph{$0$-affine} if the terminal map $\X \to \ast$ is $0$-affine. In other words, $\C$ is $0$-affine if the global sections functor induces an equivalence
\[
\Gamma\colon \QCoh(\X) \to \Mod_{\Gamma(\X)}(\sfQ(\ast)).
\]
\end{mydef}

Our goal in this subsection is to record some properties of $0$-affine morphisms, mostly drawing from Barthel--Carmeli--Schlank--Yanovski \cite[\textsection 2]{barthelcarmelischlankyanovski2024chromatic}. We begin with necessary and sufficient conditions for a morphism to be $0$-affine.

\begin{mydef}
A map $f\colon \Y \to \X$ in $\Stk$ \emph{satisfies the projection formula} if for all $\calM \in \QCoh(\X)$ and $\calN \in \QCoh(\Y)$, the natural composite
\[
f_\ast(\calN) \otimes_{\calO_\X}\calM \to f_\ast(\calN)\otimes_{\calO_\X} f_\ast f^\ast \calM \to f_\ast(\calN \otimes_{\calO_\Y} f^\ast\calM)
\]
is an equivalence.
\end{mydef}

\begin{remark}\label{rmk:unit_tensoring}
As a special case, if $f$ satisfies the projection formula, then we havethen
\[
f_*f^*\calM \simeq f_* (f^*\calO_\sfX \otimes_{\calO_\sfY} f^*\calM) \simeq f_*f^*\calO_{\sfX}\otimes_{\calO_\sfX} \calM. 
\]
So the composite $f_*f^*(-)$ is canonically identified with tensoring by $f_*f^*\calO_{\sfX}$. 
\end{remark}

\begin{prop}\label{prop:monoidalbarrbeck}
A map $f\colon \Y \to \X$ in $\Stk$ is $0$-affine if and only if the following three conditions hold:
\begin{enumerate}
\item $f_\ast$ is conservative;
\item $f_\ast$ preserves colimits;
\item $f$ satisfies the projection formula.
\end{enumerate}
Moreover, if $\X$ is $0$-affine and $\sfQ(\ast)$ is generated by the monoidal unit, then condition (2) implies condition (3).
\end{prop}
\begin{proof}
This first statement is sometimes called the monoidal Barr--Beck theorem; a proof appears in \cite[Pr.5.29]{nilpotenceanddescent}. The second claim follows as if $f_\ast$ preserves colimits and $\calN \in \QCoh(\Y)$ is fixed, then the set of $\calM\in\QCoh(\X)$ for which the projection map
\[
f_\ast(\calN) \otimes_{\calO_\X}\calM \to f_\ast(\calN \otimes_{\calO_\Y} f^\ast\calM)
\]
is an equivalence for all $\calN \in \QCoh(\Y)$ forms a localizing subcategory of $\QCoh(\X)$ containing the unit $\calO_\X$. In particular, if $\X$ is $0$-affine, then this map is an equivalence for all $\calN \in \QCoh(\Y)$ as $\QCoh(\X) \simeq \Mod_{\Ga(\O_\X)}(\sfQ(\ast))$ is again generated by $\calO_{\X}$, meaning that $f$ satisfies the projection formula.
\end{proof}

This motivates the following intermediate notion.

\begin{mydef}
A morphism $f\colon \Y \to \X$ of stacks is \emph{$0$-semiaffine} if $f_\ast\colon \QCoh(\X)\to\QCoh(\Y)$ preserves colimits and $f$ satisfies the projection formula.
\end{mydef}

Both $0$-affine and $0$-semiaffine morphisms satisfy reasonable closure properties.

\begin{prop}[{\cite[Pr.2.8]{barthelcarmelischlankyanovski2024chromatic}}]\label{prop:0affinecancellation}
Let $g\colon \sfZ \to \Y$ and $f\colon \Y \to \X$ be maps of stacks.
\begin{enumerate}
\item If $f$ is $0$-affine, then $g$ is $0$-(semi)affine if and only if $f\circ g$ is $0$-(semi)affine,
\item If $f$ and $g$ are $0$-semiaffine, then $f\circ g$ is $0$-semiaffine.
\end{enumerate}
\end{prop}
\begin{proof}
(1)~~Consider the functors
\begin{center}\begin{tikzcd}
\QCoh(\sfZ)\ar[r,"g_\ast"]&\QCoh(\Y)\ar[r,"f_\ast"]&\QCoh(\X)
\end{tikzcd}.\end{center}
By assumption, $f$ is $0$-affine, so by \Cref{prop:monoidalbarrbeck}, $f_\ast$ is conservative, preserves colimits, and satisfies the projection formula. It follows immediately that $g_\ast$ is conservative and preserves colimits if and only if the same is true of $f_\ast \circ g_\ast = (f\circ g)_\ast$, proving that $g$ is $0$-semiaffine if and only if $f\circ g$ is $0$-semiaffine. Given $\calM \in \QCoh(\X)$ and $\calN \in \QCoh(\sfZ)$, the projection map for the composite $f\circ g$ factors as a composite of projection maps for $f$ and $g$:
\begin{center}\begin{tikzcd}
(f\circ g)_\ast \calM \otimes_{\calO_\sfZ} \calN \ar[rr]\ar[d,"\simeq"]&&(f\circ g)_\ast (\calM \otimes_{\calO_\sfX} (f\circ g)^\ast \calN)\ar[d,"\simeq"]\\
f_\ast g_\ast \calM \otimes_{\calO_\sfZ} \calN \ar[r]& f_\ast(g_\ast\calM \otimes_{\calO_{\sfY}} f^\ast\calN) \ar[r]& f_\ast g_\ast(\calM \otimes_{\calO_\sfX} g^\ast f^\ast \calN)
\end{tikzcd}.\end{center}
As $f_\ast$ is conservative, $g$ satisfies the projection formula if and only if the first bottom map is an equivalence for all $\calM$ and $\calN$. As the second bottom map is an equivalence by assumption, it follows that $g$ satisfies the projection formula if and only if the same is true of $f\circ g$. This proves that $g$ is $0$-affine if and only if $f\circ g$ is $0$-affine courtesy of \Cref{prop:monoidalbarrbeck}.

(2)~~If $f_\ast$ and $g_\ast$ preserve colimits, then so too does $f_\ast \circ g_\ast = (f\circ g)_\ast$; and if $f$ and $g$ satisfy the projection formula, then the above diagram shows that the same is true of $f\circ g$.
\end{proof}

\begin{cor}
Let $\Stk_{/\X}^{0\aff}\subseteq\Stk_{/\X}$ denote the full subcategory of stacks which are $0$-affine over $\X$. Then any morphism in $\Stk_{/\sfX}^{0\aff}$ is itself $0$-affine.
\end{cor}
\begin{proof}
This is just a rephrasing of \cref{prop:0affinecancellation}(1) in the $0$-affine case.
\end{proof}

We next consider limits of $0$-(semi)affine morphisms.

\begin{mydef}\label{def:adjointable}
A commutative square
\begin{center}\begin{tikzcd}
\Y'\ar[r,"f'"]\ar[d,"h'"]&\X'\ar[d,"h"]\\
\Y\ar[r,"f"]&\X
\end{tikzcd}\end{center}
of stacks is \emph{adjointable} if the \emph{Beck--Chevalley transformation}
\[
h^* f_* \xRightarrow{\eta} f'_*f'{}^* h^* f_* \simeq f'_*h'{}^*f^*f_* \xRightarrow{\epsilon} f'_*h'{}^*
\]
filling the square
\begin{center}\begin{tikzcd}
\QCoh(\Y')\ar[r,"f'_\ast"]&\QCoh(\X')\\
\QCoh(\Y)\ar[u,"h'{ }^\ast"']\ar[r,"f_\ast"]&\QCoh(\X)\ar[u,"h^\ast"]\ar[ul,Rightarrow,shorten > = 1em, shorten < = 1em]
\end{tikzcd}\end{center}
is an equivalence. In this case, the Beck--Chevalley transformation is our canonical choice of witness for the commutativity of this square of categories.
\end{mydef}

\begin{lemma}\label{prop:limitadjointable}
Let $\jdiagram$ be a small category and $f_{(\bs)}\colon \Y_{(\bs)} \to \X_{(\bs)}$ be a natural transformation between $\jdiagram$-shaped diagrams in $\Stk$ satisfying the following conditions:
\begin{enumerate}
\item The square obtained by restricting $f_{(\bs)}$ to any arrow $[1] \to \jdiagram$ is adjointable.
\item For all $j\in \jdiagram$, the map $f_j\colon \Y_j\to\X_j$ is $0$-(semi)affine.
\end{enumerate}
Then the induced map
\[
\colim_{j\in\jdiagram}f_{j}\colon \colim_{j\in \jdiagram} \Y_j \to \colim_{j\in \jdiagram} \X_j
\]
is $0$-(semi)affine, and for all $i\in \jdiagram$ the square
\begin{equation}\label{eq:adjointablecone}\begin{tikzcd}[column sep=large]
\Y_i\ar[r,"f_i"]\ar[d]&\X_i\ar[d]\\
\colim_{j\in \jdiagram}\Y_j\ar[r,"\colim f_j"]&\colim_{j\in \jdiagram}\X_j
\end{tikzcd}\end{equation}
is adjointable.
\end{lemma}
\begin{proof}
By \cite[Cor.4.7.4.18]{ha}, adjointability guarantees that the pushforwards $f_{j\ast}$ assemble into a natural transformation $f_{(\bs)\ast}\colon \QCoh(\Y_{(\bs)}) \to \QCoh(\X_{(\bs)})$ whose limit is the pushforward of $\colim_{j\in\jdiagram}f_{j}$, and that (\ref{eq:adjointablecone}) is adjointable. The remaining claims follow as in the proof of \cite[Pr.2.9]{barthelcarmelischlankyanovski2024chromatic}.
\end{proof}

We end by considering the interaction between $0$-affineness and adjointability.

\begin{prop}\label{prop:adjointablesquare0affine}
Suppose given a commutative square of stacks
\begin{equation}\label{eq:adjointablesquare0affine}\begin{tikzcd}
\Y'\ar[r,"f'"]\ar[d,"h'"]&\X'\ar[d,"h"]\\
\Y\ar[r,"f"]&\X
\end{tikzcd}.\end{equation}
\begin{enumerate}
\item Suppose that $f$ and $f'$ are $0$-affine. Then the following are equivalent:
\begin{enumerate}
\item (\ref{eq:adjointablesquare0affine}) is adjointable.
\item The Beck--Chevalley map
\[
h^\ast f_\ast \calO_\sfY \to f'_\ast h'{ }^\ast \calO_\sfY
\]
is an equivalence.
\item (\ref{eq:adjointablesquare0affine}) is sent by $\QCoh(\bs)$ to a coCartesian square in $\CAlg(\PrLst)$.
\end{enumerate}
\item Suppose that $h$ and $h'$ are $0$-affine. Then the Beck--Chevalley map of (\ref{eq:adjointablesquare0affine}) evaluated on $\calM \in \QCoh(\sfY)$ may be identified, after applying $h_\ast$, with a map
\[
h_\ast\calO_{\X'} \otimes_{\calO_\X} f_\ast \calM \to f_\ast(h_\ast'\calO_{\Y'} \otimes_{\calO_\Y}\calM).
\]
\item Suppose that $f$, $f'$, $h$, and $h'$ are $0$-affine, and abbreviate $d = f\circ h' = h \circ f'$. Then (\ref{eq:adjointablesquare0affine}) is adjointable if and only if the square
\begin{center}\begin{tikzcd}
\calO_\X\ar[r]\ar[d]&h_\ast\calO_{\X'}\ar[d]\\
f_\ast\calO_{\Y}\ar[r]&d_\ast\calO_{\Y'}
\end{tikzcd}\end{center}
is coCartesian in $\CAlg(\QCoh(\X))$.
\end{enumerate}
\end{prop}
\begin{proof}
(1a)$\Rightarrow$(1b) is clear.

(1b)$\Leftrightarrow$(1c)~~As $f$ and $f'$ are $0$-affine, we may identify the square in $\CAlg(\PrLst)$ obtained by applying $\QCoh(\bs)$ to (\ref{eq:adjointablesquare0affine}) as 
\begin{center}\begin{tikzcd}[column sep=huge]
\Mod_{f'_\ast\calO_{\Y'}}(\QCoh(\X'))&\QCoh(\X')\ar[l,"f'_\ast\calO_{\Y'}\otimes_{\calO_{\X'}}(\bs)"']\\
\Mod_{f_\ast\calO_{\Y}}(\QCoh(\X))\ar[u,"h'{ }^\ast"']&\QCoh(\X)\ar[u,"h^\ast"']\ar[l,"f_\ast\calO_\Y\otimes_{\calO_\X}(\bs)"']
\end{tikzcd}.\end{center}
By \cite[Th.4.8.4.6]{ha}, the comparison map $\QCoh(\Y)\otimes_{\QCoh(\X)}\QCoh(\X') \to \QCoh(\Y')$ may be identified with the map
\[
\Mod_{h^\ast f_\ast \calO_{\Y}}(\QCoh(\X')) \to \Mod_{f'_\ast h'{ }^\ast\calO_\Y}(\QCoh(\X'))
\]
obtained from the Beck--Chevalley map $h^\ast f_\ast \calO_\Y \to f'_\ast h'{ }^\ast \calO_\Y$. Therefore, $\QCoh(\Y)\otimes_{\QCoh(\X)}\QCoh(\X') \to \QCoh(\Y')$ is an equivalence if and only if $h^\ast f_\ast \calO_\Y \to f'_\ast h'{ }^\ast \calO_\Y$ is an equivalence, proving (1b)$\Leftrightarrow$(1c).

(1b)$\Rightarrow$(1a)~~The symmetric monoidal functor $h^\ast\colon \QCoh(\X)\to\QCoh(\X')$ extends to a symmetric monoidal functor
\[h^\ast\colon \Mod_{f_\ast\calO_\Y}(\QCoh(\X)) \to \Mod_{h^\ast f_\ast\calO_\Y}(\QCoh(\X')),\]
which yields a commutative square
\begin{center}\begin{tikzcd}[column sep=huge]
\Mod_{h^\ast f_\ast\calO_\Y}(\QCoh(\X'))&\QCoh(\X')\ar[l,"h^\ast f_\ast \calO_\Y \otimes_{\calO_{\X'}}(\bs)"']\\
\Mod_{f_\ast\calO_\Y}(\QCoh(\X))\ar[u,"h^\ast"]&\QCoh(\X)\ar[l,"f_\ast\calO_\Y\otimes_{\calO_\X}(\bs)"']\ar[u,"h^\ast"]
\end{tikzcd}.\end{center}
As above, if $h^\ast f_\ast \calO_\Y \to f'_\ast h'{ }^\ast \calO_\Y$ is an equivalence, then this square canonically commutes, and the filler can be identified with the Beck--Chevalley transformation of  \Cref{def:adjointable}, implying (1b)$\Rightarrow$(1a).

(2)~~Since $h$ and $h'$ are $0$-affine, we obtain identifications
\[
h_\ast h^\ast f_\ast\calM\simeq h_\ast\calO_{\X'}\otimes_{\calO_\X} f_\ast\calM \quad\text{and}\quad h_\ast f'_\ast h'{ }^\ast \calM\simeq f_\ast h'_\ast h'{ }^\ast \calM \simeq f_\ast(h'_\ast\calO_{\Y'}\otimes_{\calO_\Y}\calM),
\]
using \cref{rmk:unit_tensoring} in each case.

(3)~~As $f$ and $f'$ are $0$-affine, (1) implies that the square is adjointable if and only if the Beck--Chevalley transformation evaluates to an equivalence on the $\calO_\Y$. As $h$ and $h'$ are $0$-affine, (2) identifies this Beck--Chevalley map as equivalent, after applying $h_\ast$, to the map
\[
h_\ast\calO_{\sfX'}\otimes_{\calO_\sfX} f_\ast\calO_{\sfY} \to d_\ast \calO_\sfY
\]
which is an equivalence if and only if the square of (3) is coCartesian in $\CAlg(\QCoh(\X))$.
\end{proof}

\subsubsection{Universally 0-affine morphisms}

In general, $0$-affine morphisms between stacks may fail to have good local-to-global and base change properties. We fix this by introducing the following refinement.

\begin{mydef}\label{def:univ0affine}
A morphism $f\colon \Y \to \X$ is \emph{universally $0$-affine} if for all diagrams of the form
    \begin{equation}\label{eq:defofU0A}\begin{tikzcd}
    \sfF'\ar[r]\ar[d]&\Spec \sfA'\ar[d]\\
    \sfF\ar[r]\ar[d]&\Spec \sfA\ar[d]\\
    \sfY\ar[r,"f"]&\sfX
    \end{tikzcd},\end{equation}
    where both squares are Cartesian, the map $\sfF \to \Spec \sfA$ is $0$-affine and the top square is adjointable.
\end{mydef}

\begin{remark}\label{rmk:internalandexternalU0Adef}
By \cref{prop:adjointablesquare0affine}, the top square in \cref{def:univ0affine} is adjointable if and only if, writing $f\colon \sfF \to \Spec\sfA$ and $f'\colon \sfF' \to \Spec\sfA$, the comparison map
\[
f_\ast \calO_{\sfF} \otimes_\sfA \sfA' \to f'_\ast \calO_{\sfF'}
\]
is an equivalence.
\end{remark}

Before we prove some basic facts about universally $0$-affine morphisms, let us recall the definition of an affine morphism.

\begin{mydef}\label{def:affinemorphism}
    A map of stacks $f\colon \Y \to \X$ is \emph{affine} if for each $\Spec \A \to \X$, the pullback $\Y\times_\X \Spec \A$ is equivalent to an affine stack.
\end{mydef}

\begin{prop}\label{prop:univ0affinefacts}
Fix a morphism $f\colon \Y \to \X$.
\begin{enumerate}
\item Suppose that $f$ is universally $0$-affine, and suppose we are given a Cartesian square
\begin{equation}\label{eq:un_0_aff_bc}\begin{tikzcd}
\sfY'\ar[r,"f'"]\ar[d,"h"]&\sfX'\ar[d,"h'"]\\
\sfY\ar[r,"f"]&\sfX
\end{tikzcd}\end{equation}
in $\Stk$. Then $f'$ is universally $0$-affine and this square is adjointable.
\item The following implications hold:
\[
f\text{ is affine} \quad\implies\quad f\text{ is universally $0$-affine}\quad\implies\quad f\text{ is $0$-affine}.
\]
\item If $f$ is universally $0$-affine and $g\colon \sfZ\rightarrow\sfY$ is universally $0$-affine, then $f\circ g \colon \sfZ\rightarrow\sfX$ is universally $0$-affine.
\end{enumerate}
\end{prop}
\begin{proof}
(1,2)~~First, note that universally $0$-affine morphisms are clearly stable under base change, and that affine morphisms are universally $0$-affine by the definition of a quasi-coherent context. So to prove claim (1) and (2) together, it remains to prove in the situation of (1) that $f$ is $0$-affine and (\ref{eq:un_0_aff_bc}) is adjointable.

First, consider the square in the case where $\X' = \Spec\sfA$ is affine. Note that $\sfX$ admits a presentation of the form $\sfX \simeq \colim \Spec\sfA_\bullet$ whose universal cocone contains the map $h\colon \Spec\sfA \to \sfX$. The map $\sfY \to \sfX$ may be identified as the colimit of the natural transformation
\[
\sfY\times_\sfX \Spec\sfA_\bullet \to \Spec\sfA_\bullet.
\]
As $f$ is universally $0$-affine, this natural transformation satisfies the conditions of \cref{prop:limitadjointable}, proving that $\sfY \to \sfX$ is $0$-affine and that (\ref{eq:un_0_aff_bc}) is adjointable when $\X'$ is affine.

It remains to prove that (\ref{eq:un_0_aff_bc}) is adjointable for general $\X'$. By writing $\X'$ as a colimit of affines, we see that to prove that the Beck--Chevalley transformation $\alpha\colon h^\ast f_\ast \to f_\ast' h'{ }^\ast $ is an equivalence, it suffices to prove that $i^\ast\alpha$ is an equivalence for any $i\colon \Spec\A \to \X'$. Form the Cartesian diagram
\begin{center}\begin{tikzcd}
\sfF\ar[r,"f''"]\ar[d,"j"]&\Spec\sfA\ar[d,"i"]\\
\Y'\ar[r,"f'"]\ar[d,"h'"]&\X'\ar[d,"h"]\\
\Y\ar[r,"f"]&\X
\end{tikzcd}.\end{center}
Since Beck--Chevalley transformations compose, see e.g.\ \cite[Appendix~A.6(4)] {CLLPartial}, the composite of $i^*\alpha$ and the Beck--Chevalley transformation of the top square is homotopic to the Beck--Chevalley transformation of the outer square. As $f$ and $f'$ are universally $0$-affine, both the latter transformations are equivalences. We conclude that $i^\ast\alpha$ is also an equivalence.

(3)~~Fix maps $\Spec \sfA'\rightarrow\Spec \sfA\rightarrow \sfX$, and form the Cartesian diagram
\begin{center}\begin{tikzcd}
\sfG'\ar[r]\ar[d]&\sfF'\ar[r]\ar[d]&\Spec\sfA'\ar[d]\\
\sfG\ar[r,"g'"]\ar[d]&\sfF\ar[r,"f'"]\ar[d]&\Spec\sfA\ar[d]\\
\sfZ\ar[r,"g"]&\sfY\ar[r,"f"]&\sfX
\end{tikzcd}.\end{center}
To prove that $fg$ is universally $0$-affine, we must prove that $f'g'$ is $0$-affine and that the top outer rectangle is adjointable.

As $f$ is universally $0$-affine, $f'$ is $0$-affine by definition. As $g$ is universally $0$-affine, $g'$ is universally $0$-affine and hence $0$-affine by (2). Therefore, $f'g'$ is $0$-affine by \cref{prop:0affinecancellation}.

As $f$ is universally $0$-affine, the top right square is adjointable by definition. As $g'$ is universally $0$-affine, the top left square is adjointable by (1). As a pasting of adjointable squares is adjointable, it follows that the top outer rectangle is adjointable.
\end{proof}

\subsubsection{Affine morphisms}
In some special cases, such as in \Cref{sec:nonconsag}, the slice category of stacks that are affine over a fixed stack $\X$ admits a simple description.

\begin{theorem}\label{relativespectrumtheorem}
    Let $(\csite, \tau, \sfQ)$ be a tautological quasi-coherent sheaf context (\Cref{ex:basicqcohcontext}), meaning it is of the form $\csite = \CAlg(\ccat)^\op$ for some $\ccat\in \CAlg(\PrLst)$ and $\sfQ(-) = \Mod_{(-)}(\ccat)$. Then, for any fixed stack $\X$, the lax monoidal functor
    \[\Phi\colon \Stk^\aff_{/\X} \to \CAlg(\QCoh(\X))^\op, \qquad (f\colon \Y\to \X) \mapsto f_\ast \O_\Y\]
    is an equivalence of categories. In particular, it is strong monoidal.
\end{theorem}

Lurie states and proves a version of this for connective spectral Deligne--Mumford stacks in \cite[Pr.2.5.1.2]{sag}. 

\begin{proof}
    Writing $\X \simeq \colim \Spec \A$ as a colimit of representables, the equivalences of categories
    \[\QCoh(\X) \simeq \lim \QCoh(\Spec \A), \qquad \Stk_{/\X} \simeq \lim \Stk_{/\Spec \A}\]
    by definition and \Cref{pr:stackisbasicallyatopos} respectively, induce equivalences 
    \[\Stk^\aff_{/\X} \simeq \lim \Stk^\aff_{/\Spec \A}, \qquad \CAlg(\QCoh(\X)) \simeq \lim \CAlg(\QCoh(\Spec \A))\]
    courtesy of the fact the functor $\CAlg(-)\colon \CAlg(\Pr^L_\st) \to \Pr^L_\st$ preserves limits by construction. This reduces us to the case where $\X=\Spec \A$ is affine, where the functor $\Phi$ has an obvious inverse:
    \[\CAlg(\QCoh(\Spec \A))^\op \simeq \CAlg_\A^\op \simeq (\CAlg^{\op})_{/\A}\underset{\sim}{\xrightarrow{\Spec(-)}} \Stk^\aff_{/\Spec \A}\qedhere\]
\end{proof}

\begin{remark}\label{rmk:relativespectra}
    The theorem above defines \emph{relative spectra}. In particular, given some $\calA\in \CAlg(\QCoh(\X))$, we define its relative spectrum $\relSpec_\X(\calA)$ as the unique affine stack over $\X$ whose image under $\Phi$ is $\calA$. For example, the complex-periodification of an $\E_\infty$ ring $\A$, see \Cref{maintheorem:moaoriaffineness} and \Cref{sssec:maori}, can be written as
    \[\M_\A^\ori = \relSpec_{\M_\FG^\ori}(\O_{\M_\FG^\ori} \otimes \A).\]
\end{remark}

Combining this result with standard results from \cite{ha} yields the following corollary.

\begin{cor}\label{tannakasteppingstone}
    Using the notation of \Cref{relativespectrumtheorem}, the assignment
    \[\Theta_{\QCoh}\colon (\Stk_\tau(\calC)^\aff_{/\Y})^\op \to \CAlg(\Pr^L_\st)_{\QCoh(\Y)/}, \qquad (f\colon \sfZ \to \Y) \mapsto (f^\ast\colon \QCoh(\Y) \to \QCoh(\sfZ))\]
    is fully faithful with essential image those symmetric monoidal colimit preserving functors $F\colon \QCoh(\Y) \to \QCoh(\sfZ)$, with right adjoint $G$, which satisfy the following conditions:
    \begin{enumerate}
        \item The functor $G$ preserves geometric realizations of simplicial objects.
        \item The functor $G$ is conservative.
        \item The adjunction $F\dashv G$ satisfies the projection formula.
    \end{enumerate}
\end{cor}

When $\calD = \CAlg^\cn$, Lurie proves a stronger statement, showing that the above functor is fully faithful on those morphisms $\sfZ \to \Y$ which are \emph{quasi-affine}. Many of his arguments also work for more general $\calD$ if the condition of being quasi-affine is replaced with universally $0$-affine, except for the crucial \cite[Cor.2.4.2.2]{sag}, which is the missing reference in the proof of \cite[Lm.9.2.1.2]{sag} (compare the proof of \cite[Lm.3.2.8]{dagviii}).

\begin{proof}
    This is precisely an application of \cite[Cor.4.8.5.21]{ha}. Indeed, by \Cref{relativespectrumtheorem}, we have an equivalence of categories
        \[\Phi\colon \Stk_\tau(\calD)^\aff_{/\Y} \xrightarrow{\simeq} \CAlg(\QCoh(\Y))^\op, \qquad (f\colon \sfZ\to \Y) \mapsto f_\ast \O_\sfZ.\]
    Moreover, for each affine map $f\colon \sfZ \to \Y$ we have a natural equivalence of symmetric monoidal categories $\QCoh(\sfZ) \simeq \Mod_{f_\ast \O_\sfZ}\QCoh(\Y)$ courtesy of \Cref{prop:univ0affinefacts}. The functor in question can be written as the composite
    \[(\Stk_\tau(\calD)^\aff_{/\Y})^\op \xrightarrow{\simeq} \CAlg(\QCoh(\Y)) \xrightarrow{\Theta} \CAlg(\Pr^L_\st)_{\QCoh(\Y)/},\]
    where $\Theta$ sends $\calA \in \CAlg(\QCoh(\Y))$ to $\Mod_\calA (\QCoh(\Y))$. According to \cite[Cor.4.8.5.21]{ha}, this functor $\Theta$ is fully faithful with essential image precisely as described.
\end{proof}

\subsubsection{Weak Tannaka duality}\label{sec:appendix}
In \cite{dagviii} and \cite[\textsection9]{sag}, Lurie discusses \emph{Tannaka duality}, a statement connecting maps between stacks to functors between the associated categories of quasi-coherent sheaves. Here, we state a weak version of Tannaka duality in the generality of this section. We say that a stack $\X$ has \emph{affine diagonal} if the diagonal map $\X \to \X\times \X$ is affine in the sense of \Cref{def:affinemorphism}.

\begin{theorem}\label{tannakadualityforgeometricstacks}
    Let $\sfC \in \CAlg(\PrLst)$, $\tau$ be a Grothendieck topology on $\CAlg(\sfC)^\op = \calD$ such that $(\calD,\tau, \Mod_{(-)}(\sfC))$ is a tautological quasi-coherent sheaf context (\Cref{ex:basicqcohcontext}), and $\X,\Y$ be stacks in $\Stk_\tau(\calD)$ such that $\X$ has affine diagonal. Then the assignment
    \[\Xi\colon \Map_{\Stk_\tau(\calD)}(\Y, \X) \hookrightarrow \Fun^\otimes(\QCoh(\X), \QCoh(\Y)), \quad (f\colon \Y \to \X) \mapsto (f^\ast\colon \QCoh(\X) \to \QCoh(\Y))\]
    is fully faithful.
\end{theorem}

We follow Lurie in calling this \emph{weak} Tannaka duality, as the desired functor is fully faithful without any expression of the essential image. A stronger version would give global conditions on a functor $\QCoh(\X) \to \QCoh(\Y)$ for it to be in the essential image of $\Xi$; see \cite[\textsection9]{sag} for examples of this where $\sfC = \Sp^\cn$.

\begin{proof}[Proof of \Cref{tannakadualityforgeometricstacks}]
    We start as the proof of \cite[Pr.3.3.11]{dagviii} or \cite[Pr.9.2.2.1]{sag}. For our fixed $\X$, the assignments
    \[\Y \mapsto \Map_{\Stk}(\Y, \X), \qquad \Y \mapsto \Fun^\otimes(\QCoh(\X), \QCoh(\Y))\]
    both send colimits in $\Y$ to limits of categories, so we are reduced to $\Y=\Spec\A$. Fixing two maps $f,g\colon \Spec \A\to \X$, we first note that by \cite[Rmk.4.8.5.9]{ha}, any linear symmetric monoidal natural transformation $\al\colon f^\ast \Rightarrow g^\ast$ is automatically a natural equivalence. In particular, the category $\Fun^\otimes(\QCoh(\X), \Mod_\A(\sfC))$ is actually a groupoid. To show that $\Xi$ is fully faithful, we want to show that from the commutative diagram of spaces
    \[\begin{tikzcd}
        {\Map_{\Stk_\tau(\calD)_{/\X}}(f,g)}\ar[r]\ar[d]   &   {\Map_{\CAlg(\Pr^L_\st)_{\QCoh(\X)/}}(g^\ast, f^\ast)}\ar[d]  \\
        {\Map_{\Stk_\tau(\calD)}(\Spec \A,\Spec \A)}\ar[r] &   {\Map_{\CAlg(\Pr^L_\st)_{\sfC/}}(\Mod_\A(\sfC), \Mod_\A(\sfC)),}
    \end{tikzcd}\]
    that the induced map on vertical fibers is an equivalence. It then suffices to show that the two horizontal maps above are equivalences. This is clear by two applications of \Cref{tannakasteppingstone}, once over $\X$ and once over the terminal stack of $\Stk_\tau(\calD)$. This shows that $\Xi$ is fully faithful.
\end{proof}

\subsection{Descent stacks}\label{ssec:descentstacks}
The main general construction of this article, generalizing the moduli stack of oriented formal groups, is that of a \emph{descent stack}. Again, for this section, we fix a quasi-coherent sheaf context $(\csite, \tau, \sfQ)$ and write $\Stk = \Stk_\tau(\csite)$.

\subsubsection{The descent stack factorization}

\begin{mydef}\label{def:descentstacks}
The \textit{descent stack} $\sfD_f$ of a morphism $f\colon \sfY\rightarrow\sfX$ is its image, i.e., its $(-1)$-truncation in the slice category $\Stk_{/\sfX}$. This sits in the effective epimorphism-monomorphism factorization
\begin{center}\begin{tikzcd}
\sfY\ar[r, two heads]&\sfD_f\ar[r, tail]&\sfX.
\end{tikzcd}\end{center}
\end{mydef}

\begin{lemma}\label{lem:descentstackbasicproperties}
Fix a morphism $f\colon \sfY\rightarrow\sfX$ of stacks.
\begin{enumerate}
\item\label{item:descent_\v{C}ech_nerve} The descent stack $\sfD_f$ is the geometric realization of the \v{C}ech nerve of $f$, i.e.\
\[
\sfD_f = \colim \left(\begin{tikzcd} \Y&\Y\times_\X \Y\ar[l,shift left=1mm]\ar[l,shift right=1mm]&\cdots\ar[l,shift left=1mm]\ar[l]\ar[l,shift right=1mm]\end{tikzcd}\right).
\]
\item\label{item:restrictdescentstack} Given a Cartesian square
\begin{center}\begin{tikzcd}
\Y'\ar[r,"f'"]\ar[d]&\X'\ar[d]\\
\Y\ar[r,"f"]&\X
\end{tikzcd},\end{center}
there is a unique equivalence $D_{f'}\simeq \D_f\times_\X\X'$ of stacks over $\X'$. In particular, for any stack $\sfZ$ the natural map $\D_{f\times\sfZ} \to \D_f\times \sfZ$ is an equivalence.
\item\label{item:fullyfaithfulness} The forgetful functor $\Stk_{/\sfD_f}\rightarrow\Stk_{/\sfX}$ is fully faithful, and preserves all colimits and all nonempty limits.
\item\label{item:descent_0_truncated} Given a map $\sfZ\rightarrow\sfX$, the space 
\[
\Map_\sfX(\sfZ,\sfD_f) = \left\{\begin{tikzcd}
&\sfD_f\ar[d]\\
\sfZ\ar[ur,dashed]\ar[r]&\sfX
\end{tikzcd}\right\}
\]
of lifts through $\sfD_f$ is empty or contractible.
\item\label{item:pbdescentstack} The canonical map
$
\sfZ\times_{\sfD_f}\sfZ'\to \sfZ\times_{\sfX}\sfZ'
$
is an equivalence for any $\sfZ,\sfZ'\in \Stk_{/\sfD_f}$.
\item\label{item:descentstackimagecriterion} Given $\sfZ\to\sfX$ in $\in\Stk_{/\sfX}$, the following are equivalent:
\begin{enumerate}
\item $\sfZ\to\sfX$ is in the essential image of $\Stk_{/\sfD_f}\to\Stk_{/\sfX}$.
\item $\sfZ\times_\sfX\sfY \to \sfZ$ is an effective epimorphism.
\item There exists an effective epimorphism $p\colon \sfW\to\sfZ$ and commutative diagram
\begin{center}\begin{tikzcd}
\sfW\ar[r]\ar[d,"p"]&\sfY\ar[d,"f"]\\
\sfZ\ar[r]&\sfX
\end{tikzcd}.\end{center}
\item For some (and thus any) effective epimorphism $\sfW\to\sfZ$, the composite $\sfW\to\sfZ\to\sfX$ is in the essential image of $\Stk_{/\sfD_f}\to\Stk_{/\sfX}$.
\end{enumerate}

\end{enumerate}
\end{lemma}
\begin{proof}
These are all basic properties of image factorizations in topoi that hold just as well in $\Stk$. In order:

(1)~~This follows from the fact that groupoids in $\Stk$ are effective and colimits are universal, see for example \cite[Pr.6.2.3.4]{htt}.

(2)~~This follows from the fact that both monomorphisms and effective epimorphisms are closed under pullback in $\Stk$: see \cite[Pr.5.2.8.6, Pr.6.2.3.15(5)]{htt}.

(3)~~Since $i\colon \sfD_f\to \sfX$ is a monomorphism, the square
\[\begin{tikzcd}
	{\sfD_f} & {\sfD_f} \\
	{\sfD_f} & \sfX
	\arrow[equals, from=1-1, to=1-2]
	\arrow[equals, from=1-1, to=2-1]
	\arrow[from=1-2, to=2-2]
	\arrow[from=2-1, to=2-2]
\end{tikzcd}\]
is Cartesian. This implies that the unit of the adjunction $i_!\dashv i^*$ is an equivalence, and so $i_!$ is fully faithful. It also clearly preserves colimits. For the claim that it preserves non-empty limits, it suffices to show that the limit in $\Stk_{/\sfX}$ of a non-empty diagram $F\colon I\to \Stk_{\sfD_f}$ is again in $\Stk_{\sfD_f}$. This limit admits a map to $F(i)$ for any $i\in I$, and therefore also to $\sfD_f$, and so lies in $\Stk_{\sfD_f}$.

(4)~~This is a reformulation of the statement from (2) that $\Stk_{/\sfD_f} \to \Stk_{\sfX}$ is fully faithful.

(5)~~This asserts that $\Stk_{/\sfD_f} \to \Stk_{/\sfX}$ preserves binary products, which is a special case of (2).

(6b)$\Rightarrow$(6c)$\Rightarrow$(6d)~~Clear.

(6d)$\Rightarrow$(6a)~~As $\sfD_f\to\sfX$ is a monomorphism, this is the fact that effective epimorphisms are left orthogonal to monomorphisms.

(6a)$\Rightarrow$(6b)~~By (5), there is a Cartesian diagram
\begin{center}\begin{tikzcd}
\sfZ\times_{\sfX}\sfY\ar[r]\ar[d]&\sfY\ar[d]\\
\sfZ\ar[r]&\sfD_f
\end{tikzcd}.\end{center}
As $\sfY\to\sfD_f$ is an effective epimorphism, so is $\sfZ\times_{\sfX}\sfY\to\sfZ$.
\end{proof}

We have introduced the descent stack as a tool for studying \emph{Adams descent} for quasi-coherent sheaves. This connection is provided by the following.

\begin{prop}\label{prop:qcohdescentstack}
Suppose that $f\colon \sfY\rightarrow\sfX$ is universally $0$-affine, and write $q\colon \sfY\rightarrow\sfD_f$ for the projection and $j \colon \sfD_f\rightarrow \sfX$ for the inclusion.
\begin{enumerate}
\item $q\colon \sfY\rightarrow\sfD_f$ is universally $0$-affine.
\item The maps
\begin{align*}
\QCoh(\sfD_f)&\rightarrow \lim\left(
\begin{tikzcd}[ampersand replacement=\&]
\QCoh(\sfY)\ar[r,shift left=1mm]\ar[r,shift right=1mm]\&\QCoh(\sfY\times_\sfX \sfY)\ar[r,shift left=1mm]\ar[r]\ar[r,shift right=1mm]\&\cdots\end{tikzcd}
\right)\\
&\leftarrow \lim\left(
\begin{tikzcd}[ampersand replacement=\&]
\QCoh(\sfY)\ar[r,shift left=1mm]\ar[r,shift right=1mm]\&\QCoh(\sfY)\otimes_{\QCoh(\sfX)}\QCoh(\sfY)\ar[r,shift left=1mm]\ar[r]\ar[r,shift right=1mm]\&\cdots\end{tikzcd}
\right)\\
&\rightarrow \lim\left(
\begin{tikzcd}[ampersand replacement=\&]
\Mod_{f_\ast\calO_{\sfY}}(\QCoh(\sfX))\ar[r,shift left=1mm]\ar[r,shift right=1mm]\&\Mod_{f_\ast\calO_{\sfY}\otimes_{\calO_{\sfX}}f_\ast\calO_{\sfY}}(\QCoh(\sfX))\ar[r,shift left=1mm]\ar[r]\ar[r,shift right=1mm]\&\cdots
\end{tikzcd}\right)
\end{align*}
are equivalences.
\item The restriction $q^\ast\colon \QCoh(\sfD_f)\rightarrow\QCoh(\sfY)$ is comonadic and the natural transformation $q^\ast q_\ast \to f^\ast f_\ast$ is an equivalence, thus realizing $\QCoh(\sfD_f)$ as the category of coalgebras for the descent comonad on $\QCoh(\sfY)$ associated to the adjunction $f^\ast \dashv f_\ast$.
\item\label{item:descentstacknilpotentcompletion} If $\calM \in \QCoh(\sfX)$, then the unit $\calM\rightarrow j_\ast j^\ast \calM$ realizes $j_\ast j^\ast$ as the functor of $f_\ast \calO_{\sfY}$-nilpotent completion inside $\QCoh(\sfX)$. In particular, $j_\ast \calO_{\sfD_f}\simeq (\calO_{\sfX})_{\calO_{\sfY}}^\wedge$.
\end{enumerate}
\end{prop}
\begin{proof}
(1)~~This is immediate given \cref{lem:descentstackbasicproperties}(\ref{item:pbdescentstack}).

(2)~~Given the description of $\sfD_f$ as a \v{C}ech nerve in \cref{lem:descentstackbasicproperties}, this follows by combining \cref{prop:univ0affinefacts} and \cref{prop:adjointablesquare0affine}.

(3)~~This follows from (2) and the (levelwise opposite variant of) \cite[Th.4.7.5.2]{ha}.

(4)~~Given $\calM \in \QCoh(\sfX)$, we can use the cobar resolution of $j^\ast\calM$ for the desired computation:
\begin{align*}
j_\ast j^\ast \calM \simeq \Tot \left(j_\ast (q^\ast q_\ast)^{\bullet+1} j^\ast\calM\right)&\simeq \Tot \left((f^\ast f_\ast)^{\bullet+1}\calM\right)\\
&\simeq \Tot \left((f_\ast \calO_\sfY)^{\otimes_{\calO_\sfX}\bullet+1}\otimes_{\calO_X}\calM\right) \simeq \calM_{f_\ast \calO_\sfY}^\wedge.\qedhere
\end{align*}
\end{proof}

\subsubsection{Adams flat morphisms}
The descent stack is particularly well-behaved in the presence of some additional assumptions.

\begin{mydef}\label{df:affinecover}
    A morphism $f\colon \Y \to \X$ of stacks is an \emph{affine cover} if it is affine and for every map $\Spec\A\to\X$, the singleton $\{\Y\times_\X\Spec\A\to\Spec\A\}$ is a cover.
\end{mydef}

The following lemma is immediate from this definition.

\begin{lemma}\label{lem:coverssatisfydescent}
Suppose given a Cartesian diagram
\begin{center}\begin{tikzcd}
\Y'\ar[r,"f'"]\ar[d]&\X'\ar[d,"h"]\\
\Y\ar[r,"f"]&\X
\end{tikzcd}\end{center}
of stacks. If $f$ is an affine cover, then $f'$ is an affine cover, and the converse holds if $h$ is an effective epimorphism.
\end{lemma}

\begin{mydef}\label{def:adamsflat}
    A morphism $f\colon \sfY\rightarrow\sfX$ is \textit{Adams flat} if either projection $\sfY\times_\sfX\sfY \to \sfY$ is an affine cover.
\end{mydef}

The inspiration for this definition comes from stable homotopy theory, where one says that a homotopy commutative ring spectrum $E$ is {Adams flat} if either the left or the right unit maps $E \to E\otimes_\Sph E$ induce flat maps of graded rings.

We will now state a useful criterion for the map $q\colon \sfY \to \sfD_f$ to be an affine cover. For this, we require an additional definition.

\begin{mydef}
We say a site $(\csite, \tau)$ is \emph{generated by singleton coverings} if, for every covering $\{\Spec \A_i \to \Spec\A : i\in I\}$, there exists a subset $I_0 \subseteq I$ such that $\{\Spec \prod_{i\in I_0}\A_i \to \Spec\A\}$ and $\{\Spec\A_i \to \Spec\prod_{i\in I_0}\A_i : i \in I_0\}$ are coverings.
\end{mydef}

\begin{remark}
The fpqc topology on non-connective $\E_\infty$ rings, defined in \cref{def:fpcq_top}, is clearly generated by singleton coverings. 
\end{remark}

\begin{prop}\label{localpropertiesofDpandADAMSflatness}
Let $f\colon \sfY\rightarrow\sfX$ be a morphism of stacks. Furthermore, suppose that $(\csite, \tau)$ is generated by singleton coverings. Then the map $\sfY\rightarrow\sfD_f$ is an affine cover if and only if $f$ is Adams flat.
\end{prop}
\begin{proof}
First, suppose that $\sfY \to \sfD_f$ is an affine cover. By base change, it follows that $\sfY \times_{\sfD_f}\sfY \to \sfY$ is an affine cover. As $\sfY\times_{\sfD_f}\sfY\simeq\sfY\times_{\sfX}\sfY$, it follows that $f$ is Adams flat.

Conversely, suppose that $f$ is Adams flat. To prove that $\sfY \to \sfD_f$ is an affine cover, we must prove that if $h\colon \Spec\A\to\sfD_f$ is arbitrary then
\[
\Spec\A \times_{\sfD_f}\sfY \simeq \Spec\A\times_{\sfX}\sfY \to \Spec\A
\]
is an affine cover. As the topology on $\csite$ is generated by singleton coverings, \cref{lem:descentstackbasicproperties} implies that there exists a covering $\{p\colon \Spec\B\to\Spec\A\}$ and a commutative diagram
\begin{center}\begin{tikzcd}
\Spec\B\ar[r]\ar[d,"p"]&\Y\ar[d,"f"]\\
\Spec\A\ar[r]&\X
\end{tikzcd}.\end{center}
By \cref{lem:coverssatisfydescent}, it now suffices to prove that
\[
\Spec\B\times_\X\Y\to\Spec\B
\]
is an affine cover. As $\Spec\B\to\X$ lifts through $\Y$, we may identify
\[
\Spec\B\times_\X\Y\simeq\Spec\B\times_\Y(\Y\times_\X\Y),
\]
so this follows from \cref{lem:coverssatisfydescent} and the fact that $\Y\times_\X\Y\to\Y$ is an affine cover.
\end{proof}

\subsubsection{Locally descendable morphisms}
Even if $f\colon \sfY\to \sfX$ is affine, the inclusion $\sfD_f\to\sfX$ of its image need not even be $0$-semiaffine. By \cref{prop:qcohdescentstack}, this is essentially due to the failure of nilpotent completion to be a finitary construction. Our next goal is to single out a class of particularly well-behaved universally $0$-affine morphisms for which these issues do not arise.

\begin{mydef}\label{def:ldgeneral}
Let $\cstablecat$ be a symmetric monoidal stable $\infty$-category with unit $\cmonoidalunit$. An object $A \in \cstablecat$ is \emph{locally descendable} if the Bousfield localization $L_A\cmonoidalunit$ lies in $\Thick^\otimes_\cstablecat(A)$. If in addition $\cmonoidalunit$ is $A$-local, then $A$ is simple \emph{descendable}.
\end{mydef}

\begin{remark}
Descendability was first studied by Balmer \cite[Df.3.16]{balmer2016separable} (there called ``nil-faithfulness'') and Mathew \cite[Df.3.18]{akhilgalois}, see \cite{mathew2018examples} for a detailed survey. The extension to local descendability was studied in \cite[\textsection A.1]{balderrama2024total} with a slightly different definition, which is seen to be equivalent to \cref{def:ldgeneral} by \cref{lem:ldequiv} below.
\end{remark}

\begin{mydef}\label{def:ld}
A morphism $f\colon \sfY\to\sfX$ of stacks is \emph{(locally) descendable} if $f$ is universally $0$-affine and $f_\ast\calO_\Y\in\QCoh(\X)$ is (locally) descendable.
\end{mydef}

The main observation of this section is the following.

\begin{theorem}\label{thm:locallydescendable}
Suppose given a (locally) descendable morphism $f\colon \Y\to\X$. Then the inclusion $j\colon \sfD_f\to\sfX$ is (locally) descendable. Moreover, there are equivalences
\[
j_\ast\calO_{\sfD_f}\simeq L_{f_\ast\calO_\sfY},\qquad \QCoh(\sfD_f)\simeq L_{f_\ast\calO_\sfY}(\QCoh(\sfX)).
\]
\end{theorem}

The proof of \cref{thm:locallydescendable} requires several preliminaries. We begin by recalling some of the general theory of descendability. Let $\cstablecat$ denote a presentably symmetric monoidal stable $\infty$-category with unit $\cmonoidalunit$. Say that $A \in \cstablecat$ is an \emph{algebra} if it is equipped with a unital product. Given an algebra $A \in \cstablecat$, define the tower $T(A)$ by
\begin{equation}\label{eq:adamstower}
T(A)_n = \cof\left(\ol{A}^{\otimes n} \to \cmonoidalunit\right)\quad\text{where}\quad \ol{A} = \fib\left(\monoidalunit_\cstablecat \to A\right).
\end{equation}
This is the tower under $\cmonoidalunit$ dual to the Adams tower of $A$ \cite[\textsection5]{bousfieldlocalisationexists}. If $A$ admits a coherently associative multiplication, then it is equivalent to the tower of partial totalizations of the cosimplicial object $A^{\otimes \bullet+1}$ \cite[Pr.2.14]{nilpotenceanddescent}.

\begin{lemma}\label{lem:ldequiv}
Fix an algebra $A \in \cstablecat$, and let $L_A$ denote the functor of $A$-localization on $\cstablecat$. The following are equivalent:
\begin{enumerate}
\item $A$ is locally descendable, i.e.\ $L_A\cmonoidalunit \in \Thick^\otimes(A) = \Thick(\{A\otimes X : X \in \cstablecat\})$.
\item The tower $T(A)$ is pro-constant.
\end{enumerate}
Moreover, in this case, $A$-localization is smashing and agrees with $A$-nilpotent completion, i.e.\
\[
L_A X \simeq L_A\cmonoidalunit \otimes X\simeq \lim (T(A) \otimes X) 
\]
for any $X \in \cstablecat$.
\end{lemma}
\begin{proof}
We first claim that, in general, if $M \in \Thick^\otimes(A)$ then the canonical maps
\[
M \rightarrow L_A M \leftarrow L_A\cmonoidalunit \otimes M
\]
are equivalences. As $\Thick^\otimes(A) = \Thick(\{A\otimes X : X \in \cstablecat\})$ and these are natural transformations between exact functors in $M$, it suffices to consider just the case where $M = A \otimes X$. Now $L_A(A\otimes X)\simeq A\otimes X$ because $A$ is an algebra, and $L_A\cmonoidalunit \otimes A \otimes X\simeq A\otimes X$ because $L_A\cmonoidalunit \otimes A \simeq \cmonoidalunit \otimes A \simeq A$ by definition of $A$-localization.

(1)$\Rightarrow$(2): As each object in the tower $T(A)$ lives in $\Thick^\otimes(A)$, we may identify
\[
T(A)\simeq L_A\cmonoidalunit \otimes T(A).
\]
As $L_A\cmonoidalunit \in \Thick^\otimes(A) = \Thick(\{A\otimes X\} : X \in \cstablecat)$ and the family of pro-constant towers form a thick subcategory of all towers, it follows that $T(A)$ is pro-constant provided that $A \otimes T(A)$ is pro-constant, which holds for any $A$.

(2)$\Rightarrow$(1): Suppose that $T(A)$ is pro-constant. We first prove the extra claim that $A$-localization is smashing and agrees with $A$-nilpotent completion.

As $T(A)$ is pro-constant, we may identify
\[
X_A^\wedge = \lim (T(A) \otimes X) \simeq (\lim T(A))\otimes X
\]
for any $X \in \cstablecat$, so it suffices to prove that $X \to X_A^\wedge$ exhibits $X_A^\wedge$ as the $A$-localization of $X$. As $X_A^\wedge$ is $A$-local, being a limit of $A$-modules, we must verify that $X \to X_A^\wedge$ is an $A$-equivalence. This holds as $A \otimes T(A)$ is pro-equivalent to $A$ for any $A$, and therefore $A \otimes X_A^\wedge\simeq \lim A \otimes T(A)\otimes X \simeq A \otimes X$.

This proves that if $T(A)$ is pro-constant, then it is pro-equivalent to the constant tower on $L_A\cmonoidalunit$. This implies that $L_A\cmonoidalunit \to T(A)_n$ admits a retraction for some $n$. As each $T(A)_n \in \Thick^\otimes(A)$, it follows that $L_A\cmonoidalunit \in \Thick^\otimes(A)$ as claimed.
\end{proof}

Say that a map $A \to B$ in $\CAlg(\cstablecat)$ is (locally) descendable if it realizes $B$ as a (locally) descendable object of $\Mod_A(\cstablecat)$. Given $A \in \CAlg(\cstablecat)$ and $E \in \Mod_A(\cstablecat)$, write $L_{E/A}$ for the functor of $E$-localization inside $\Mod_A(\cstablecat)$. If $E$ is an algebra in $\Mod_A$, then we write $T_A(E)$ for the tower of (\ref{eq:adamstower}) formed in $\Mod_A(\cstablecat)$. This satisfies the base change formula $B\otimes_A T_A(E) = T_B(B\otimes_A E)$ for any map $A \to B$ in $\CAlg(\cstablecat)$.

\begin{lemma}\label{lem:ldclosure}
Local descendability has the following closure properties.
\begin{enumerate}
\item If $F\colon \cstablecat\to\dstablecat$ is an exact symmetric monoidal functor between presentably symmetric monoidal stable categories and $A \in \cstablecat$ is a (locally) descendable algebra, then $F(A)$ is (locally) descendable and $F(L_A X)\simeq L_{F(A)}F(X)$ for all $X \in \cstablecat$.
\item Consider a pushout square 
\begin{center}\begin{tikzcd}
A\ar[r]\ar[d]&B\ar[d]\\
C\ar[r]&D
\end{tikzcd}\end{center}
in $\CAlg(\cstablecat)$.  If $A \to C$ is (locally) descendable, then $B \to D$ is (locally) descendable and $L_{D/B}(B\otimes_A M) = B\otimes_A L_{C/A}M$ for any $M \in \Mod_A(\cstablecat)$.
\item In a pushout square as above, suppose that $A \to B$ is locally descendable and $L_{B/A}C = C$. Then $A\to C$ is locally descendable if and only if $B \to D$ is locally descendable.
\item In a pushout square as above, suppose that $C \in \Thick^\otimes_A(B)$. Then $L_{B/A}C = C$ and $C \to D$ is descendable. This holds, for example, if $A \to B$ is descendable or $A \to C$ factors through $B$.
\item Consider a composable pair $f\colon A \to B$ and $g\colon B \to C$. If $f$ is locally descendable and $L_{C/B}(B\otimes_A C) = B\otimes_A C$, then $g\circ f$ is locally descendable. In particular, this holds if $g$ is descendable.
\end{enumerate}
\end{lemma}
\begin{proof}
(1)~~If $A \in \cstablecat$ is locally descendable, then $T(A)$ is pro-constant. As pro-constant towers are preserved by any functor, it follows that $F(T(A))$ is pro-constant. As $F$ is exact and symmetric monoidal, we have $F(T(A))\simeq T(F(A))$, and therefore $F(A)$ is locally descendable.

If $X \in \sfC$, then as $X \otimes T(A)$ is pro-constant, its limit is preserved by any functor. It follows that $F(L_AX) \simeq F(\lim X\otimes T(A))\simeq \lim F(X\otimes T(A))\simeq \lim F(X) \otimes T(F(A))\simeq L_{F(A)}F(X)$ as claimed. 

Taking $X = \cmonoidalunit$, this implies that $F(L_A\cmonoidalunit) = L_{F(A)}\dmonoidalunit$. Therefore, if $A$ is descendable, then $L_{F(A)}\dmonoidalunit = F(L_A\cmonoidalunit) = F(\cmonoidalunit) = \dmonoidalunit$, implying that $F(A)$ is descendable.

(2)~~Apply (1) to the base change functor $B\otimes_A({-})\colon \Mod_A(\cstablecat) \to \Mod_B(\cstablecat)$.

(3)~~Suppose that $A \to B$ is locally descendable and that $L_{B/A}C = C$. To prove that $A \to C$ is locally descendable, we must show that $T_A(C)$ is pro-constant. As $L_{B/A}C = C$ and $L_{B/A}$-localization is smashing, it follows that $L_{B/A}M = M$ for all $M \in \Thick^\otimes_A(C)$, and so $T_A(C) = L_{B/A}T_A(C) = L_{B/A}A \otimes_A T_A(C)$. As $A \to B$ is locally descendable, $L_{A/B}A$ is in $\Thick_{A}^\otimes(B)$. Since pro-constant towers form a thick subcategory of all towers, it follows that $T_A(C)$ is pro-constant if and only if $B \otimes_A T_A(C) = T_B(B\otimes_A C) = T_B(D)$ is pro-constant, which in turn holds if and only if $B\to D$ is locally descendable as claimed.

(4)~~Suppose that $C \in \Thick^\otimes_A(B)$. Clearly, $L_{B/A}C = C$. To show that $C \to D$ is descendable, we must show that $T_C(D)$ is pro-constant. As before, as $C \in \Thick^\otimes_A(B)$, it suffices to prove that $B\otimes_A T_C(D)$ is pro-constant. Indeed, $B\otimes_A T_C(D) = D\otimes_C T_C(D) = T_D(D\otimes_C D)$ is pro-constant as $D \to D \otimes_C D$ admits a retraction.

(5)~~Consider the diagram
\begin{center}\begin{tikzcd}
A\ar[r]\ar[d]&B\ar[r]\ar[d]&C\ar[d]\\
C\ar[r]&B\otimes_A C \ar[r]&C\otimes_A C
\end{tikzcd}.\end{center}
As $A \to B$ is locally descendable and $A \to C$ factors through $B$, it follows from (3) and (4) that $A \to C$ is locally descendable if and only if $B \to B\otimes_A C$ is locally descendable. As $C \to C\otimes_A C$ is descendable, it follows from (3) that if $L_{C/B}(B\otimes_A C) = B \otimes_A C$ then $B \to B\otimes_A C$ is locally descendable as needed.
\end{proof}

This immediately implies corresponding closure properties for locally descendable morphisms of stacks. We now give the simplest of these.

\begin{prop}\label{prop:ldstackclosure}
Locally descendable morphisms of stacks have the following closure properties.
\begin{enumerate}
\item Suppose given a Cartesian square
\begin{center}\begin{tikzcd}
\sfY'\ar[r,"f'"]\ar[d,"h'"]&\sfX'\ar[d,"h"]\\
\sfY\ar[r,"f"]&\sfX
\end{tikzcd}\end{center}
of stacks. If $f$ is (locally) descendable, then $f'$ is (locally) descendable.
\item Suppose that in the above Cartesian square, $h$ is descendable and $f$ is universally $0$-affine. Then $f$ is (locally) descendable if and only if $f'$ is (locally) descendable.
\item Suppose given morphisms $f\colon \sfY \to \sfX$ and $g\colon \sfZ \to \sfY$ of stacks. If $g$ is descendable and $f$ is universally $0$-affine, then $f$ is (locally) descendable if and only if $g\circ f$ is (locally) descendable.
\end{enumerate}
\end{prop}
\begin{proof}
(1)~~As $f$ is universally $0$-affine, so is $f'$. By \cref{prop:univ0affinefacts}(1), this square is adjointable, meaning that
\[
h^\ast f_\ast \calO_{\sfY} \to f'_\ast h'{ }^\ast \calO_{\sfY} \simeq f'_\ast\calO_{\sfY'} 
\]
is an equivalence. As $f$ is (locally) descendable, $f_\ast\calO_\sfY\in \QCoh(\sfX)$ is (locally) descendable, and hence so is $h^\ast f_\ast\calO_\sfY \in \QCoh(\sfX')$ by \cref{lem:ldclosure}(1). Thus $f'$ is (locally) descendable.

(2)~~As $f$ is universally $0$-affine, we must only prove that $f_\ast\calO_\sfY \in \QCoh(\sfX)$ is (locally) descendable if and only if $f'_\ast\calO_{\sfY'} \in \QCoh(\sfX')$ is (locally) descendable. By \cref{prop:univ0affinefacts}(1) and \cref{prop:adjointablesquare0affine}(3), there is a coCartesian square
\begin{center}\begin{tikzcd}
h_\ast f'_\ast \calO_{\sfY'}&h_\ast\calO_{\sfX'}\ar[l]\\
f_\ast\calO_\sfY\ar[u]&\calO_\sfX\ar[u]\ar[l]
\end{tikzcd}\end{center}
in $\QCoh(\sfX)$. As $h$ is universally $0$-affine, $f'_\ast\calO_{\sfY'}\in\QCoh(\sfX')$ is (locally) descendable if and only if $h_\ast f'_\ast\calO_{\sfY'} \to h_\ast \calO_{\sfX'}$ is (locally) descendable. Therefore, the claim follows from \cref{lem:ldclosure}(3).

(3)~~This follows from \cref{lem:ldclosure} as in the proof of (2).
\end{proof}

We now give the following.

\begin{proof}[Proof of \cref{thm:locallydescendable}]
\cref{prop:qcohdescentstack} and local descendability implies $j_\ast \calO_{\sfD_f}\simeq (\calO_{\sfX})_{j_\ast\calO_{\sfY}}^\wedge\simeq L_{j_\ast\calO_{\sfY}}\calO_{\sfX}$, so we must only prove that $j$ is universally $0$-affine. Fix morphisms $\Spec\sfA'\to \Spec\sfA\to\sfX$. By \cref{lem:descentstackbasicproperties}(\ref{item:restrictdescentstack}), we may form a Cartesian diagram
\begin{center}\begin{tikzcd}
\sfF'\ar[r]\ar[d]&\sfD_{p''}\ar[r]\ar[d]&\Spec\sfA'\ar[d]\\
\sfF\ar[r]\ar[d]&\sfD_{p'}\ar[r]\ar[d]&\Spec\sfA\ar[d]\\
\sfY\ar[r]&\sfD_f\ar[r]&\sfX
\end{tikzcd},\end{center}
where $\sfD_{p'}$ and $\sfD_{p''}$ are the descent stacks of the composites $p'\colon \sfF \to \Spec\A$ and $p''\colon \sfF' \to \Spec\A''$. We must prove that $\sfD_{p'}\to\Spec\A$ is $0$-affine and that the top right square is adjointable. As $f$ is locally descendable, so is $p'$ by \cref{prop:ldstackclosure}(1). If we set $\B = p'_\ast\calO_{\sfF}$, then \cref{prop:qcohdescentstack} implies
\begin{align*}
\QCoh(\sfD_{p'}) &\simeq \lim \left(\begin{tikzcd}[ampersand replacement=\&]
\Mod_\B(\Mod_\A)\ar[r,shift left=1mm]\ar[r,shift right=1mm]\&\Mod_{\B\otimes_\A \B}(\Mod_\A)\ar[r,shift left=1mm]\ar[r]\ar[r,shift right=1mm]\&\cdots\end{tikzcd}\right)\\
&\simeq \lim \left(\begin{tikzcd}[ampersand replacement=\&]
\Mod_\B(\Mod_\A)\ar[r,shift left=1mm]\ar[r,shift right=1mm]\&\Mod_{\B\otimes_{L_\B \A} \B}(\Mod_\A)\ar[r,shift left=1mm]\ar[r]\ar[r,shift right=1mm]\&\cdots\end{tikzcd}\right) \simeq \Mod_{L_\B\A}(\Mod_\A),
\end{align*}
where the final equivalence holds by \cite[Pr.3.22]{akhilgalois} as $L_\B\A \to \B$ is descendable. This proves that $\sfD_{p'} \to \Spec\A$ is $0$-affine. Adjointability of the top right square therefore asks that the canonical map
\[
\A'\otimes_\A L_{\Gamma(\calO_{\sfF})} \A \to L_{\Gamma(\calO_{\sfF'})}\A'
\]
is an equivalence. By the second part of \cref{lem:ldclosure}(2), this holds provided $\A'\otimes_\A p'_\ast(\calO_{\sfF}) \to p''_\ast(\calO_{\sfF'})$ is an equivalence, which follows from adjointability of the top outer rectangle.
\end{proof}

\begin{example}
Let $G$ be a finite group and $f\colon \Spec\sfB \to \Spec\sfA$ be a faithful $G$-Galois extension in the sense of Rognes \cite{rognesgalois}; see \cite[Df.6.12]{akhilgalois} for the generalization to any $\sfC \in \CAlg(\PrLst)$. Assuming that $\tau$ is a Grothendieck topology on $\CAlg(\sfC)^\op$ for which $\Mod_{(-)}$ is a $\tau$-sheaf, the induced morphism $j\colon \Spec\sfB//G \to \Spec\sfA$ is a descendable monomorphism of stacks; where $-//G$ indicates the colimit of the diagram $BG \to \Stk$ defined by the $G$-action on $\sfB$. Indeed, the Galois property allows us to rewrite the \v{C}ech nerve of $f$ as the simplicial resolution of $\Spec \B//G$,
\[
 \colim \left(\begin{tikzcd} {\Spec \B}& {\coprod_G\Spec \B}\ar[l,shift left=1mm]\ar[l,shift right=1mm]&\cdots\ar[l,shift left=1mm]\ar[l]\ar[l,shift right=1mm]\end{tikzcd}\right)
 \simeq \colim \left(\begin{tikzcd} {\Spec \B}& {\Spec \B\otimes_\A \B}\ar[l,shift left=1mm]\ar[l,shift right=1mm]&\cdots\ar[l,shift left=1mm]\ar[l]\ar[l,shift right=1mm]\end{tikzcd}\right),
\]
which identifies $\Spec \B//G$ with $D_f$; this is the same identification as in \cite[\textsection4.1]{osyn}. Again, the fact that $\A \to \B$ is a faithful $G$-Galois extension then implies that $\A \to \B$ is descendable, see \cite[Pr.6.13]{akhilgalois} for example. By \cref{thm:locallydescendable}, we then immediately see that $\Spec \sfB//G \to \Spec \sfA$ is descendable.
\end{example}

\section{Nonconnective spectral algebraic geometry}\label{sec:nonconsag}
We now specialize the previous section to $\csite = \Aff = \CAlg^\op$, the opposite category of $\E_\infty$ rings, equipped with the fpqc topology, and $\sfQ(\A) = \Mod_\A$. The associated category of stacks, which we denote by $\Stk$, is our candidate for a workable category of stacks in \emph{nonconnective spectral algebraic geometry}. We focus here on the connection between $\Stk$, its connective variant $\Stk^\cn$, and its classical variant $\Stk^\heartsuit$, and how this interacts with flatness. We also define quasi-affine morphisms in $\Stk$ and show that such morphisms are universally $0$-affine; see \Cref{ex:quasiaffineareU0A}.

\subsection{Stacks on \texorpdfstring{$\E_\infty$}{E-infinity} rings with the fpqc topology}\label{ssec:stackswithfpqctopology}
\subsubsection{Basic definitions, connective covers, and truncations}\label{sssec:basicdefinitions}
Recall that a map of $\E_\infty$ rings $\A\to \B$ is said to be \emph{(faithfully) flat} if the underlying map on $\pi_0$ is (faithfully) flat and the natural map $\pi_\ast\A \otimes_{\pi_0 \A}\pi_0 \B \to \pi_\ast\B$ is an isomorphism; see \cite[Df.7.2.2.10]{ha}. One can now define the fpqc topology on $\Aff = \CAlg^\op$ following \cite[Rmk.D.6.1.1]{sag}. 

\begin{mydef}\label{def:fpcq_top}
A collection of maps of $\E_\infty$ rings $\{\A \to \A_i\}$ is an \emph{fpqc covering} for $\A$ if it contains a finite collection of maps $\{\A \to \A_j\}$ such that $\A \to \prod \A_j$ is faithfully flat. This generates the \emph{fpqc topology} on $\Aff = \CAlg^\op$. Similarly, we have fpqc topologies on $\Aff^\cn = (\CAlg^{\cn})^\op$ and $\Aff^\heartsuit = (\CAlg^\heartsuit)^\op$, the categories of connective and discrete $\mathbb{E}_\infty$ rings.
\end{mydef}

\begin{lemma}
    The pairs $(\Aff,\fpqc)$, $(\Aff^\cn,\fpqc)$, and $(\Aff^\heartsuit,\fpqc)$ are locally small sites. 
\end{lemma}

\begin{proof}
    The facts that these categories are all locally small and admit finite limits are clear. The fact that $\fpqc$ covers define a Grothendieck pretopology on these categories is classical for $\Aff^\heartsuit$, and follows by standard collapsing K\"unneth spectral sequence arguments for $\Aff^\cn$ and $\Aff$.
\end{proof}

This allows us to define our categories of stacks in spectral algebraic geometry.

\begin{mydef}\label{def:categoryoffpqcstacks}
    Define the categories of \emph{stacks}, \emph{connective stacks}, and \emph{classical stacks} as
    \[\Stk = \Stk_\fpqc(\Aff), \qquad \Stk^\cn = \Stk_\fpqc(\Aff^\cn), \qquad \Stk^\heartsuit = \Stk_\fpqc(\Aff^\heartsuit).\]
\end{mydef}

\begin{cor}
    There are natural equivalences of categories
    \[\Stk \simeq \colim_\kappa \Shv_\fpqc(\Aff_\kappa), \qquad \Stk^\cn \simeq \colim_\kappa \Shv_\fpqc(\Aff^\cn_\kappa), \qquad \Stk^\heartsuit \simeq \colim_\kappa \Shv_\fpqc(\Aff_\kappa^\heartsuit).\]
\end{cor}

\begin{proof}
    These equivalences follow from \Cref{cor:accessiblesites} as the categories $\csite^\op = \CAlg$, $\CAlg^\cn$, and $\CAlg^\heartsuit$ are all accessible, and $(\csite_\kappa,\fpqc)$ are subsites of $(\csite, \fpqc)$ for all regular cardinals $\kappa$.
\end{proof}

There are natural adjoint functors between the categories underlying these sites, given by either the natural inclusions of subcategories, the connective cover functor $\tau_{\geq 0}$, or the truncation functor $\tau_{\leq 0}$:
\begin{equation}\label{eq:adjunctionsoncategories}\begin{tikzcd}
    {\CAlg^\heartsuit} & {\CAlg^\cn} & {\CAlg.}
	\arrow[shift right=2, from=1-1, to=1-2]
	\arrow[shift right=2, from=1-2, to=1-1, "{\tau_{\leq 0}}", swap]
	\arrow[shift left=2, from=1-3, to=1-2, "{\tau_{\geq 0}}"]
	\arrow[shift left=2, from=1-2, to=1-3]
\end{tikzcd}\end{equation}

These adjunctions, where left adjoints are always written on top, induce functors between the associated categories of stacks.

\begin{prop}\label{th:coversandtruncations}
    The adjunctions of (\ref{eq:adjunctionsoncategories}) induces the adjunctions between categories of stacks
\[\begin{tikzcd}
    {\Stk^\heartsuit} & {\Stk^\cn} & {\Stk,}
	\arrow[shift left=2, from=1-2, to=1-1, "{\tau_{\leq 0}}"]
	\arrow[shift left=2, from=1-1, to=1-2]
	\arrow[shift right=2, from=1-2, to=1-3]
	\arrow[shift right=2, from=1-3, to=1-2, "{\tau_{\geq 0}}", swap]
\end{tikzcd}\]
    where the unlabelled functors are fully faithful.
\end{prop}

The functor $\tau_{\geq0}\colon \Stk \to \Stk^\cn$ is simply notation for $(\tau_{\geq0})_\ast$ of \Cref{prop:accessiblemorphismofsites}

\begin{proof}
    First, notice that all of the functors of (\ref{eq:adjunctionsoncategories}) are morphisms of sites by a standard collapsing K\"unneth spectral sequence argument. The (opposites of the) adjunctions of (\ref{eq:adjunctionsoncategories}) then induce the desired adjunctions and fully faithfulness between categories of stacks by \Cref{cor:adjunctionsandmapsofsites}.
\end{proof}

\begin{mydef}
    Define the \emph{underlying stack functor} $(-)^\heartsuit\colon \Stk \to \Stk^\heartsuit$ as the composite
    \[\Stk \xrightarrow{\tau_{\geq 0}} \Stk^\cn \xrightarrow{\tau_{\leq0}} \Stk^\heartsuit.\]
\end{mydef}

It follows from \Cref{th:coversandtruncations} that connective covers and underlying stack functors preserve presentations of stacks in terms of colimits of representables.

\begin{cor}\label{cor:coversandunderlyingascolimits}
    Let $\X$ be a stack written as a colimit of representables $\X \simeq \colim \Spec \A$. Then the natural maps of connective stacks and classical stacks
    \[\colim \Spec \tau_{\geq 0} \A \xrightarrow{\simeq} \tau_{\geq 0} \X,\qquad \colim \Spec \pi_0 \A \xrightarrow{\simeq} \X^\heartsuit,\]
    respectively, are equivalences. In particular, all of the functors of \Cref{th:coversandtruncations} preserve affine objects.
\end{cor}

\begin{proof}
    Both $\tau_{\geq 0}\colon \Stk \to \Stk^\cn$ and $\tau_{\leq0} \colon \Stk^\cn \to \Stk^\heartsuit$, as well as their composite, preserve colimits.
\end{proof}

\begin{remark}
The above corollary is one of the only ways in which we can compute the functors $\tau_{\geq0}\colon \Stk \to \Stk^\cn$ and $\tau_{\leq0}\colon \Stk^\cn \to \Stk^\heartsuit$. One should also not confuse these functors with restriction along the inclusions $i\colon \Aff^\cn \to \Aff$ or $j\colon \Aff^\heartsuit \to \Aff^\cn$. For example, the composite $(ij)^\ast\colon \Stk \to \Stk^\heartsuit$ applied to $\M_\FG^\ori$ is represented by the zero ring by \cite[Rmk.4.3.11]{ec2}, as the only discrete ring $\A$ which admits a map $\Spec\A\to\M_\FG^\ori$ is the zero ring, whereas the underlying stack construction instead satisfies $(\M_\FG^\ori)^\heartsuit \simeq \M_\FG^\heartsuit$; see \Cref{underlyingclassicalstackofDMUP} below.
\end{remark}

\subsubsection{The canonical quasi-coherent sheaf context}
There is an obvious candidate for a quasi-coherent sheaf context on $\Stk$. Let
\[\Mod_{(-)} \colon \CAlg \to \CAlg(\Pr^L_\st)\]
be the fully faithful functor sending $\A$ to $\Mod_\A$ of \cite[Pr.7.1.2.7]{ha}.

\begin{lemma}\label{lm:canonicalqcohcontext}
    The triple $(\Aff,\fpqc, \Mod_{(-)})$ forms a tautological quasi-coherent sheaf context.
\end{lemma}

\begin{proof}
    By \Cref{ex:basicqcohcontext}, it suffices to show that $\Mod_{(-)}$ is an fpqc sheaf, which follows from \cite[Cor.D.6.3.3]{sag}.
\end{proof}

\begin{mydef}\label{def:thecanonicalqcohcontext}
    Write $\QCoh(-)\colon \Stk^\op \to \CAlg(\PrLst)$ for the quasi-coherent sheaf functor associated to the context of \Cref{lm:canonicalqcohcontext}.
\end{mydef}

By construction, we have $\QCoh(\Spec \A) \simeq \Mod_\A$ and if $\X \simeq \colim \Spec \A$, then $\QCoh(\X) \simeq \lim \Mod_\A$.

We will see in \Cref{ex:qcohofspDM} that this notion of quasi-coherent sheaf agrees with Lurie's for nonconnective spectral Deligne--Mumford stacks \cite[Df.2.2.2.1]{sag}.

\begin{remark}
    Restricting the functor $\Mod_{(-)} \colon \CAlg \to \CAlg(\PrLst)$ to either $\CAlg^\cn$ or $\CAlg^\heartsuit$ also produce quasi-coherent sheaf contexts for the $\fpqc$ topology. However, there is no need to consider these contexts separately, as the associated quasi-coherent sheaf functors on $\Stk^\cn$ and $\Stk^\heartsuit$ are the restriction of the above $\QCoh(-)$ along the fully faithful embeddings $\Stk^\heartsuit \to \Stk^\cn \to \Stk$. This follows from the fact that these embeddings preserve affine objects and colimits.
\end{remark}

\subsection{Nonconnective spectral Deligne--Mumford stacks}\label{ssec:spdm}

We now describe how our category $\Stk$ of fpqc stacks relates to the category $\SpDMnc$ of nonconnective spectral Deligne--Mumford stacks (\cite[Df.1.4.4.2]{sag}) developed extensively by Lurie in \cite{sag}.

\subsubsection{The embedding into étale stacks}
Recall that a map of $\E_\infty$ rings $\A \to \B$ is \emph{étale} if it is flat and the associated map on $\pi_0$ is étale. This then leads to an étale topology on $\Aff$. The setup of \cref{ssec:stackssection} defines a category
\[
\Stk^\et = \Stk_\et(\Aff)
\]
of \emph{étale stacks}, whose inclusion of affines we write as
\[
\Spet\colon \CAlg^\op\to\Stk^\et.
\]
Our first goal of this subsection is to show that $\SpDMnc$ embeds nicely into the category $\Stk^\et$ of \'etale stacks. This is a minor variation of the fact, proved in \cite{rezk2022spectral}, that the functor of points embedding $\SpDM^{\nc}\to\Fun(\CAlg,\spaces)$ is fully faithful. 

Given a regular cardinal $\kappa$, write $\Shv^\et_\kappa = \Shv(\Aff_\kappa,\text{\'etale})$, where $\Aff_\kappa$ denotes the opposite category of $\kappa$-compact $\E_\infty$ rings. In particular, we have $\Stk^\et\simeq\colim_\kappa\Shv^\et_\kappa$ by \Cref{cor:accessiblesites}.

\begin{lemma}
Fix regular cardinals $\kappa < \lambda$. Then in the geometric morphism
\[
i_* : \Shv_\kappa^\et \rightleftarrows \Shv_\lambda^\et : i^\ast
\]
associated to the inclusion $i\colon \Aff_\kappa\subseteq\Aff_\lambda$, the left adjoint $i_*$ is fully faithful and the right adjoint $i^\ast$ preserves colimits.
\end{lemma}
\begin{proof}
To show that $i_*$ is fully faithful, it suffices to show that the unit $\sfX \to i^\ast i_* \sfX$ is an equivalence for all $\sfX \in \Shv_\kappa^\et$. This is clearly an equivalence if $\sfX \in \Aff_\kappa$ as $\Aff_\kappa$ generates $\Shv_\kappa^\et$ under colimits, so it suffices to prove that $i^\ast$ preserves colimits.

For an $\E_\infty$ ring $\A$, write $\CAlg_A^\et \subseteq\CAlg_\A$ for the category of \'etale $\A$-algebras and
\[
\Shv^\et_\A\subseteq\Fun(\CAlg_\A^\et,\spaces)
\]
for the small \'etale topos of $\A$. As in \cite[\textsection15.6]{rezk2022spectral}, $\Shv_\kappa^\et$ is equivalent to the category of sections of the Cartesian fibration associated to
\[
\CAlg_\kappa^\op\to\Cat_\infty,\qquad \A \mapsto \Shv_\A^\et.
\]
In particular, colimits in $\Shv_\kappa^\et$ are computed levelwise as colimits in each $\Shv_\A^\et$. This description is independent of $\kappa$, and therefore $i^\ast$ preserves colimits.
\end{proof}

Let $\SpDMncet$ denote the wide subcategory of $\SpDMnc$ spanned by étale morphisms (\cite[Df.1.4.10.1]{sag}), and $\SpDM^{\nc,\kappa}\subseteq\SpDM^{\nc}$ denote the full subcategory generated under colimits in $\SpDMncet$ by $\Spet \A$ for $\kappa$-compact $\E_\infty$ rings $\A$.

\begin{lemma}\label{lem:kappadmstacksintoetalestacks}
The functor of points embedding
\[
h_{(\bs)}^{\et,\kappa}\colon \SpDM^{\nc,\kappa} \to \Shv^\et_\kappa\subseteq\Fun(\CAlg_\kappa,\spaces),\qquad \sfX \mapsto \Map_{\SpDM^{\nc,\kappa}}(\Spet(\bs),\sfX)
\]
is fully faithful, and preserves colimits when restricted to the wide subcategory $\SpDM^{\nc,\kappa}_\et\subseteq\SpDM^{\nc,\kappa}$ of \'etale morphisms. Moreover, if $\kappa < \lambda$ then the diagram 
\begin{center}\begin{tikzcd}
\SpDM^{\nc,\kappa}\ar[r,tail]\ar[d,tail,"h_{(\bs)}^{\et,\kappa}"]&\SpDM^{\nc,\lambda}\ar[d,"h_{(\bs)}^{\et,\lambda}",tail]\\
\Shv_\kappa^\et\ar[r,"j_\ast",tail]&\Shv_\lambda^\et
\end{tikzcd}\end{center}
commutes.
\end{lemma}
\begin{proof}
With $\kappa$ omitted, the first statement is proved in \cite[\textsection15.6]{rezk2022spectral}. The same proof applies for each fixed $\kappa$, so we just verify that the diagram commutes. If $\sfX \in \SpDM^{\nc,\kappa}$ then $j^\ast h_\sfX^{\et,\lambda}\simeq h_\sfX^{\et,\kappa}$, so there is a natural transformation
\[
j_\ast h_{(\bs)}^{\et,\kappa} \to h_{(\bs)}^{\et,\lambda}.
\]
As both sides preserves colimits when restricted to $\SpDM^{\nc,\kappa}_\et\subseteq\SpDM^{\nc,\kappa}$ to show that this is a natural isomorphism, it suffices to show that it is a natural isomorphism on restriction to $\Aff_\kappa\subseteq\SpDM^{\nc,\kappa}$, which is clear. 
\end{proof}

This leads us to our first main result of this subsection.

\begin{theorem}\label{thm:dmstacksintoetalestacks}
There is a unique fully faithful embedding
\[
h_{(\bs)}^\et\colon \SpDMnc\to\Stk^\et
\]
extending the inclusion of affines, i.e.\ satisfying
\[
h_{\Spet\A}\simeq \Spet \A
\]
naturally in $\A\in\CAlg$. Moreover, $h_{(\bs)}^\et$ preserves colimits when restricted to the wide subcategory $\SpDMncet\subseteq\SpDM^{\nc}$ spanned by the \'etale morphisms.
\end{theorem}
\begin{proof}
The functor $h_{(\bs)}^{\et}$ is obtained from \cref{lem:kappadmstacksintoetalestacks} by taking a colimit along $\kappa$. Uniqueness is a consequence of the fact that the affines $\Aff\subseteq\Stk^\et$ form a dense subcategory by \cref{prop:limitextension}.
\end{proof}

\subsubsection{The functor into fpqc stacks}

We now use this discussion to relate $\SpDMnc$ to our category of stacks. As \'etale covers are fpqc covers, \cref{prop:accessiblemorphismofsites} applied to the identity on $\Aff$ provides an fpqc sheafification functor
\[
L^{\fpqc}\colon \Stk^{\et} \to \Stk.
\]
The composite functor
\[
h_{(\bs)} = L^{\fpqc}\circ h_{(\bs)}^\et\colon \SpDMnc\to\Stk^\et\to\Stk.
\]
is generally \emph{not} fully faithful.

\begin{mydef}\label{df:flatdescentDMstacks}
    A nonconnective spectral Deligne--Mumford stack $\sfX \in \SpDMnc$ \emph{satisfies flat descent} if its functor of points
\[
\Map_{\SpDMnc}(\Spet(\bs),\sfX)\colon \CAlg\to\spaces
\]
satisfies flat descent. Write
\[
\SpDM^{\nc,\flat}\subseteq\SpDM^{\nc}
\]
for the full subcategory of nonconnective spectral Deligne--Mumford stacks that satisfy flat descent.
\end{mydef} 

\begin{prop}\label{prop:dmstacksintoflatstacks}
The functor
\[
h_{(\bs)}\colon \SpDMnc\to\Stk
\]
has the following properties:
\begin{enumerate}
\item There are natural equivalences $h_{\Spet \A}\simeq \Spec \A$ for $\A\in\CAlg$;
\item $h_{(\bs)}$ preserves colimits when restricted $\SpDMncet$;
\item If $\sfX$ satisfies flat descent, then
\[
\Map_{\SpDMnc}(\sfY,\sfX) \to \Map_{\Stk}(h_\sfY,h_\sfX)
\]
is an equivalence for $\sfY \in \SpDMnc$. In particular, $h_{(\bs)}$ is fully faithful when restricted to $\SpDM^{\nc,\flat}$.
\item Any \'etale sheaf on $\SpDM^\nc$ whose restriction to $\Aff$ satisfies flat descent factors uniquely through $h_{(\bs)}$.
\end{enumerate}
\end{prop}
\begin{proof}
(1)~~Immediate from the definitions.

(2)~~This follows from \cref{thm:dmstacksintoetalestacks}(2) as $L^{\fpqc}$ preserves all colimits.

(3)~~By writing $\sfY$ as a colimit of affines along \'etale maps, we reduce to the case where $\sfY = \Spet \B$ for some $\B\in\CAlg$. Choose a regular cardinal $\kappa$ for which $\B$ is $\kappa$-compact and $\sfX \in \SpDM^{\nc,\kappa}$. For each $\lambda \geq \kappa$, consider the square
\begin{center}\begin{tikzcd}
\Shv_\kappa^\et\ar[r,"j_\ast"]\ar[d,"L^{\fpqc}"]&\Shv_\lambda^\et\ar[d,"L^{\fpqc}"]\\
\Shv_\kappa^\fpqc\ar[r,"j_\ast"]&\Shv_\lambda^{\fpqc}
\end{tikzcd}.\end{center}
By \cref{lem:kappadmstacksintoetalestacks}, we have $j_\ast h_\sfX^{\et,\kappa}\simeq h_\sfX^{\et,\lambda}$. As $\sfX$ satisfies flat descent, this is already an fpqc sheaf. It follows that
\begin{align*}
\Map_{\Stk}(\Spec \B,h_\sfX)&\simeq\colim_{\lambda>\kappa}\Map_{\Shv_\lambda^\fpqc}(\Spec\B,L^{\fpqc} h_\sfX^{\et,\lambda})\\
&\simeq \colim_{\lambda > \kappa}\Map_{\presheaves(\Aff_\lambda)}(\Spec\B,h_\sfX^{\et,\lambda})\simeq \Map_{\SpDM^{\nc}}(\Spet\B,\sfX)
\end{align*}
as claimed.

(4)~~As $F$ is an \'etale sheaf, it extends uniquely to a limit-preserving presheaf on $\Stk^{\et,\op}$. This factors (necessarily uniquely) through $L^{\fpqc}\colon \Stk^{\et}\to\Stk$ if and only if $F$ satisfies flat descent.
\end{proof}

\begin{example}\label{ex:qcohofspDM}
The functor
\[
\QCoh\colon \SpDM^{\nc,\op} \to \CAlg(\PrLst)
\]
sending a nonconnective spectral Deligne--Mumford stack to its category of quasi-coherent sheaves satisfies flat descent, and therefore factors uniquely through
\[
h_{(\bs)}\colon \SpDMnc\to \Stk.
\]
In other words, for any $\sfX \in \SpDMnc$ there is an equivalence
\[
\QCoh(\sfX)\simeq\QCoh(h_\sfX).
\]
Thus, if one is interested in studying properties of categories of quasi-coherent sheaves, then no information is lost in working with $\Stk$ instead of $\SpDMnc$.
\end{example}

The nonconnective spectral Deligne--Mumford stacks that arise in practice frequently satisfy flat descent. For example \cite[Th.1.6.2.1]{ec3} implies that the nonconnective spectral Deligne--Mumford stack associated to a nonconnective spectral scheme always satisfies flat hyperdescent. Non-schematic examples are obtained from the following.

\begin{prop}\label{prop:qaffdiagonalmeansdescent}
Let $\sfX$ be a nonconnective spectral Deligne--Mumford stack. If $\sfX$ has quasi-affine diagonal, meaning that for each map $\Spet \A \to \sfX\times \sfX$, the pullback along the diagonal map $\sfX \to \sfX\times \sfX$ is quasi-affine \'{a} la \cite[Df.2.4.1.1]{sag}, then $\sfX$ satisfies flat descent.
\end{prop}

In particular, nonconnective spectral schemes (\cite[Df.1.1.2.8]{sag}) and algebraic spaces (\cite[Df.1.6.8.1]{sag}) all satisfy flat descent.

Our argument closely follows Lurie's proof of \cite[Pr.9.1.4.3]{sag}, which unfortunately is only stated in the \emph{connective} case. In fact, our proof will show that $\sfX$ satisfies flat \emph{hyperdescent}.

\begin{proof}
    The inclusion of the full subcategory $\widehat{\Shv}_\fpqc$ of $\Fun(\CAlg, \widehat{\Spc})$ of large fpqc sheaves into large presheaves admits a left exact left adjoint $L$ after possibly enlarging the universe. Let us write $X = h_\sfX$ and $\al\colon X\to Y=L(X)$ for the natural unit map. First, we want to show that the value of this unit map $X(\A) \to Y(\A)$ on an $\E_\infty$ ring $\A$ is $(-1)$-truncated. To see this, it suffices to show that for each pair $x,x'\in X(\A)$, the natural map of (large) spaces
    \[\be\colon \{x\} \times_{X(\A)} \{x'\} \to \{\eta(\A)(x)\} \times_{Y(\A)} \{\eta(\A)(x')\}\]
    is an equivalence. By the Yoneda lemma, the points $x,x'$ both represent maps of nonconnective spectral Deligne--Mumford stacks from $\Spet \A$ to $\sfX$. Let $\sfX' = \Spet \A \times_\sfX \Spet \A$ be the fiber product over these two maps. To see that $\be$ is an equivalence, it then suffices to show that the natural map of (large) presheaves $\ga\colon h_{\sfX'} \to \Spec \A \times_Y \Spec \A$ is an equivalence. The fact that $L$ is left exact shows that $\ga$ induces an equivalence $L(h_{\sfX'}) \simeq \Spec \A \times_Y \Spec \A$ of fpqc-hypersheaves. It then suffices to see that $h_{\sfX'}$ itself is an fpqc-hypersheaf, which is precisely the content of \cite[Pr.2.4.3.1]{sag}, as by assumption $\sfX'$ is quasi-affine.

    To see that the map $\al$ is actually an equivalence, it now suffices to show that $\al$ is an effective epimorphism in the \'{e}tale topology on $\CAlg$, as both sides are \'{e}tale sheaves. In other words, for each point $\eta\in Y(\A)$, we want to show the existence of an \'{e}tale faithfully flat map of $\E_\infty$ rings $\A \to \B$ such that the image of $\eta$ in $Y(\B)$, denoted by $\eta_\B$, lies in the essential image of $X(\B) \to Y(\B)$. From the definition $Y=LX$, we see that associated with such an $\eta\in Y(\A)$ is a finite collection of flat maps $\A \to \B_i$ such that the product $\A \to \prod \B_i$ is faithfully flat and that each $\eta_{\B_i}$ lies in the image of $X(\B_i) \to Y(\B_i)$. Writing $\B^0 = \prod \B_i$ and $\B^\bullet$ for the \v{C}ech nerve of $\A \to \B^0$, we consequently see that each $\eta_{\B^n}$ lies in the essential image of the fully faithful embedding $X(\B^n) \to Y(\B^n)$. In particular, there is a unique choice of lift of each $\eta_{\B^n}$ to $X(\B^n)$, and this lift is classified by a map $\Spet \B^n \to X$, defining a simplical object over $X$. Let $\Spet \sfC \to \X$ be an \'{e}tale map from an affine. For each $n\geq 0$, let $\sfV_n$ denote the pullback $\sfU\times_\sfX \Spet \B^n$. The fact that $\sfX$ has quasi-affine diagonal means that $\sfV_n$ is itself a quasi-affine nonconnective spectral Deligne--Mumford stack over $\Spet \B^n$. As the category of quasi-affine nonconnective spectral Deligne--Mumford stacks satisfies fpqc-descent by \cite[Pr.2.4.3.2]{sag}, there exists a unique quasi-affine $\sfV$ over $\Spet \A$ such that $\sfV_\bullet \simeq \sfV \times_{\Spet \A} \Spet \B^\bullet$. Choosing an \'{e}tale surjection $\Spet \A' \to \sfV$ to the quasi-affine nonconnective spectral scheme $\sfV$, we now consider the commutative diagram of nonconnective spectral Deligne--Mumford stacks
    \[\begin{tikzcd}
        {\Spet \B'}\ar[r, "p'"]\ar[d, "g'"] &   {\sfV_0}\ar[r, "f'"]\ar[d]    &   {\Spet \B^0}\ar[d, "g"]  \\
        {\Spet \A'}\ar[r, "p"]  &   {\sfV}\ar[r, "f"]       &   {\Spet \A,}
    \end{tikzcd}\]
    the right-hand square is a pullback by the definition of $\sfV$ and the left-hand square is defined as a pullback. By construction, $f'$ is \'{e}tale, and $p'$ is \'{e}tale by base change, so we see that the upper-composite is \'{e}tale. Moreover, the map $g$ is faithfully flat by assumption, hence $g'$ is also faithfully flat. To show that $f$ \'{e}tale, it suffices to show that its precomposition with the \'{e}tale cover $p$ is \'{e}tale, as \'{e}tale morphisms are \'{e}tale local on the target by \cite[Ex.2.8.1.8]{sag}. We are now reduced to show that the map of $\E_\infty$ rings $\A \to \A'$ is \'{e}tale. A collapsing K\"unneth spectral sequence shows that $\pi_0 \B'\simeq \pi_0 \B^0 \otimes_{\pi_0 \A} \pi_0 \A'$, so we see that $\pi_0 \A \to \pi_0 \A'$ is \'{e}tale by classical fpqc-descent for \'{e}tale morphisms of commutative rings; see \cite[\href{https://stacks.math.columbia.edu/tag/0BTL}{0BTL}]{stacks}. It remains to show that $\pi_\ast \A\otimes_{\pi_0\A} \pi_0 \A' \to \pi_\ast \A'$ is an isomorphism of $\pi_0\A'$-modules. As $\pi_0\A'\to \pi_0 \B'$ is fully faithful, it suffices to show this map is an equivalence after base change, in which case we see that
    \[\pi_\ast \A \otimes_{\pi_0 \A} \pi_0 \A' \otimes_{\pi_0 \A'} \pi_0 \B' \simeq \pi_\ast \A \otimes_{\pi_0\A} \pi_0 \B' \xrightarrow{\simeq} \pi_\ast \B' \simeq \pi_\ast \A' \otimes_{\pi_0\A'} \pi_0 \B'\]
    is an equivalence, the second isomorphism courtesy of the fact that $\A \to \B^0 \to \B'$ is flat, and the third from the expression for $\pi_0\B'$ above.
\end{proof}

\begin{remark}\label{rmk:fpqcclassicalcheck}
    The above proposition and its proof apply just as easily in the connective and classical settings. Indeed, \cite[Pr.9.1.4.3]{sag} is precisely stated in the connective setting and the classical statement is a formal consequence.
\end{remark}

\subsection{Quasi-affine morphisms between stacks}\label{ssec:quasiaffine}
Quasi-affine morphisms give us many examples of $0$-affine stacks, and play a key role in distinguishing between reconstruction and affineness.

The goal of this subsection is to define quasi-affine morphisms and show that they provide an example of universally $0$-affine morphisms; see \Cref{ex:quasiaffineareU0A}.

\subsubsection{Basic definitions and properties}

\begin{mydef}\label{def:openimmersion}
    If $\Spec \A$ is an affine and $U\subseteq |\Spec \A|$ an open subset of its underlying space, then we write $(\Spec \A)|_U$ for the (nonconnective) spectral scheme defined by restriction to $U$. A morphism of stacks $f\colon \X \to \Spec \A$ is an \emph{qc-open immersion} if for some quasi-compact open subset $U\subseteq |\Spec \A|$, the map $f$ factors as $\X \to (\Spec \A)|_U \to \Spec \A$ and the first map is an equivalence.
\end{mydef}

\begin{mydef}
    A stack $\X$ is \emph{quasi-affine} if there exists a qc-open immersion $\X \to \Spec \A$ for an $\E_\infty$ ring $\A$. Note this forces $\X$ to be equivalent to a spectral scheme. A map of stacks $\Y \to \X$ is \emph{quasi-affine} if for all maps $\Spec \A \to \X$, the pullback $\Spec \A \times_\X \Y$ is quasi-affine.
\end{mydef}

By \Cref{prop:qaffdiagonalmeansdescent}, nonconnective spectral Deligne--Mumford stacks with quasi-affine diagonal embed into $\Stk$. In particular, nonconnective spectral schemes embed into $\Stk$, allowing us to use the above language of quasi-affine morphisms inside $\Stk$.

\begin{remark}\label{rmk:quasi-affine-canonical}
    By \cite[Pr.2.4.1.3]{sag}, a stack $\sfX$ is quasi-affine if and only if the canonical map $\X \to \Spec \Ga(\X)$ is a qc-open immersion. This can be upgraded to a relative statement using \Cref{rmk:relativespectra}: a map of stacks $f\colon \Y \to \X$ is quasi-affine if and only if the canonical map $\Y \to \relSpec_\X(f_\ast\O_\Y)$ is a qc-open immersion.
\end{remark}

\begin{mydef}
    A map of stacks $\Y \to \X$ is a \emph{relative nonconnective spectral Deligne--Mumford stack} if for each map $\Spec \A \to \X$ the pullback $\Spec \A\times_\X \Y$ is represented by a nonconnective spectral Deligne--Mumford stack which satisfies flat descent (\Cref{df:flatdescentDMstacks}).
\end{mydef}

The following facts about quasi-affine morphisms are standard.

\begin{prop}\label{pr:basicqaffinefacts}
    Let $f\colon \Y \to \X$ and $\X' \to \X$ be maps of stacks.
    \begin{enumerate}
        \item If $f$ is quasi-affine, then the base change $\X'\times_\X \Y \to \Y$ is quasi-affine.
        \item If $\X=\Spec \A$ is affine, then $f$ is quasi-affine if and only if $\Y$ is quasi-affine.
        \item Suppose that both $\Y$ and $\X$ are connective and $f$ is a relative nonconnective Deligne--Mumford stack. If the underlying map of classical stacks $f^\heartsuit\colon \Y^\heartsuit \to \X^\heartsuit$ is affine (resp.\ quasi-affine), then $f$ is affine (resp.\ quasi-affine).
    \end{enumerate}
\end{prop}

\begin{proof}
    Parts 1 and 2 follow from the definition of quasi-affine morphisms. For part 3, suppose that both $\Y$ and $\X$ are connective. Any map $\Spec \A \to \X$ factors through $\Spec \tau_{\geq 0} \A \to \X$, so in the pullback diagram of stacks
    \[\begin{tikzcd}
        {\Y''}\ar[r]\ar[d, "{f''}"] &   {\Y'}\ar[d, "{f'}"]\ar[r]    &   {\Y}\ar[d, "{f}"]   \\
        {\Spec \A}\ar[r]   &   {\Spec \tau_{\geq 0}\A}\ar[r]    &   {\X,}
    \end{tikzcd}\]
    the right-hand square is Cartesian in both $\Stk$ and $\Stk^\cn$. As the functor $\tau_{\leq 0} \colon \Stk^\cn \to \Stk^\heartsuit$ is a right adjoint, we see that the diagram of classical stacks
    \[\begin{tikzcd}
        {\Y'^\heartsuit}\ar[r]\ar[d, "{f'^\heartsuit}"] & {\Y^\heartsuit}\ar[d, "{f^\heartsuit}"]   \\
        {\Spec \pi_0 \A}\ar[r] &   {\X^\heartsuit}
    \end{tikzcd}\]
    is also Cartesian. If $f^\heartsuit$ is affine, then the above diagram immediately shows that $\Y'^\heartsuit$ is affine. By \cite[Cor.1.4.7.3]{sag}, we see that $\Y'\simeq \Spec \B$ is itself affine, as by assumption it is a nonconnective spectral Deligne--Mumford stack. It then follows that $\Y''\simeq \Spec (\A\otimes_{\tau_{\geq0} \A} \B)$ is also affine, as desired. If $f^\heartsuit$ is quasi-affine, we see that $f'^\heartsuit$ is quasi-affine by base change, and that $f'$ is quasi-affine courtesy of \cite[Rmk.2.4.4.2]{sag}. It follows that $f''$ is also quasi-affine by parts 1 and 2 applied to the rectangle above.
\end{proof}

\subsubsection{Quasi-affine morphisms are universally 0-affine}

\begin{theorem}\label{ex:quasiaffineareU0A}
    A quasi-affine morphism of stacks $\Y \to \X$ is universally $0$-affine.
\end{theorem}

Our proof of this theorem requires a little preparation. We start by recalling some statements appearing in \cite{sag} concerning qc-open immersions; they will be generalized to quasi-affine morphisms as a consequence of \Cref{prop:univ0affinefacts} and \Cref{ex:quasiaffineareU0A}.

\begin{lemma}[{\cite[Pr.2.4.1.5(1)]{sag}}]\label{preservationofcolimitsopens}
    Let $f\colon \X \to \Spec \A$ be a qc-open immersion. Then the functor
    \[f_\ast \colon \QCoh(\X) \to \QCoh(\Spec \A)\]
    preserves colimits.
\end{lemma}

\begin{lemma}[{\cite[Pr.2.4.1.5(3)]{sag}}]\label{qcohpullpushopen}
    A Cartesian diagram of stacks
    \[\begin{tikzcd}
        {\X'}\ar[r, "{f'}"]\ar[d, "{g'}"']    &   {\Spec \A'}\ar[d, "g"]\\
        {\X}\ar[r, "{f}"]    &   {\Spec \A,}
    \end{tikzcd},\]
    where $f$, and hence $f'$, is a qc-open immersion is adjointable.
\end{lemma}

\begin{lemma}[{\cite[Cor.2.4.1.6]{sag}}]\label{openimmersionsareff}
    Let $f\colon \X\to \Spec \A$ be a qc-open immersion. Then the functor $f_\ast\colon \QCoh(\X) \to \Mod_\A$ is fully faithful.
\end{lemma}

Together, these statements provide us with another source of universally $0$-affine maps.

\begin{lemma}\label{lm:openimmersionsareU0A}
    A qc-open immersion $\X \to \Spec \A$ is universally $0$-affine.
\end{lemma}

\begin{proof}
    Consider the diagram (\ref{eq:defofU0A}) in this case
    \[\begin{tikzcd}
    \sfF'\ar[r, "{i'}"]\ar[d, "{g'}"]&\Spec \sfB'\ar[d, "g"]\\
    \sfF\ar[r, "i"]\ar[d]&\Spec \sfB\ar[d]\\
    \sfX\ar[r,"f"]&{\Spec \sfA}
    \end{tikzcd}\]
    where both squares are Cartesian. The map $i$ is a qc-open immersion, so the pushforward functor $i_\ast \colon \QCoh(\sfF) \to \Mod_{\B}$ preserves colimits by \Cref{preservationofcolimitsopens} and is conservative by \Cref{openimmersionsareff}. In particular, we see that $\sfF$ is $0$-affine by \Cref{prop:monoidalbarrbeck}. The upper square is adjointable by \Cref{qcohpullpushopen}, which concludes the proof.
\end{proof}

\begin{proof}[Proof of \Cref{ex:quasiaffineareU0A}]
    Consider the diagram of (\ref{eq:defofU0A}) of Cartesian squares in this case
    \[\begin{tikzcd}
        {\sfF'}\ar[r]\ar[d]   &   {\Spec\B}\ar[r]\ar[d] &   {\Spec \A'}\ar[d, "g"] \\
        {\sfF}\ar[r, "i"]\ar[d]   &   {\Spec \Ga(\O_\sfF)}\ar[r] &   {\Spec\A}\ar[d] \\
        {\Y}\ar[rr, "{f}"]   &&   {\X}
    \end{tikzcd}\]
    where $f$ is quasi-affine, the given map $\sfF \to \Spec \A$ is factored through the unit map $i\colon \sfF \to \Spec \Ga(\O_\sfF)$, and the upper squares are defined to be Cartesian. By definition, the stack $\sfF$ is quasi-affine, and the map $i$ is a qc-open immersion, see \cref{rmk:quasi-affine-canonical}. In particular, $i$ is universally $0$-affine by \Cref{lm:openimmersionsareU0A}. This implies that $\sfF$ is $0$-affine by part 2 of \cref{prop:univ0affinefacts}. This also yields the desired property of global sections of $\sfF'$
    \[\Ga(\O_\sfF) \otimes_\A \A' \simeq \B \simeq \B \otimes_{\Ga(\O_\sfF)} \Ga(\O_\sfF) \simeq \Ga(\O_{\sfF'}),\]
    where the first equivalence comes from the upper-right Cartesian square above, and the last equivalence from the fact that $i$ is universally $0$-affine, hence has this base change property. By \Cref{rmk:internalandexternalU0Adef}, this shows that $f$ is universally $0$-affine.
\end{proof}

\begin{remark}\label{rmk:BisactuallyGa(Z)}
    As a consequence of our proof above, we also see that $\B \simeq \Ga(\O_{\sfF'})$.  
\end{remark}

\subsection{Flat morphisms between geometric stacks}\label{ssec:flatnessandgeometricstacks}
To relate spectral algebraic geometry to classical algebraic geometry, we study \emph{flat} morphisms between nicely behaved stacks, the so-called \emph{geometric stacks}.

\subsubsection{Geometric stacks}\label{sssec:geometricstacks}
Stacks can always be written as a colimit of representables by \Cref{pr:stackisbasicallyatopos}. We give a particular name to those stacks which can be written as a colimit of representables and only flat maps between these representables. Recall the definition of an affine map of stacks from \Cref{def:affinemorphism}.

\begin{mydef}\label{def:affineandflattogether}
    An affine maps of stacks $\Y \to \X$ is \emph{flat} if for each $\Spec \A \to \Y$, the associated map $\Spec \B\simeq \Y\times_\X\Spec \A \to \Spec \A$ induces a flat map of $\E_\infty$ rings $\A\to \B$. A collection of maps $\{\Spec \A \to \X\}$ to a fixed stack $\X$ is an \emph{atlas} if $\coprod \Spec \A \to \X$ is an effective epimorphism and each map $\Spec \A \to \X$ is affine and flat. A \emph{geometric stack} is a stack $\X$ which admits an atlas. We write $\GeoStk\subset\Stk$ for the full subcategory of geometric stacks.
\end{mydef}

This notion of geometric stacks is more general than that often seen in the literature. For example, Gregoric \cite[Df.3.1.5]{rokeven} assumes that $\X$ has affine diagonal, and that the coproduct $\coprod \Spec \A$ is affine, further restricting ourselves to the quasi-compact case. Lurie \cite[Df.9.3.0.1]{sag} also makes these assumptions and further restricts to the connective case. We prefer our generalized definition as it has more pleasant categorical properties.

\begin{example}\label{ex:SPDMaregeometric}
    Let $\X$ be a nonconnective spectral Deligne--Mumford stack with affine diagonal, regarded as a stack via \cref{prop:dmstacksintoflatstacks} and \cref{prop:qaffdiagonalmeansdescent}. Then $\X$ is geometric. Indeed, $\X$ can be written as a colimit of \'etale maps $\Spet\A\to\X$ in $\SpDMncet$, all of which are affine as $\X$ has affine diagonal. As $\SpDMncet\to\Stk$ preserves colimits by \Cref{prop:dmstacksintoflatstacks}, this implies that $\X$ is a geometric stack.
\end{example}

Some facts about affine maps which are also flat follow immediately from the definitions.

\begin{prop}\label{prop:simpleflataffinefacts}
    Let $f\colon \Y \to \X$ and $g\colon \sfZ \to \Y$ be maps of stacks.
    \begin{enumerate}
        \item If $g$ and $f$ are both affine and flat, then their composite $f\circ g$ is affine and flat.
        \item For any other map $\X' \to \X$, if $f$ is affine and flat, then the pullback $\Y \times_\X \X' \to \X'$ is affine and flat.
    \end{enumerate}
\end{prop}

\begin{proof}
    Both of these properties follow from the analogous facts for affine maps of stacks and flat maps of $\E_\infty$ rings.
\end{proof}

By definition, geometric stacks admit an atlas, and so can be written as a colimit of the \v{C}ech nerve of this atlas. In fact, a weaker condition characterizes geometric stacks.

\begin{lemma}\label{lm:geometrichasflatdiagram}
    A stack $\X$ is geometric if and only if there is some presentation $\X \simeq \colim \Spec \A$ of $\X$ as a (small) colimit of representables, such that the maps $\Spec \A \to \X$ are affine and flat. Moreover, if $\X$ is a geometric stack, then the diagram $\colim \Spec \A \simeq \X$ can be chosen such that each map $\A \to \A'$ in this diagram is flat.
\end{lemma}

Due to this lemma, we will often abuse notation between an atlas for a geometric stack and a colimit presentation $\X \simeq \colim \Spec \A$, where all of the maps $\Spec \A \to \X$ are affine and flat.

\begin{proof}
    If $\X$ is geometric, then we can take such a presentation of $\X$ via the Čech nerve of the map $\coprod_{I} \Spec \A_i \to \X$. Indeed, writing $\A_{i_1,\ldots, i_n}$ for the affine representing the $k$-fold pullback of the flat affine maps $\Spec \A_{i_k} \to \X$, we see the Čech nerve takes the form
    \[
    \begin{tikzcd} {\coprod_{I} \Spec \A_i} & {\coprod_{I^2} \Spec \A_{i,j}} \ar[l,shift left=1mm]\ar[l,shift right=1mm]&{\coprod_{I^3} \Spec \A_{i,j,k}} \ar[l,shift left=2mm]\ar[l]\ar[l,shift right=2mm] & \cdots\ar[l,shift left=3mm]\ar[l,shift left=1mm]\ar[l,shift right=1mm]\ar[l,shift right=3mm]\end{tikzcd}
    \]
    which shows that $\X \simeq \colim \Spec \A_{i_1,\ldots, i_n}$. This also shows the ``in particular'' statement. Conversely, if we can write $\X$ as $\X \simeq \colim \Spec \A$ where each $\Spec \A \to \X$ is flat and affine, then taking the $0$-truncation of the small indexing category for this colimit, we obtain the desired effective epimorphism $\coprod \Spec \A \to \X$.
\end{proof}

\subsubsection{Flat maps}\label{sssec:flatmapsbetweengeometricstacks}
We will only consider flat maps between geometric stacks, as this gives us an atlas to test flatness against.

\begin{mydef}\label{def:flatmapsofgeometricstacks}
A morphism of geometric stacks $f\colon \Y \to \X$ is \emph{flat} if for an atlas $\coprod \Spec \A \to \X$ of $\X$ and an atlas $\coprod \Spec \B \to \Y$ of $\Y$, the base change of $\Spec \B \to \Y \to \X$ along each $\Spec \A \to \X$
\[\Spec \sfC \simeq \Spec \B \times_\X \Spec \A \to \Spec \A\]
induces a flat map of $\E_\infty$ rings $\A \to \B$. Such a morphism $f$ is \emph{faithfully flat} if $f$ is flat and an epimorphism. 
\end{mydef}

These maps are a little different from those of \Cref{df:affinecover}, as they are not necessarily affine.

\begin{example}
    If $f\colon \Y \to \X$ is a map of nonconnective spectral Deligne--Mumford stacks both with affine diagonal, then $f$ is flat as a map of geometric stacks if and only if $f$ is flat as a map of nonconnective spectral Deligne--Mumford stacks as defined in \cite[Df.2.8.2.1]{sag}. Indeed, this is because a separated nonconnective spectral Deligne--Mumford stack $\X$ has an atlas where all maps to $\X$ are \'{e}tale.
\end{example}

Next, we observe some simple properties of flat maps between geometric stacks.

\begin{prop}\label{prop:simplepropertiesofflatness}
    Consider a pullback of geometric stacks
    \[\begin{tikzcd}
        {\Y'}\ar[r, "h'"]\ar[d, "f'"]   &   {\Y}\ar[d, "f"] \\
        {\X'}\ar[r, "h"]               &   {\X,}
    \end{tikzcd}\]
    where $h$ is affine.
    \begin{enumerate}
        \item If $f$ is flat, then $f'$ is flat.
        \item If $h$ is faithfully flat and $f'$ is flat, then $f$ is flat.
    \end{enumerate}
\end{prop}

\begin{proof}
    These two arguments are both diagram chases. For both statements, consider the following commutative diagram of stacks
\[\begin{tikzcd}
	{\Spec \sfC'} && {\Y'_\A} && {\Spec \A'} \\
	& {\Spec \sfC} && {\Y_\A} && {\Spec \A} \\
	{\Spec \B'} && {\Y'} && {\X'} \\
	& {\Spec \B} && \Y && \X
	\arrow[from=1-1, to=1-3]
	\arrow[from=1-1, to=2-2]
	\arrow[from=1-1, to=3-1]
	\arrow[from=1-3, to=1-5]
	\arrow[from=1-3, to=2-4]
	\arrow[from=1-3, to=3-3]
	\arrow[from=1-5, to=2-6]
	\arrow[from=1-5, to=3-5]
	\arrow[from=2-2, to=2-4]
	\arrow[from=2-2, to=4-2]
	\arrow[from=2-4, to=2-6]
	\arrow[from=2-4, to=4-4]
	\arrow["g", from=2-6, to=4-6]
	\arrow[from=3-1, to=3-3]
	\arrow[from=3-1, to=4-2, "h''"]
	\arrow[from=3-3, to=3-5]
	\arrow[from=3-3, to=4-4, "h'"]
	\arrow["h"', from=3-5, to=4-6]
	\arrow[from=4-2, to=4-4, "i"]
	\arrow["f"', from=4-4, to=4-6]
\end{tikzcd}\]
    where $g$ is part of an atlas for $\X$, $i$ is part of an atlas for $\Y$, and each square is a pullback. Suppose that $f$ is flat. Notice that by pullback, that $\Spec \A' \to \X$ is affine and forms part of an atlas for $\X'$, and similarly that $\Spec \B' \to \Y'$ is affine and forms part of an atlas for $\Y'$. By definition, to see that $f'\colon \Y' \to \X'$ is flat, it suffices to show that $\A' \to \sfC'$ is flat. By the pasting lemma, we have $\sfC'\simeq \sfC\otimes_\A \A'$, so in particular $\A' \to \sfC'$ is flat as $\A \to \sfC$ is flat by assumption.

    Suppose now that $h$ is faithfully flat and that $f'$ is flat. Then we similarly have $\sfC'\simeq \sfC\otimes_\A \A'$ with $\A'\to \sfC'$ flat and $\A \to \A'$ faithfully flat. It follows from \cite[Pr.2.8.4.2(6)]{sag} that flat morphisms are detected fpqc locally, so that $\A \to \sfC$ is flat as desired.
\end{proof}

It follows that the definition of flatness is invariant under the choice of atlas.

\begin{cor}\label{cor:invariantofflatness}
    A morphism of geometric stacks $f\colon \Y \to \X$ is flat (resp.\ faithfully flat) if for all atlases $\coprod \Spec \A \to \X$ of $\X$ and for all atlases $\coprod \Spec \B \to \Y$ of $\Y$, the base change of the composite $\Spec \B \to \Y \to \X$ along each $\Spec \A \to \X$ is flat (resp.\ faithfully flat) as a map of $\E_\infty$ rings.
\end{cor}

\subsubsection{From derived stacks to underlying stacks}
Flat morphisms between geometric stacks give a connection between stacks, their connective covers, and their underlying classical stacks. We begin with the affine case.

\begin{lemma}\label{heartsandbasechangeAFFINE}
    Consider a square of affine stacks and the associated squares of connective and classical stacks
	\begin{equation}\label{eq:affineflatness}\begin{tikzcd}
	{\Spec \B'}\ar[r]\ar[d, "f'"]	&	{\Spec\A'}\ar[d, "f"]	\\
	{\Spec \B}\ar[r]		&	{\Spec \A}
	\end{tikzcd}\qquad 
    \begin{tikzcd}
	{\Spec \tau_{\geq 0}\B'}\ar[r]\ar[d, "{\tau_{\geq 0} f'}"]	&	{\Spec\tau_{\geq 0}\A'}\ar[d, "{\tau_{\geq 0} f}"]	\\
	{\Spec \tau_{\geq 0}\B}\ar[r]		    &	{\Spec \tau_{\geq 0}\A,}
	\end{tikzcd}\qquad 
    \begin{tikzcd}
	{\Spec \pi_0\B'}\ar[r]\ar[d, "{\pi_0 f'}"]	&	{\Spec\pi_0\A'}\ar[d, "{\pi_0 f}"]	\\
	{\Spec \pi_0\B}\ar[r]		    &	{\Spec \pi_0\A,}
	\end{tikzcd}\end{equation}
    and assume that both $f$ and $f'$ are flat. Then the left-hand square is Cartesian if and only if the middle square is Cartesian if and only if the right-hand square is Cartesian.
\end{lemma}

\begin{proof}
    Consider the commutative diagram of graded abelian groups
    \[\begin{tikzcd}
        {\pi_\ast\B \otimes_{\pi_0 \A} \pi_0 \A' \simeq \pi_\ast \B \otimes_{\pi_0 \B} (\pi_0 \B \otimes_{\pi_0\A} \pi_0 \A')}\ar[r]\ar[d, "{\simeq}"]    &   {\pi_\ast \B \otimes_{\pi_0 \B} \pi_0 \B'}\ar[d, "{\simeq}"]    \\
        {\pi_\ast (\B \otimes_{\A} \A')}\ar[r]    &   {\pi_\ast \B',}
    \end{tikzcd}\]
    where the vertical maps are equivalences courtesy of Künneth spectral sequences concentrated in filtration $0$ from our flatness assumptions. This shows that the upper horizontal map is an equivalence if and only if the lower horizontal map is an equivalence, i,e.\ that the right-hand square of (\ref{eq:affineflatness}) is Cartesian if and only if the left-hand square of (\ref{eq:affineflatness}) is Cartesian. If we consider the nonnegatively graded components of the above diagram, we see that the middle square is Cartesian if and only if the right-hand square is Cartesian.
\end{proof}

This globalizes as follows.

\begin{prop}\label{flatnessislovely}
    Consider a square of geometric stacks and the associated squares of connective and classical stacks
	\begin{equation}\label{eq:globalflatness}\begin{tikzcd}
	{\Y'}\ar[r]\ar[d, "f'"]	&	{\X'}\ar[d, "f"]	\\
	{\Y}\ar[r]		&	{\X}
	\end{tikzcd}\qquad \qquad
    \begin{tikzcd}
	{\tau_{\geq 0}\Y'}\ar[r]\ar[d, "{\tau_{\geq 0} f'}"]	&	{\tau_{\geq 0}\X'}\ar[d, "{\tau_{\geq 0} f}"]	\\
	{\tau_{\geq 0}\Y}\ar[r]		    &	{\tau_{\geq 0}\X,}
	\end{tikzcd}\qquad \qquad
    \begin{tikzcd}
	{\X'^\heartsuit}\ar[r]\ar[d, "{f'^\heartsuit}"]	&	{\X'^\heartsuit}\ar[d, "{f^\heartsuit}"]	\\
	{\Y^\heartsuit}\ar[r]		    &	{\X^\heartsuit,}
	\end{tikzcd}\end{equation}
    and assume that both $f$ and $f'$ are flat. If the left-hand square of stacks is Cartesian, then both the middle and right-hand squares are also Cartesian.
\end{prop}

We note that although $f$ is flat, it does not immediately follow that $f'$ is also flat. Nevertheless, this is often the case, such as when all stacks in sight are affine, or when the map $\Y \to \X$ is affine; see \Cref{prop:simplepropertiesofflatness} and \Cref{heartsandbasechangeAFFINE}.

\begin{proof}
    We will only show that the right-hand square is Cartesian; showing that the middle square is Cartesian follows from the same arguments.

    First, suppose that both $\X=\Spec \A$ and $\Y=\Spec \B$ are affine. Fix an atlas for $\X'$, so a presentation $\X' \simeq \colim \Spec \A'$ where each $\Spec \A' \to \X'$ is affine and flat. 
     Pulling each map $\Spec \A' \to \X' \to \Spec \A$ back along $\Spec \B \to \Spec \A$, we obtain a similar atlas for $\Y'$:
    \[\Y' \simeq \Spec \B \underset{\Spec \A}{\times} \X' \simeq \Spec \B \underset{\Spec \A}{\times} \colim \Spec \A' \simeq \colim \Spec \B \underset{\Spec \A}{\times} \Spec \A' = \colim \Spec \B'\]
    courtesy of \Cref{pr:stackisbasicallyatopos}. Using \Cref{cor:coversandunderlyingascolimits}, we also obtain equivalences $\X'^\heartsuit \simeq \colim \Spec \pi_0 \A'$ and $\Y'^\heartsuit \simeq \colim \Spec \pi_0 \B'$. Another application of \Cref{pr:stackisbasicallyatopos} together \Cref{heartsandbasechangeAFFINE} shows that the right-hand square above is Cartesian:
    \[\Y'^\heartsuit \simeq \colim \Spec \pi_0 \B' \simeq \colim \Spec \pi_0 \B \underset{\Spec \pi_0 \A}{\times} \Spec \pi_0 \A'\]
    \[\simeq \Spec \pi_0\B \underset{\Spec \pi_0 \A}{\times} \colim \Spec \pi_0\A' \simeq \Spec \pi_0\B \underset{\Spec \pi_0 \A}{\times} \X'^\heartsuit.\]

    Next, suppose that $\X=\Spec \A$ is affine and $\Y$ is a general geometric stack, with some chosen atlas presentation $\Y \simeq \colim \Spec \B$. By base change, this also induces a presentation $\Y' \simeq \colim \Y'_\B$ of $\Y'$. For each affine and flat map $\Spec \B \to \Y$ in this presentation, consider the Cartesian squares of stacks
    \[\begin{tikzcd}
        {\Y'_\B}\ar[r]\ar[d]    &   {\Y'}\ar[d]\ar[r]   &   {\X'}\ar[d] \\
        {\Spec \B}\ar[r]        &   {\Y}\ar[r]          &   {\Spec \A.}
    \end{tikzcd}\]
    The outer rectangle is a pullback, so from the previous paragraph, it induces an equivalence of underlying classical stacks
    \[\Y'^\heartsuit_\B \simeq \X'^\heartsuit \underset{\Spec \pi_0\A}{\times} \Spec \pi_0 \B.\]
    Moreover, \Cref{cor:coversandunderlyingascolimits} induces the equivalences of stacks
    \[\Y^\heartsuit \simeq \colim \Spec \pi_0\B, \qquad \Y'^\heartsuit \simeq \colim \Y'^\heartsuit_\B.\]
    Combining all of these equivalences together with \Cref{pr:stackisbasicallyatopos} yields the desired expression:
    \[\Y'^\heartsuit \simeq \colim \Y'^\heartsuit_\B \simeq \colim \X'^\heartsuit \underset{\Spec \pi_0\A}{\times} \Spec \pi_0 \B\]
    \[\simeq \X'^\heartsuit \underset{\Spec \pi_0\A}{\times} \colim \Spec \pi_0 \B \simeq \X'^\heartsuit \underset{\Spec \pi_0\A}{\times} \Y^\heartsuit.\]

    Finally, let $\X$ and $\Y$ be general, and fix an atlas presentation $\X \simeq \colim \Spec \A$ of $\X$. By \Cref{pr:stackisbasicallyatopos}, this induces an equivalence of categories
    \[\Stk_{/\X} \xrightarrow{\simeq} \lim \Stk_{/\Spec \A}, \qquad (\Y \to \X) \mapsto (\Y\underset{\X}{\times} \Spec \A \to \Spec \A).\]
    This equivalence sends the left-hand square of this proposition to the collection of Cartesian squares of stacks on the left-hand side below,
    \[\begin{tikzcd}
        {\Y'_\A}\ar[r]\ar[d]    &   {\X'_\A}\ar[d]  \\
        {\Y_\A}\ar[r]           &   {\Spec \A,}
    \end{tikzcd}\qquad  \overset{(-)^\heartsuit}{\longmapsto} \qquad 
    \begin{tikzcd}
        {\Y'_{\A^\heartsuit}}\ar[r]\ar[d]    &   {\X'_{\A^\heartsuit}}\ar[d]  \\
        {\Y_{\A^\heartsuit}}\ar[r]           &   {\Spec \pi_0\A.}
    \end{tikzcd}\]
    This, in turn, induces the right-hand collection of Cartesian squares of underlying classical stacks by the previous paragraph. As all of these squares of classical stacks are Cartesian, and the inverse of this collections of squares under the equivalence of categories
    \[\Stk^\heartsuit_{/\X^\heartsuit} \xrightarrow{\simeq} \lim \Stk^\heartsuit_{/\Spec \pi_0\A}\]
    is precisely the desired right-hand square of this proposition, we see that this desired square is Cartesian.
\end{proof}

\begin{remark}\label{rmk:adamsflatness}
Suppose we have a map $f\colon \sfY \to \sfX$ between geometric stacks such that the fiber product $\Y \times_\X \Y$ is also geometric. Furthermore, suppose that $f$ is Adams flat à la \Cref{def:adamsflat}, meaning that either one of the projections $\Y \times_\X \Y \to \Y$ is affine and for each $\Spec \A \to \Y$, the pullback $\Y \times_\X \Spec \A \to \Y$ is a faithfully flat map between affine stacks. In particular, the projections are both flat maps between geometric stacks, and so are the projections between the iterated fiber products as affine flat maps are closed under base-change. In this case, the underlying classical stack of the associated descent stack $\D_f$ is a groupoid object in $\Stk^\heartsuit$. To see this, we use \Cref{lem:descentstackbasicproperties} to write $\D_f \simeq \colim \Y^{\times_\X(\bullet+1)}$, and then \Cref{cor:coversandunderlyingascolimits} to write
\[\D_f^\heartsuit \simeq \left(\colim \Y^{\times_\X(\bullet+1)}\right)^\heartsuit \simeq \colim (\Y^{\times_\X(\bullet+1)})^\heartsuit.\]
Finally, the Adams flatness condition combines with \Cref{flatnessislovely} to show that all of the terms $(\Y^{\times_\X(\bullet+1)})^\heartsuit$ in the above colimit are naturally equivalent to the $\bullet$-fold product of $(\Y \times_\X \Y)^\heartsuit$ over $\X^\heartsuit$, in other words, this simplicial diagram of classical stacks satisfies the Segal condition. In particular, this witnesses $\D_f^\heartsuit$ as the object associated to the groupoid pair $(\Y^\heartsuit, (\Y \times_\X \Y)^\heartsuit)$. 
\end{remark}

We will see an application of this in chromatic homotopy theory in \Cref{underlyingclassicalstackofDMUP}.

As a useful application of the above statement, we can often show that properties of morphisms of stacks descend to properties of morphisms between the underlying classical stacks as soon as the map is also flat. The following corollary is an example.

\begin{cor}\label{pr:heartsofquasiaffinenadaffine}
    Let $f\colon \Y \to \X$ be a flat map between geometric stacks. If $f$ is affine (resp.\ quasi-affine), then $f^\heartsuit$ is affine (resp.\ quasi-affine).
\end{cor}

\begin{proof}
    Let $\amalg \Spec \A \to \X$ be an atlas for $\X$, inducing an atlas $\amalg\Spec \pi_0 \A \to \X^\heartsuit$ for $\X^\heartsuit$. By \Cref{flatnessislovely}, the left-hand Cartesian square of stacks below induces the right-hand Cartesian square of underlying classical stacks:
    \[\begin{tikzcd}
        {\Y_\A}\ar[r, "{f_\A}"]\ar[d] &   {\Spec \A}\ar[d]    \\
        {\Y}\ar[r, "f"]     &   {\X,}
    \end{tikzcd}\qquad \begin{tikzcd}
        {\Y_\A^\heartsuit}\ar[r, "{f_\A^\heartsuit}"]\ar[d]          &   {\Spec \pi_0 \A}\ar[d]    \\
        {\Y^\heartsuit}\ar[r, "{f^\heartsuit}"] &   {\X^\heartsuit}
    \end{tikzcd}\]
    First, suppose that $f$ is affine, so that $\Y_\A \simeq \Spec \B$ is affine. To show that $f^\heartsuit$ is affine, it suffices to show this fpqc-locally on the base following \cite[\href{https://stacks.math.columbia.edu/tag/02L5}{02L5}]{stacks}. In other words, it suffices to show that $\Spec \pi_0 A \times_{\X^\heartsuit} \Y^\heartsuit$ is affine. But as observed above, this pullback is equivalent to $\Y_\A^\heartsuit = (\Spec\B)^\heartsuit \simeq \Spec \pi_0 \B$, and so affine.

    Similarly, being quasi-affine is fpqc local on the base \cite[\href{https://stacks.math.columbia.edu/tag/02L7}{02L7}]{stacks}. If $f$ is quasi-affine, then $\Y_\A$ is quasi-affine. We conclude that $\Y_\A^\heartsuit$ is also a quasi-affine classical scheme,as open immersions of spectral schemes induce open immersions on underlying schemes. So $f_\A^\heartsuit$ is quasi-affine for all $\A$, proving the statement.
\end{proof}

\subsubsection{From underlying stacks back to derived stacks}
To obtain some useful converses to \Cref{flatnessislovely} and \Cref{pr:heartsofquasiaffinenadaffine}, we place some more assumptions on our maps of stacks.

\begin{lemma}
    Let $f\colon \Y \to \X$ be a flat map of geometric stacks which is moreover a relative nonconnective spectral Deligne--Mumford stack. Then $f$ is an equivalence if and only if $f^\heartsuit$ is an equivalence.
\end{lemma}

\begin{proof}
    The ``only if'' statement is clear. So assume that $f^\heartsuit$ is an equivalence. Clearly, $f$ being an equivalence can be tested fpqc locally. So consider an atlas presentation $\X \simeq \colim \Spec \A$, and for each affine and flat map $\Spec \A \to \X$ in this presentation, consider the pullback square of stacks
    \[\begin{tikzcd}
        {\Y_\A}\ar[r, "{f_\A}"]\ar[d] &   {\Spec \A}\ar[d]    \\
        {\Y}\ar[r, "f"]          &   {\X,}
    \end{tikzcd}\qquad \overset{(-)^\heartsuit}{\longmapsto} \qquad 
    \begin{tikzcd}
        {\Y_\A^\heartsuit}\ar[r, "{f_\A^\heartsuit}"]\ar[d]  &   {\Spec \pi_0\A}\ar[d]    \\
        {\Y^\heartsuit}\ar[r, "{f^\heartsuit}"]      &   {\X^\heartsuit,}
    \end{tikzcd}\]
    as well as the associated commutative diagram of underlying classical stacks. In particular, $\colim \Y_\A \simeq \Y$ gives a presentation for $\Y$. By \Cref{prop:simplepropertiesofflatness}, $f_\A$ is also flat as $\Spec \A \to \X$ is affine, and by \Cref{flatnessislovely}, we see that the right-hand square of classical stacks above is also a pullback. By assumption, $f^\heartsuit$ is an equivalence, and so $f_\A^\heartsuit$ is also an equivalence. The stack $\Y_\A$ is a nonconnective spectral Deligne--Mumford stack by assumption whose underlying stack is affine as $f_\A^\heartsuit$ is an equivalence. Lurie shows in \cite[Cor.1.4.7.3]{sag} that if a nonconnective spectral Deligne--Mumford stack $\sfZ$ has affine underlying classical stack, then $\sfZ$ itself is affine. In our case, this shows that $\Y_\A = \Spec \sfC$ is affine. The map $f_\A \colon \Spec \sfC \to \Spec \A$ is then a flat map between affines inducing an isomorphism on $\pi_0$, so is an equivalence by definition of flatness.
\end{proof}

\begin{lemma}\label{lem:connectivecoversofflatareflat}
    Let $f\colon \Y \to \X$ be a flat map of geometric stacks. Then $\tau_{\geq 0}f$, viewed as a map of geometric stacks in $\Stk$, is flat.
\end{lemma}

\begin{proof}
    First, fix atlas presentations for $\X$ and $\Y$:
    \[\X \simeq \colim \Spec \A, \qquad \Y \simeq \colim \Spec \B.\]
    By \Cref{cor:coversandunderlyingascolimits}, we have presentations for the connective covers of these stacks too:
    \[\tau_{\geq 0}\X \simeq \colim \Spec \tau_{\geq 0}\A, \qquad \tau_{\geq 0}\Y \simeq \colim \Spec \tau_{\geq 0}\B.\]
    By assumption, we know that in the Cartesian diagram of stacks
    \[\begin{tikzcd}
        {\Spec \sfC}\ar[d]\ar[rr]   &&  {\Spec \A}\ar[d]    \\
        {\Spec \B}\ar[r]   &    {\Y}\ar[r]  &   {\X,}
    \end{tikzcd}\]
    the map of $\E_\infty$ rings $\A \to \sfC$ is flat. By \Cref{flatnessislovely}, the associated diagram upon applying the functor $\tau_{\geq 0}$ is also Cartesian. By definition, the desired map $\tau_{\geq 0}\A \to \tau_{\geq0} \sfC$ is then also flat.
\end{proof}

\begin{lemma}\label{lem:connectiveflatsquareCartesian}
    Let $f\colon \Y \to \Spec \A$ be a flat map of geometric stacks. Then the natural square of stacks
    \[\begin{tikzcd}
        {\Y}\ar[r]\ar[d, "f"]   &   {\tau_{\geq 0}\Y}\ar[d, "{\tau_{\geq0}f}"]  \\
        {\Spec \A}\ar[r]        &   {\Spec \tau_{\geq 0} \A}
    \end{tikzcd}\]
    is a pullback.
\end{lemma}

\begin{proof}
    Let $\Y \simeq \colim \Spec \B$ be an atlas presentation for $\Y$. In particular, by definition we see that each composite $\Spec \B \to \Spec \A$ is also flat, hence we have a Cartesian square of affine stacks
    \[\begin{tikzcd}
        {\Spec \B}\ar[r]\ar[d]  &   {\Spec \tau_{\geq 0}\B}\ar[d]   \\
        {\Spec \A}\ar[r]        &   {\Spec \tau_{\geq 0}\A.}
    \end{tikzcd}.\]
    Using the above atlas presentation for $\Y$ and the associated presentation for $\tau_{\geq 0}\Y$ from \Cref{cor:coversandunderlyingascolimits}, we obtain the desired expression again using \Cref{pr:stackisbasicallyatopos}:
    \[\Y \simeq \colim \Spec \B \simeq \colim \Spec \A \underset{\Spec \tau_{\geq0} \A}{\times} \Spec \tau_{\geq 0} \B\]
    \[\simeq \Spec \A \underset{\Spec \tau_{\geq0} \A}{\times} \colim \Spec \tau_{\geq0}\B  \simeq \Spec \A \underset{\Spec \tau_{\geq0} \A}{\times} \tau_{\geq0}\Y \qedhere\]
\end{proof}

One corollary of the above statement is the following well-known fact, implicitly providing the justification for many of Lurie's definitions in \cite{ec2}; see \cite[Lm.2.1]{elltempcomp}, for example. For an $\E_\infty$ ring $\A$, write $\GeoStk_{/\Spec \A}^\flat$ for the full subcategory of geometric stacks over $\Spec \A$ such that the structure map to $\Spec \A$ is flat.

\begin{cor}
    Let $\A$ be an $\E_\infty$ ring. Then the connective cover functor $\tau_{\geq0}\colon \Stk_{/\Spec \A} \to \Stk_{/\Spec \tau_{\geq 0}}$ induces an equivalence of categories
    \[\GeoStk_{/\Spec \A}^\flat \xrightarrow{\simeq} \GeoStk_{/\Spec \tau_{\geq 0}\A}^\flat.\]
\end{cor}

\begin{proof}
    The connective cover adjunction of \Cref{th:coversandtruncations} induces an adjunction
    \[\tau_{\geq 0}\colon \Stk_{/\Spec\A} \rightleftarrows \Stk_{/\Spec \tau_{\geq 0}\A} \colon ({-})\times_{\Spec \tau_{\geq 0}\A} \Spec \A,\]
    see \cite[Ex.3.26]{solosil} for example.
    By definition, connective covers of flat maps are flat and the base change of flat maps along an affine map like $\Spec \A \to \Spec \tau_{\geq0}\A$ are flat by \Cref{prop:simplepropertiesofflatness}. In particular, this also shows that both the right and left adjoints above preserve geometric stacks. We now consider the induced adjunction between subcategories of geometric stacks that are flat over their base. The unit of this adjunction is clearly an equivalence, and the counit is an equivalence by \Cref{lem:connectiveflatsquareCartesian}.
\end{proof}

Another immediate corollary of the techniques used to prove \Cref{lem:connectiveflatsquareCartesian} is the following analogue connecting connective stacks with their underlying stacks.

\begin{lemma}\label{lem:connectivetoHEARTflatsquareCartesian}
    Let $f\colon \Y \to \Spec \A$ be a flat map of connective geometric stacks. Then the natural square of stacks
    \[\begin{tikzcd}
        {\Y^\heartsuit}\ar[r]\ar[d, "f^\heartsuit"]   &   {\Y}\ar[d, "{f}"]  \\
        {\Spec \pi_0\A}\ar[r]        &   {\Spec \A}
    \end{tikzcd}\]
    is a pullback.
\end{lemma}

\begin{proof}
    This follows along the same lines as the proof of \Cref{lem:connectiveflatsquareCartesian}; we leave the details to the reader.
\end{proof}

From these statements, one can check if certain morphisms between geometric stacks are (quasi-) affine by inspecting the associated underlying map of stacks. This is a partial converse to \Cref{pr:heartsofquasiaffinenadaffine}.

\begin{cor}\label{cor:quasiaffinenessfromtheheart}
    Let $f\colon \Y \to \X$ be a flat map between geometric stacks which is also a relative nonconnective spectral Deligne--Mumford stack and where $\X$ has affine diagonal. If $f^\heartsuit$ is affine (resp.\ quasi-affine), then $f$ is affine (resp.\ quasi-affine).
\end{cor}

\begin{proof}
    For each map $\Spec \A \to \X$, consider the pullback $\Y_\A = \Y \times_\X \Spec \A$. As $f$ is flat, \Cref{flatnessislovely} yields the equivalence of underlying classical stacks
    \[\Y^\heartsuit_\A \simeq \Y^\heartsuit\underset{\X^\heartsuit}{\times} \Spec \pi_0 \A.\]
    Suppose now that $f^\heartsuit$ is affine, so in particular, $\Y_\A^\heartsuit$ is affine. By \cite[Cor.1.4.7.3]{sag}, this immediately implies that $\Y_\A$ is itself affine, hence $f$ is affine.

    Suppose now that $f^\heartsuit$ is merely quasi-affine, so $\Y_\A^\heartsuit$ is quasi-affine. By part 3 of \Cref{pr:basicqaffinefacts}, this means that $\tau_{\geq 0} \Y_\A$, or equivalently $\tau_{\geq 0} f_\A\colon \tau_{\geq 0} \Y_\A \to \Spec \tau_{\geq 0} \A$, is quasi-affine. As $\X$ has affine diagonal, the map $\Spec \A \to \X$ is affine, hence by \Cref{prop:simplepropertiesofflatness}, the base change map $f_\A \colon \Y_\A \to \Spec \A$ is also flat. This means we can use \Cref{lem:connectiveflatsquareCartesian} to see that $f_\A$ is the pullback of the quasi-affine map $\tau_{\geq 0} f_\A$. In particular, the base change property of quasi-affine map (\Cref{pr:basicqaffinefacts}), shows that $f_\A$ is quasi-affine, meaning that $\Y_\A$ is quasi-affine, as desired.
\end{proof}

\section{Complex-periodic geometry}\label{sec:complexperiodicgeometry}
In this section, we again work within the category $\Stk$ of \Cref{def:categoryoffpqcstacks}. The connection between the nonconnective spectral algebraic geometry of \Cref{sec:nonconsag} and chromatic homotopy theory comes from restricting our attention to \emph{complex-periodic stacks}. 

Recall that an $\E_\infty$ ring $\A$ is \emph{complex-periodic} if:
\begin{enumerate}
    \item $\A$ is \emph{weakly $2$-periodic}, meaning that the shift $\Sigma^2 \A$ is a locally free $\A$-module of rank $1$, or equivalently, $\pi_n \A \otimes_{\pi_0 \A} \pi_2 \A \to \pi_{n+2}\A$ is an isomorphism for all $n \in\Z$, and
    \item $\A$ is \emph{complex-orientable}, meaning the unit $\Sph \to \A$ admits a (non-fixed) factorization through $\Sph \simeq \Sigma^{\infty-2} S^2 \to \Sigma^{\infty-2} K(\Z,2)$.
\end{enumerate}

One salient feature of the adjective ``complex-periodic'' is that if $\A$ if complex-periodic, and $\A \to \B$ is a map of $\E_\infty$ rings, then $\B$ is complex-periodic; see \cite[Rmk.4.1.10]{ec2}. This is not true for the adjective ``even'', for example. This section is concerned with the following collection of stacks.

\begin{mydef}\label{def:locallycomplexperiodic}
    A stack $\X$ is \emph{complex-periodic} if for each map $\Spec \A\to \X$, the $\E_\infty$ ring $\A$ is complex-periodic.
\end{mydef}

 We call the study of such objects \emph{complex-periodic geometry}. We begin by comparing complex periodic stacks with stacks over $\M_\FG^\ori$, defined below.

\subsection{The moduli stack of oriented formal groups}\label{ssec:MFGOR}
The classical moduli stack of formal groups $\M_\FG^\heartsuit$ is defined as the presheaf of groupoids on the category of commutative rings which sends a commutative ring $R$ to the groupoid of formal groups and isomorphisms over $R$. To see this is an fpqc sheaf, we refer to \cite[Th.2.30]{goerssquasicoherent}.

In this subsection, we explore Lurie's moduli stack of \emph{oriented} formal groups $\M_\FG^\ori$ and study how it fits into our framework of \Cref{sec:nonconsag}.

\subsubsection{The points of $\M_{\FG}^\ori$}
In \cite{ec2}, Lurie defines \emph{oriented formal groups} over $\E_\infty$ rings.

\begin{mydef}
    A \emph{formal group} over an $\E_\infty$ ring $\A$ is an abelian group object in formal hyperplanes over $\A$, defined to be the subcategory of $\Fun(\CAlg^\cn_{\tau_{\geq 0} \A}, \Spc)$ spanned by cospectra of smooth coalgebras over $\tau_{\geq 0}\A$; see \cite[Df.1.5.10 \& Rmk.1.6.5]{ec2}. A \emph{preorientation} of a formal group $\widehat{\G}$ over $\A$ is a map of pointed spaces $e\colon S^2 \to \widehat{\G}(\tau_{\geq 0}\A)$. A preorientation $e$ is an orientation if $\A$ is complex-periodic and the map of formal groups $\widehat{\G}^\QQ_\A \to \widehat{\G}$ associated to $e$, where $\widehat{\G}^\QQ_\A$ is the \emph{Quillen formal group} of $\A$, see \cite[Con.4.1.13]{ec2}, is an equivalence.
\end{mydef}

\begin{mydef}\label{defofmfg}
    Let $\M_\FG\colon \CAlg \to \Spc$ be the \emph{moduli stack of formal groups}, so the functor sending an $\E_\infty$ ring to the space $\FG(\A)^\simeq$ of formal groups over $\A$. Let $\M_\FG^\ori\colon \CAlg \to \Spc$ be \emph{moduli stack of oriented formal groups}, so the functor sending an $\E_\infty$ ring $\A$ to the space $\FG^\ori(\A)^\simeq$ of oriented formal groups over $\A$.
\end{mydef}

We will see shortly in \Cref{affinemappingin} that this presheaf satisfies fpqc descent.

This stack implicitly appears in \cite{ec2}, explicitly in \cite[Con.4.2]{luriestheorem}, and plays a leading role in work of Gregoric \cite{rokmup,rokfiltration,rokeven}.

Computing the points of $\M_\FG^\ori$ is quite formal; see also \cite[Rmk.2.2.5]{rokmup}.

\begin{prop}\label{affinemappingin}
    Let $\A$ be an $\E_\infty$ ring. Then we have an equivalence of spaces
    \[\FG^\ori(\A)^\simeq = \Map_{\Stk}(\Spec \A, \M_\FG^\ori) \simeq 
    \begin{cases}
        \ast    &       \text{if $\A$ is complex-periodic,} \\
        \varnothing &   \text{if $\A$ is not complex-periodic.}
    \end{cases}\]
    In particular, the unique oriented formal group associated to a complex-periodic $\E_\infty$ ring is the \emph{Quillen formal group} $\widehat{\G}^\QQ_\A$ of \cite[\textsection4.1.3]{ec2}.
\end{prop}

\begin{proof}
    By \cite[Pr.4.3.23]{ec2}, if an oriented formal group lives over $\A$, then $\A$ is necessarily complex-periodic. In particular, if $\A$ is not complex-periodic, then $\FG^\ori(\A)=\varnothing$. Conversely, suppose that $\A$ is complex-periodic. By \cite[Pr.4.3.21]{ec2}, the category of preoriented formal groups over $\A$ agrees with the slice category $\FG(\A)_{\widehat{\G}^\QQ_\A/}$. By \cite[Pr.4.3.23]{ec2}, the subcategory $\FG^\ori(\A)$ of oriented formal groups over $\A$ is identified with those $\widehat{\G}^\QQ_\A \to \widehat{\G}$ which are equivalences, in other words with the subcategory of initial objects of $\FG(\A)_{\widehat{\G}^\QQ_\A/}$. Therefore, $\FG^\ori(\A)$ is contractible and hence so is $\FG^\ori(\A)^\simeq$.
\end{proof}

We will shortly see that $\M_\FG^\ori$ is an object of $\Stk$ by identifying it with a descent stack, but at the moment we consider it just as an object of $\Fun(\CAlg,\Spc^{\leq -1})\subset\Fun(\CAlg,\widehat{\Spc})$.

\begin{remark}\label{rmk:homotopygroupsofcomplexperiodicrings}
    The arguments used to prove \Cref{affinemappingin} give another interpretation for the homotopy groups of complex-periodic $\E_\infty$ rings. Indeed, the Quillen formal group $\widehat{\G}^\QQ_\A$ is canonically oriented over $\A$, which by definition yields a canonical equivalence
    \[\omega_{\widehat{\G}^\QQ_\A} \simeq \A[-2]\]
    between the shift of $\A$ and the \emph{dualizing line} of $\widehat{\G}^\QQ_\A$; see \cite[Df.4.2.1]{ec2}. Taking homotopy groups, we then have
    \[\omega_{\widehat{\G}^\heartsuit_\A} \simeq \pi_0 \omega_{\widehat{\G}^\QQ_\A} \simeq \pi_0 \A[-2] \simeq \pi_2 \A,\]
    where $\widehat{\G}^\heartsuit_\A$ is the classical Quillen formal group over $\pi_0\A$ and $\omega_{\widehat{\G}^\heartsuit_\A}$ is its module of invariant differentials. Tensoring yields $\omega_{\widehat{\G}^\heartsuit_\A}^{\otimes k} \simeq \pi_{2k}\A$ for all integers $k$.
\end{remark}

\subsubsection{Identification with the descent stack of $\MUP$}
Recall that $\MUP$ denotes the $\E_\infty$ ring of periodic complex cobordism, defined as the Thom spectrum of the map $\mathrm{BU}\times \mathbb{Z} \to \mathrm{pic}(\Sph)$ of infinite loop spaces.\footnote{For the purposes of this article, it suffices to take any complex-periodic $\E_\infty$ ring $\A$ such that $\Spec \A \to \M_\FG^\ori$ is a cover. As the proof of \cref{thm:pointsofDMUP} shows, other candidates include Snaith's construction $\Sph[\mathrm{BU}][\be^{-1}]$, $\MU^{t S^1}$, or any $\E_\infty$ form of periodic complex bordism which is isomorphic to $\MUP$ as a ring in the homotopy category of spectra. We recall that the first alternative and $\MUP$ are not equivalent as $\E_\infty$ rings by the main theorem of \cite{Hahn20Exotic}.}

Following \Cref{def:descentstacks}, we write $\D_\MUP$ for the descent stack of the map of stacks $\Spec(\MUP) \to \Spec(\spherespectrum)$. To relate $\D_\MUP$ to $\M_\FG^\ori$, it suffices to compute its $\A$-points for each $\E_\infty$ ring $\A$.

\begin{theorem}\label{thm:pointsofDMUP}
    Let $\A$ be an $\E_\infty$ ring. Then we an equivalence of spaces
    \[\D_\MUP (\A) \simeq 
    \begin{cases}
        \ast    &       \text{if $\A$ is complex-periodic,} \\
        \varnothing &   \text{if $\A$ is not complex-periodic.}
    \end{cases}\]
\end{theorem}

To prove this theorem, we will use the following two lemmata, which are surely well-known to experts.

\begin{lemma}\label{complexperiodicisfpqclocal}
An $\E_\infty$ ring $\A$ is complex-periodic if and only if there exists a faithfully flat map of $\E_\infty$ rings $\A \to \B$ and an $\E_\infty$ map $\MUP \to \B$.
\end{lemma}

\begin{proof}
    If $\A$ is complex-periodic, then $\A \to \A \otimes \MUP$ is faithfully flat from the usual computations of $\MUP$-homology of complex-oriented ring spectra
    \[\pi_\ast \A \otimes \MUP \simeq \pi_\ast \A[u^\pm, b_1, b_2, \ldots];\]
    see \cite[Lm.II.4.5]{bluebook}. Conversely, suppose that a faithfully flat map of $\E_\infty$ rings $\A \to \B$ exists such that $\B$ admits a map from $\MUP$. Note that $\MUP$ is obviously complex periodic, and therefore so is $\B$. So $\B$ is weakly $2$-periodic, which we recall means that $\B[2]$ is a locally free $\B$-module of rank 1. The same is then true for $\A$ as this property can be checked after faithfully flat base change. To finish the proof, it suffices to show that $\A$ is complex-orientable, which follows from \Cref{complexorientationisfpqclocal} below.
\end{proof}

\begin{lemma}\label{complexorientationisfpqclocal}
    Let $f\colon \A \to \B$ be a map of homotopy ring spectra which induces an injection on $\pi_\ast$. Then $\A$ is complex-orientable if and only if $\B$ is complex-orientable.
\end{lemma}

\begin{proof}
    If the unit $\Sph \to \A$ factors through $(\Sigma^\infty B^2\Z) [-2]$, then so does the unit into $\B$. Conversely, let $\sfC$ be any homotopy commutative ring spectrum, and denote by $x_\sfC$ the generator of $H^2(\CP^\infty;\sfC)$, given as the image of the canonical generator $x$ of $H^2(\CP^\infty;\Z)$ under the unit map $\Z \to \pi_0 \sfC$. Consider the multiplicative $\pi_\ast \sfC$-linear Atiyah--Hirzebruch spectral sequence (AHSS)
    \[E_2 = H^\ast(\CP^\infty; \pi_\ast \sfC) \simeq (\pi_\ast\sfC)[x_\sfC] \implies \sfC^\ast(\CP^\infty).\]
    Recall that a homotopy commutative ring spectrum $\sfC$ is complex-orientable if and only if $x_\sfC$ is a permanent cycle in this spectral sequence. In particular, if $\B$ is complex-orientable, then $x_\B$ is a permanent cycle, and the AHSS above for $\sfC = \B$ supports no differentials. The map $f\colon \A \to \B$ induces a map of AHSSs which is an injection on $E_2$-pages such that $f(x_\A) = x_\B$. It follows that $x_\A$ is also a permanent cycle. Indeed, if there is an $r$ such that $d_r(x_\A) \neq 0$, then this would be the first differential in the AHSS for $\A$, hence $f$ would induce an injection on $E_r$-pages. The nonvanishing of $d_r(x_\A)$ would then imply
    \[0\neq f(d_r(x_\A)) = d_r(f(x_\A)) = d_r(x_\B),\]
    contradicting the fact that $\B$ is complex-orientable. In particular, $d_r(x_\A)=0$ for all $r\geq 2$, and $x_\A$ is a permanent cycle.
\end{proof}

\begin{proof}[Proof of \Cref{thm:pointsofDMUP}]
    By \Cref{lem:descentstackbasicproperties}(\ref{item:descent_0_truncated}), for an $\E_\infty$ ring $\A$, the space $\D_\MUP(\A)$ is either empty or contractible. By the affine case of \Cref{lem:descentstackbasicproperties}(\ref{item:descent_0_truncated}), we see that $\D_\MUP(\A)$ is nonempty if and only if there exists a commutative diagram of the form
    \[\begin{tikzcd}
        {\Spec \B}\ar[r, "p"]\ar[d]  &   {\Spec \A}\ar[d]    \\
        {\Spec \MUP}\ar[r]  &   {\Spec \Sph}
    \end{tikzcd}\]
    where $p$ is faithfully flat. By \Cref{complexperiodicisfpqclocal}, this occurs if and only if $\A$ is complex-periodic.
\end{proof}

This immediately leads to the following statement.

\begin{theorem}\label{identicationofMFGOR}
    There is a unique equivalence of stacks
    \[\D_\MUP \simeq \M_\FG^\ori.\]
\end{theorem}

\begin{proof}
    By \Cref{cor:bigsheaves}, it suffices to construct this equivalence in $\Fun(\CAlg, \widehat{\Spc})$. This is then a consequence of \Cref{thm:pointsofDMUP} and \Cref{affinemappingin}, as both stacks are the same subobject of the terminal presheaf.
\end{proof}

We will see another computation of $\D_\A$ in \Cref{pointsofleqvech}. Other examples also follow as in \Cref{identicationofMFGOR}.

\begin{example}\label{dMOP}
    Let $\MOP$ denote the periodic unoriented bordism $\E_\infty$ ring. Using Hopkins' theorem that a homotopy commutative ring spectrum of characteristic $2$ is equivalent to a direct sum of suspensions and desuspensions of $H \F_2$, a similar argument as in the proof of \cref{thm:pointsofDMUP} shows that
    \[
    \D_{\MOP}(\A) \simeq \begin{cases}\ast&\text{if $\A$ is weakly $1$-periodic and $2=0$ in $\pi_0 A$},\\
    \emptyset&\text{otherwise}.
    \end{cases}
    \]
    Let $\M_\FG^{2{-}\text{tor},\heartsuit}$ denote the moduli stack of $2$-torsion formal groups. Then as with \cref{underlyingclassicalstackofDMUP}, Quillen's theorem for unoriented bordism gives an identificiation
    \[
    \D_{\MOP}^\heartsuit\simeq\M_\FG^{2{-}\text{tor},\heartsuit}
    \]
    of classical stacks.
\end{example}

\begin{remark}
Say that a map $\A\to\B$ of $\E_\infty$ rings is \emph{$\pi_\ast$-(faithfully) flat} if the induced map $\pi_\ast\A\to\pi_\ast\B$ of graded rings is (faithfully) flat. Following \cref{def:fpcq_top}, this notion of flatness leads to the \emph{$\pi_\ast$-fpqc topology} on $\Aff$. We have learnt through Marius Nielsen that forthcoming work will show that the functor $\Mod_{(-)}\colon \CAlg \to \CAlg(\PrLst)$ satisfies $\pi_\ast$-fpqc descent, allowing the general setup of \cref{sec:sheaves} to be applied with the $\pi_\ast$-fpqc topology on $\Aff$ in place of the fpqc topology. An analogous argument to \cref{thm:pointsofDMUP} shows that, if we write $\tilde{\D}_\MU$ for the descent stack of $\MU$ with respect to the $\pi_\ast$-fpqc topology, then $\tilde{\D}_\MU(\A)$ is contractible if $\A$ is complex orientable and empty otherwise. Similarly, $\tilde{\D}_\MO(\A) \simeq \tilde{\D}_{H\F_2}(\A)$ is contractible if $2=0$ in $\pi_0\A$ and empty otherwise.
\end{remark}

\subsubsection{Basic properties}
By identifying the descent stack of $\MUP$ with $\M_\FG^\ori$, the latter inherits many pleasant properties of the former. Many of the results here concerning $\M_\FG^\ori$ also appear in \cite{rokmup}.

\begin{cor}\label{coverisaffine}
    The unique map $\Spec \MUP \to \M_\FG^\ori$ is a faithfully flat map between geometric stacks. In particular, the natural map of stacks
    \[\colim \Spec \MUP^{\otimes (\bullet+1)} \to \M_\FG^\ori\]
    is an equivalence.
\end{cor}

\begin{proof}
    First, we notice that $\Spec \MUP \to \M_\FG^\ori$ is affine, as by \Cref{lem:descentstackbasicproperties}, the diagonal map $\M_\FG^\ori \to \M_\FG^\ori \times \M_\FG^\ori$ is affine; in fact, it is an equivalence. The first claim follows from \Cref{localpropertiesofDpandADAMSflatness} as $\Sph \to \MUP$ is Adams flat, while the second follows from \cref{lem:descentstackbasicproperties}.
\end{proof}

\begin{cor}\label{globalsectionofMFOR}
    The unit map $\Sph \to \Ga(\M_\FG^\ori)$ is an equivalence of $\E_\infty$ rings.
\end{cor}

\begin{proof}
    By \Cref{coverisaffine}, we see that $\colim \Spec \MUP^{\otimes(\bullet+1)}$ is a presentation of $\M_\FG^\ori$. In particular, we can compute $\Ga(\M_\FG^\ori)$ as the limit
    \[\Ga(\M_\FG^\ori) \simeq \lim \MUP^{\otimes(\bullet+1)},\]
    which is naturally equivalent to $\Sph$ as the sphere is $\MUP$-nilpotent complete \cite[Th.3.1]{bousfieldlocalisationexists}.
\end{proof}

Finally, recall that by Quillen's theorem the classical moduli stack of formal groups $\M_\FG^\heartsuit$ is equivalent to the stack associated to the Hopf algebroid $(\pi_0 \MUP, \pi_0 \MUP\otimes \MUP)$.

\begin{cor}\label{underlyingclassicalstackofDMUP}
    There are equivalences of classical stacks $\D_\MUP^\heartsuit \simeq (\M_\FG^\ori)^\heartsuit \simeq \M_\FG^\heartsuit$.
\end{cor}

\begin{proof}
    Although the map $\Sph \to \MUP$ is not flat, both of the units $\MUP \to \MUP \otimes \MUP$ are. In particular, the canonical map $f \colon \Spec \MUP \to \Spec \Sph$ is \emph{Adams flat}, as in \Cref{rmk:adamsflatness}, and this remark shows that $\D_f^\heartsuit$ is the groupoid object associated with the pair $(\Spec \pi_0 \MUP, \Spec \pi_0 \MUP \otimes \MUP)$. By Quillen's theorem, this is precisely the presentation of the classical moduli stack of formal groups; see for example \cite[Ex.3.1.19]{peterson_chromaticbook}.
\end{proof}

Both of the following results are a consequence of \Cref{lem:descentstackbasicproperties} and the identification $\M_\FG^\ori \simeq \D_\MUP$.

\begin{cor}\label{forgettingdataisfaithful}
    The forgetful functor $\Stk_{/\M_\FG^\ori} \to \Stk$ is fully faithful with essential image precisely those complex-periodic stacks.
\end{cor}

\begin{cor}\label{cor:complexperiodification}
    Let $\Aff^\cp$ denote the subcategory of $\Aff$ spanned by those complex-periodic affine stacks and write $\Stk^\cp$ for the category of stacks on $\Aff^\cp$ equipped with the fpqc-topology. Then the inclusion $\Aff^\cp \to \Aff$ induces a fully faithful functor $\Stk^\cp \to \Stk$ with essential image precisely the complex-periodic stacks. In particular, there is a unique identification $\Stk^\cp \simeq \Stk_{/\M_\FG^\ori}$ as subcategories of $\Stk$.
\end{cor}

\subsubsection{Fiber products}
Recall \Cref{lem:descentstackbasicproperties}(\ref{item:pbdescentstack}), which states that for a map $f\colon \Y \to \X$, and stacks $\sfZ_0,\sfZ_1$ defined over $\D_f$, then the natural map
\[\sfZ_0 \times_{\D_f} \sfZ_1 \to \sfZ_0 \times_\X \sfZ_1\]
is an equivalence. Combined with the equivalence $\M_\FG^\ori \simeq \D_\MUP$ by \Cref{identicationofMFGOR}, this immediately implies the following.

\begin{cor}\label{fiberofaffine}
    Given a span of stacks $\X \to \M_\FG^\ori \gets \Y$, then the natural map
\[\X \times \Y \to \X \times_{\M_\FG^\ori} \Y\]
    is an equivalence. In particular, 
    $\M_\FG^\ori \simeq \M_\FG^\ori \times \M_\FG^\ori$.
\end{cor}

There are many classical instances of this observation; see \cite[Pr.2.4]{akhilandlennart} or \cite[Pr.6.1]{tmfbook}, for example. 

The following is a further variation of this fact which we will need later.

\begin{lemma}
    Let $p\colon \Y \to \X$ be a map of stacks and $f\colon \sfZ \to \D_p$ be a map into the associated descent stack. If the composite $h = i f\colon \sfZ \to \X$ is universally $0$-affine, then the counit
    \[i^\ast i_\ast f_\ast (-) \xrightarrow{\simeq} f_\ast(-) \colon \QCoh(\sfZ) \to \QCoh(\D_p)\]
    is an equivalence.
\end{lemma}

This argument holds in any quasi-coherent sheaf context $(\csite, \tau,\sfQ)$.

\begin{proof}
    To show that the counit $i^\ast i_\ast f_\ast (-) \to f_\ast(-)$ is an equivalence, it suffices to pullback along the cover $q\colon \Y \to \D_p$. By \Cref{lem:descentstackbasicproperties}(\ref{item:pbdescentstack}), we have a Cartesian diagram of stacks
    \[\begin{tikzcd}
        {\sfZ\times_\X \Y}\ar[r, "{f'\simeq h'}"]\ar[d, "{q'\simeq p'}"]    &   {\Y}\ar[d, "q"] \\
        {\sfZ}\ar[r, "f"]                                                   &   {\D_p,}
    \end{tikzcd}\]
    which also identifies the pullbacks of morphisms $f'\simeq h'$ and $q'\simeq p'$. Notice that as $h$ is universally $0$-affine, then so is $f$ by \Cref{lem:descentstackbasicproperties}(\ref{item:pbdescentstack}). In particular, as both $f$ and $h$ are universally $0$-affine, part 1 of \Cref{prop:univ0affinefacts} implies that the above Cartesian square is adjointable. We then have the desired natural equivalences
    \[q^\ast i^\ast i_\ast f_\ast (-) \simeq p^\ast h_\ast (-) \simeq h'_\ast p'^\ast(-) \simeq f'_\ast q'^\ast(-) \simeq q^\ast f_\ast(-),\]
    the first by definition, the second by adjointability, the third from the identifications of pullbacks of morphisms above, and the fourth by adjointability again.
\end{proof}

This immediately implies the following.

\begin{cor}\label{relativespectrumfiberofaffine}
    Let $\A$ be a complex-periodic $\E_\infty$ ring and $\X$ a complex-periodic stack. Then there is a natural equivalence of quasi-coherent sheaves on $\X$ of the form
    \[\O_\X \otimes \A \simeq g^\ast f_\ast \O_{\Spec \A},\]
    where $f\colon \Spec \A \to \M_\FG^\ori$ and $g\colon \X\to \M_\FG^\ori$ are the canonical maps.
\end{cor}

\begin{proof}
    Writing $i\colon \M_\FG^\ori \to \Spec \Sph$ for the unique map, for each $M\in \Mod_\A$ we have a natural equivalence
    \[i^\ast i_\ast f_\ast(M) \simeq \O_{\M_\FG^\ori} \otimes M.\]
    Moreover, as $g^\ast$ is symmetric monoidal and colimit preserving, we also have the desired natural equivalences
    \[g^\ast f_\ast \O_{\Spec \A} \simeq g^\ast i^\ast i_\ast f_\ast \A \simeq g^\ast(\O_{\M_\FG^\ori} \otimes \A) \simeq \O_\X \otimes \A.\qedhere\]
\end{proof}

\subsubsection{Relative properties of complex-periodic stacks}\label{ssec:relativeqaffineness}
Each complex-periodic stack $\X$ admits a unique map $\X \to \M_\FG^\ori$, so asking this map to have various properties is a property of the stack $\X$ itself. 

\begin{mydef}
    Let $\P$ be a property of morphisms of stacks. A stack $\X$ is said to be \emph{chromatically $\P$} if $\X$ is complex-periodic and the associated map $\X \to \M_\FG^\ori$ is $\P$.
\end{mydef}

This applies to affine, $0$-semiaffine, quasi-affine, and universally $0$-affine, for example.

By \Cref{prop:univ0affinefacts} and \Cref{ex:quasiaffineareU0A}, we have the implications
\[
\text{chromatically affine}\implies \text{chromatically quasi-affine}\implies \text{chromatically universally $0$-affine}.
\]

The following is a useful criterion to check if a complex-periodic stack is $0$-affine.

\begin{prop}\label{cocontinuousandaffineis0affine}
    Let $\X$ be a chromatically affine stack. Then $\X$ is $0$-affine if and only if $\X$ is $0$-semiaffine, ie, the global sections functor $\Ga\colon \QCoh(\X) \to \Sp$ commutes with colimits.
\end{prop}

For example, quasi-compact separated schemes often have a cocontinuous global sections functor; see \cite[Cor.3.25]{akhilandlennart} or \Cref{ex:chromsemiafinne}.

\begin{proof}
    Part 3 of \Cref{prop:monoidalbarrbeck} asserts that a stack $\X$ is $0$-affine if and only if $\Ga\colon \QCoh(\X)\to\Sp$ is conservative and commutes with colimits. It therefore suffices to show that if $\X$ is chromatically affine and $\Ga$ commutes with colimits, then $\Ga$ is conservative. Let $\calM\in \QCoh(\X)$ be such that $\Ga(\calM)=0$. As $\Ga$ commutes with colimits, it follows that $\Ga(\calM)\otimes \MUP\simeq \Ga(\calM\otimes \MUP)$ also vanishes. Consider now the canonical Cartesian diagram of stacks
    \[\begin{tikzcd}
        {\X_\MUP}\ar[r, "{f'}"]\ar[d, "{g'}"]   &   {\Spec \MUP}\ar[d, "g"] \\
        {\X}\ar[r, "f"]                         &   {\M_\FG^\ori}
    \end{tikzcd}\]
    where all maps are affine, and $g$ and $g'$ are both faithfully flat. We can now further rewrite $\calM\otimes \MUP$ as
    \[\calM\otimes \MUP \simeq \calM \otimes_{\O_\X} (\O_\X\otimes \MUP)\simeq \calM\otimes_{\O_\X} (f^\ast g_\ast \MUP) \simeq \calM \otimes_{\O_\X} g'_\ast \O_{\X_\MUP} \simeq g'_\ast g'^\ast \calM\]
    courtesy of \Cref{relativespectrumfiberofaffine}, \Cref{prop:univ0affinefacts}, and the projection formula, respectively. As $f'$ is affine, $\X_\MUP$ is affine and hence the global sections functor $\Ga\colon \QCoh(\X_{\MUP}) \to \Sp$ is conservative. So the vanishing of $\Ga(g'_\ast(g'^\ast \calM)) \simeq \Ga(\calM \otimes \MUP)$ implies that $g'^\ast \calM$ vanishes. However, $g'$ is faithfully flat, hence $\calM$ also vanishes, as desired.
\end{proof}

\subsection{Examples of chromatically (quasi-) affine stacks}\label{ssec:classicalmoduliofformalgroups}
In this subsection, we demonstrate some techniques for proving that complex-periodic stacks are flat or (quasi-) affine over $\M_\FG^\ori$. At the end, we give a selection of examples. The reader who is willing to take these examples for granted can skip this subsection.

\subsubsection{Landweber exact stacks}\label{sssec:left}
Classically, we would say that a ring $R$ equipped with a formal group $\widehat{\G}$ is \emph{Landweber exact} if the associated map $\Spec R \to \M_\FG$ is flat; see \cite{naumann07} for the translation between this condition and the classical Landweber criterion of \cite{landweberoriginal}. We now discuss an analogous notion intrinsic to spectral algebraic geometry. 

\begin{mydef}\label{def:leftaffine}
    A \emph{Landweber exact stack} $\X$ is a complex-periodic geometric stack such that the unique map $\X \to \M_\FG^\ori$ is a flat map between geometric stacks. In particular, an $\E_\infty$ ring $\A$ is \emph{Landweber exact} if $\Spec \A$ is a Landweber exact stack.
\end{mydef}

Landweber exact stacks have a close relationship to their underlying stack.

\begin{cor}\label{equivdefniitionsofrelativelyquasiaffine}
    Let $\X$ be a nonconnective spectral Deligne--Mumford stack with affine diagonal which we further assume to be Landweber exact. Then $\X$ is chromatically affine (resp.\ quasi-affine) if and only if the associated map of classical stacks $f^\heartsuit\colon \X^\heartsuit \to \M_\FG^\heartsuit$ is affine (resp.\ quasi-affine).
\end{cor}

\begin{proof}
    This follows from \Cref{pr:heartsofquasiaffinenadaffine} and \Cref{cor:quasiaffinenessfromtheheart}.
\end{proof}

There is a close relationship between this spectral notion and the classical notion of Landweber exactness discussed above.

\begin{lemma}\label{classicalvsspectralpullbackv1}
    Let $\A$ and $\B$ be complex-periodic $\E_\infty$ rings where $\A$ is Landweber exact. Then the commutative diagram of classical sheaves
    \[\begin{tikzcd}
        {\Spec \pi_0 (\A \otimes \B)}\ar[d]\ar[r]   &   {\Spec \pi_0\A}\ar[d]  \\
        {\Spec \pi_0 \B}\ar[r]   &   {\M_\FG^\heartsuit}
    \end{tikzcd}\]
    is Cartesian. Moreover, the natural map of graded $\pi_\ast \B$-algebras
    \[\pi_\ast \B \otimes_{\pi_0 \B} \pi_0 (\A\otimes \B) \xrightarrow{\simeq} \pi_\ast (\A \otimes \B)\]
    is an isomorphism.
\end{lemma}

This will often be applied in the case where $\A=\MUP$. Note that $\MUP$ is Landweber exact by \cref{localpropertiesofDpandADAMSflatness}, as $\Spec(\MUP)\to \Spec(\spherespectrum)$ is Adams flat.

\begin{proof}
	First, notice that by \Cref{fiberofaffine}, the commutative diagram of stacks
	\[\begin{tikzcd}
		{\Spec \A \otimes \B}\ar[r]\ar[d]	&	{\Spec \A}\ar[d]	\\
		{\Spec \B}\ar[r]					&	{\M_\FG^\ori}
	\end{tikzcd}\]
	is Cartesian. In particular, the map $\Spec \B \to \M_\FG^\ori$ is affine. As the map $\Spec \A \to \M_\FG^\ori$ is flat as a map of geometric stacks by assumption, we can apply \Cref{flatnessislovely} to see that the desired diagram of classical stacks is Cartesian. For the ``moreover'' statement, we use the Cartesian diagram of stacks above together to see that the map $\B \to \A \otimes \B$ is flat by base change, from which the desired statement immediately follows.
\end{proof}

\begin{lemma}\label{eqdefLEFTaffine}
    A complex-periodic $\E_\infty$ ring $\A$ is Landweber exact if and only if the classical Quillen formal group $\widehat{\G}^\QQ_{\pi_0 \A}$ over $\pi_0 \A$ is Landweber exact.
\end{lemma}

\begin{proof}
    If $\Spec \A \to \M_\FG^\ori$ is flat, then by definition the base change $\MUP \to \MUP\otimes \A$ is flat. To see that $\Spec \pi_0 \A \to \M_\FG^\heartsuit$ is flat, it suffices to work after base change along the fpqc cover $\Spec \pi_0 \MUP \to \M_\FG^\heartsuit$. In this case, the classical pullback is the classical scheme $\Spec \pi_0(\MUP \otimes \A)$ by \Cref{classicalvsspectralpullbackv1}. As the map $\MUP \to \MUP\otimes \A$ was assumed to be flat, we see that $\pi_0 \MUP \to \pi_0(\MUP \otimes \A)$ is flat, as desired.

    Conversely, suppose that $\Spec \pi_0 \A \to \M_\FG^\heartsuit$ is flat. To see that $\Spec \A \to \M_\FG^\ori$ is flat, we need to show that its base change along the fpqc cover $\Spec \MUP \to \M_\FG^\ori$ is flat. By \Cref{fiberofaffine}, this base change is of the form $\Spec \MUP\otimes \A \to \Spec \MUP$, so it suffices to see that the map of $\E_\infty$ rings $\MUP \to \MUP \otimes \A$ is flat. By \Cref{classicalvsspectralpullbackv1}, the map $\pi_0 \MUP \to \pi_0 \MUP \otimes \A$ is flat, as this is the base change of a flat map. All that is left is the base change condition on higher homotopy groups. Using the classical computations of the $\MUP$-homology of complex-oriented ring spectra from \cite[\textsection II]{bluebook}, although here we are using $\MUP$ as opposed to $\MU$, we also obtain the isomorphisms
    \[\pi_\ast \MUP \otimes_{\pi_0 \MUP} \pi_0 (\MUP \otimes \A) \simeq (\pi_0 \MUP \otimes \A)[\be^\pm] \simeq \pi_\ast (\MUP \otimes \A)\]
    which shows that $\MUP \to \MUP\otimes \A$ is flat.
\end{proof}

This leads to an alternative characterization of Landweber exact stacks.

\begin{cor}\label{eqdefofLEFT}
    Let $\X$ be a complex-periodic geometric stack. Then the following are equivalent:
    \begin{enumerate}
        \item $\X$ is Landweber exact.
        \item For every affine and flat map $\Spec \A \to \X$, the complex-periodic $\E_\infty$ ring $\A$ is Landweber.
        \item The map of classical sheaves $\X^\heartsuit \to \M_\FG^\heartsuit$ defined by the classical Quillen formal group $\widehat{\G}_{\X^\heartsuit}^\QQ$ is flat.
    \end{enumerate}
\end{cor}

\begin{proof}
    Parts 1 and 2 are equivalent by \Cref{cor:invariantofflatness}, which shows that flatness can be checked on any atlas for $\X$. Given our flat presentation $\X\simeq \colim \Spec \A$ of $\X$, consider the following diagram of stacks
    \[\begin{tikzcd}
        {\Spec \MUP\otimes \A}\ar[r]\ar[d]  &   {\X_{\MUP}}\ar[r, "{f'}"]\ar[d]   &   {\Spec \MUP}\ar[d]  \\
        {\Spec \A}\ar[r]                    &   {\X}\ar[r, "f"]                 &   {\M_\FG^\ori,}
    \end{tikzcd}\]
    where the right-hand square is Cartesian by definition and the rectangle is Cartesian by \Cref{fiberofaffine}, hence the left-hand square is also Cartesian. The lower-horizontal composite is flat if and only if the underlying map of classical stacks $\Spec \pi_0 \A \to \M_\FG^\heartsuit$ is flat by \Cref{eqdefLEFTaffine}, and this latter condition is equivalent to the underlying map of classical stacks $\X^\heartsuit \to \M_\FG^\heartsuit$ being flat by \Cref{flatnessislovely} and descent. Hence, parts 2 and 3 are also equivalent.
\end{proof}

We will see examples of Landweber exact stacks in \Cref{ssec:examplesofLEFTandQAFF}.

\begin{nonexample}
    Consider $\Z^{tS^1}$, the $S^1$-Tate construction on $\Z$. This is an $\E_\infty$ ring with homotopy groups $\Z[t^\pm]$ for $|t| = -2$. This is clearly a complex-periodic $\E_\infty$ ring, it is an $\E_\infty$ $\Z$-algebra, so the associated formal group law is the additive formal group law. The associated map of classical stacks $\Spec \Z \to \M_\FG^\heartsuit$ is \textbf{not} flat. This is easily seen from the Landweber criterion, as multiplication by $v_0=p$ on $\Z$ has quotient $\mathbf{F}_p$, and the $p$-series of the additive formal group law over $\mathbf{F}_p$ vanishes, hence $(v_0,v_1,\ldots)$ is not a regular sequence.
\end{nonexample}

\subsubsection{Comparison to even-periodic refinements}
Although we do not find the definition \emph{Landweber exact} for spectral stacks in the literature, it has been a staple of chromatic homotopy theory for many years under a different guise.

\begin{mydef}
    Following \cite[Df.2.3]{akhilandlennart}, an $\E_\infty$ ring $\A$ is \emph{even-periodic}\footnote{This differs slightly to Lurie's definition of even-periodic from \cite[Ex.4.1.9]{ec2}, where he assumes that $\pi_\ast\A\simeq \pi_0 \A[t^\pm]$ with $|t|=2$. His definition clearly implies ours.} if $\A$ is complex-periodic and $\pi_{2k+1}\A=0$ for all integers $k$.
\end{mydef}

\begin{remark}
We note that complex-periodic $\E_\infty$ rings are not necessarily even-periodic; take $\KU^{S^1_+}$ or $\F_2^{tC_2}$ as examples.
\end{remark}

\begin{remark}\label{rmk:mapsbetweenevenperiodicthings}
    The category of even-periodic $\E_\infty$ rings has some nice properties. For example, a map of even-periodic $\E_\infty$ rings $\A \to \B$ is an equivalence (resp.\ flat/\'{e}tale) if and only if the underlying map of classical rings $\pi_0 \A \to \pi_0 \B$ is an equivalence (resp.\ flat/\'{e}tale). Indeed, in all cases, it suffices to show that the natural map $\pi_n \A \otimes_{\pi_0 \A} \pi_0 \B \to \pi_n \B$ is an equivalence for each integer $n$. For odd integers $n$, both sides of this map vanish, so we are left to check this for even $n=2m$. In the other case, we identify the above map with
    \[\omega^{\otimes m}_{\widehat{\G}^\QQ_\A} \otimes_{\pi_0 \A} \pi_0 \B \xrightarrow{\simeq} \omega^{\otimes m}_{\widehat{\G}^\QQ_\B}\]
    using \Cref{rmk:homotopygroupsofcomplexperiodicrings}. This map is an equivalence by base change, as the classical Quillen formal group over $\B$ is precisely that over $\A$ base changed along \emph{any} map $\A \to \B$. This last fact is clear geometrically: any map $\Spec \B \to \Spec \A$ is uniquely a map over $\M_\FG^\ori$.
\end{remark}

In the following, we will assume that our classical Deligne--Mumford stacks have quasi-affine diagonal so that their represented functors are fpqc-sheaves; see \Cref{rmk:fpqcclassicalcheck}.

\begin{mydef}
    Let $\X_0$ be a classical Deligne--Mumford stack equipped with a flat map of classical sheaves $f\colon \X_0 \to \M_\FG^\heartsuit$. Recall from \cite[Df.2.5]{akhilandlennart} that an \emph{even-periodic refinement} is a sheaf of $\E_\infty$ rings $\O^\top$ on the small \'{e}tale site $\X_0^\et$ equipped with an isomorphism $\pi_0 \O^\top \simeq \O_{\X_0}$ such that when restricted to affines $\Spec \A \to \X_0^\heartsuit$, the $\E_\infty$ ring $\O^\top(\Spec \A)$ is even periodic and Landweber exact (\'{a} la \Cref{def:leftaffine} and \Cref{eqdefLEFTaffine}) and for each integer $k$, the sheaf $\pi_k \O^\top$ on $\X_0^\et$ is quasi-coherent.
\end{mydef}

This concept of an even-periodic refinement was also used as the basis for \cite{osyn}.

Recall from \cite[Pr.1.3.1.7]{sag}, that a sheaf of $\E_\infty$ rings $\O^\top$ on the site $\X_0^\et$ is equivalently a sheaf of $\E_\infty$ rings on the topos $\Shv(\X_0^\et)$.

\begin{prop}\label{directcomparisontoMATHEWMEIER}
    Let $\X_0$ be a classical Deligne--Mumford stack with affine diagonal equipped with a flat map of classical stacks $f\colon \X_0 \to \M_\FG^\heartsuit$, and let $\O^\top$ be an even-periodic refinement thereof. Then $(\Shv(\X_0^\et), \O^\top)$ defines a $1$-localic nonconnective spectral Deligne--Mumford stack $\X$ and the associated stack $h_{\X}$ in $\Stk$ is Landweber exact. Moreover, every Landweber exact complex-periodic stack represented by a $1$-localic nonconnective spectral Deligne--Mumford stack with affine diagonal whose structure sheaf has vanishing homotopy groups in odd degrees occurs from such a construction.
\end{prop}

The affine diagonal assumption assures that $(\Shv(\X_0^\et), \O^\top)$ defines a \emph{geometric stack} \'{a} la \Cref{def:affineandflattogether}.

\begin{proof}
    By \cite[Lm.B.2]{tmfwls}, the spectrally ringed $\infty$-topos $\X = (\Shv(\X_0^\et), \O^\top)$ is a $1$-localic nonconnective spectral Deligne--Mumford stack and all such $1$-localic nonconnective spectral Deligne--Mumford stacks arise in such a way. As $\O^\top$ is an even-periodic refinement, then we see that $\X$ admits a representation $\X \simeq \colim \Spet \A$ in $\SpDMncet$, where all of the maps $\Spet \A \to \X$ are \'{e}tale and each $\A$ is even-periodic. In particular, as $h_{(-)} \colon \SpDMncet \to \Stk$ preserves colimits by \Cref{prop:dmstacksintoflatstacks}, we also have $h_{\X} \simeq \colim \Spec \A$ in $\Stk$, which immediately shows that $h_{\X}$ is complex-periodic:
    \[\Map_{\Stk}(h_{\X}, \M_\FG^\ori) \xrightarrow{\simeq} \lim \Map_{\Stk}(\Spec \A^\bullet, \M_\FG^\ori) \simeq \lim \ast \simeq \ast.\]
    Also notice that $h_{\X}$ is a geometric stack. Indeed, it suffices to show that $h_{\X}$ has affine diagonal, as we have already have an expression $h_{\X} \simeq \Spec \A$. It is clear that the nonconnective spectral Deligne--Mumford stack $\X$ has affine diagonal, as given $\Spet \B \to \X\times \X$, the pullback over $\X$ is affine if and only if the underlying Deligne--Mumford stack is affine; see \cite[Cor.1.4.7.3]{sag}. It follows from \Cref{prop:qaffdiagonalmeansdescent} that $h_{\X}$ satisfies flat descent, from which \Cref{prop:dmstacksintoflatstacks} also shows that $h_{\X}$ has affine diagonal.

    We now use \Cref{eqdefofLEFT}, which says that the complex-periodic geometric stack $h_{\X}$ is Landweber exact if and only if the map of classical sheaves is flat, which is true by assumption. The ``moreover'' statement follows from \cite[Lm.B.2]{tmfwls}, which shows that every $1$-localic nonconnective spectral Deligne--Mumford stack occurs in such a way, and the discussion above.
\end{proof}

\subsubsection{Examples of Landweber exact and chromatically (quasi-) affine stacks}\label{ssec:examplesofLEFTandQAFF}
This subsection is focused on examples of Landweber exact and chromatically (quasi-) affine stacks.

Keep in mind that \Cref{directcomparisontoMATHEWMEIER,equivdefniitionsofrelativelyquasiaffine} shows that all \emph{even-periodic refinements} of \cite{akhilandlennart} define $1$-localic nonconnective spectral Deligne--Mumford stacks whose structure sheaf is Landweber exact. In particular, this shows that $\M_\Ell^\ori$ and $\Spec \KU / C_2$ are chromatically affine, for example. For the rest of this section, we will be interested in \emph{not} appealing to \Cref{directcomparisontoMATHEWMEIER} and showing how one can work purely in complex-periodic geometry.

\begin{example}\label{leftisleft}
    Any complex-periodic $\E_\infty$ ring $\A$, whose underlying complex-orientable spectrum satisfies the classical criterion of Landweber, see \cite[Th.2.6]{landweberoriginal} and \cite[\textsection5]{naumann07}, is Landweber exact by \Cref{eqdefLEFTaffine}. For example, $\KU$, $\MUP$, and all Lubin--Tate theories $E_h$ are Landweber exact in this way. Other examples can be constructed from \emph{Lurie's theorem} \cite{luriestheorem}, a vast generalization of the Goerss--Hopkins--Miller theorem producing an $\E_\infty$-structure on Lubin--Tate theories. 
\end{example}

Recall that $\KU$ has an $\E_\infty$-action of $C_2$ given by complex conjugation. Let us write $\M_\Tori^\ori = \Spec \KU/C_2$ for the \emph{moduli stack of oriented tori}; in \cite[\textsection A]{normsontmf} we will discuss a moduli description for this stack. 

\begin{example}\label{modulioftori}
By construction, $\M_\Tori^\ori$ is a nonconnective spectral Deligne--Mumford stack and admits a $2$-sheeted \'{e}tale cover by $\Spec \KU$, so it is clearly complex-periodic. As the classical map $\widehat{\G}_m\colon \Spec \Z \to \M_\FG^\heartsuit$ induced by this cover, defined by the multiplicative formal group over $\Z$, is flat by \Cref{leftisleft}, we see that $\M_\Tori^\ori$ is also Landweber exact by \Cref{eqdefofLEFT}. Moreover, by \cite[\textsection6.2]{akhilandlennart}, for each geometric point $x\colon \Spec \kappa \to \Spec \Z$, the action of $C_2$ on $x^\ast \widehat{\G}_m$ is faithful, we see that the classical map of stacks $\Spec \Z / C_2 \to \M_\FG^\heartsuit$ is affine. By \Cref{equivdefniitionsofrelativelyquasiaffine}, we see that $\M_\Tori^\ori$ is also chromatically affine. Note that $\Ga(\O_{\M_\Tori^\ori}) \simeq \KO$.
\end{example}

\begin{remark}\label{modulioftoriwithtateyyyy}
    The above also shows that $\Spec \KU\llbracket q\rrbracket/C_2$ and $\Spec \KU\llpar q \rrpar/C_2$ are also Landweber exact and chromatically affine; here, $C_2$ acts on $\KU$ and trivially on $q$, see \cite[\textsection6]{globaltate} for more discussion. In fact, the analogous statement holds for any Landweber exact $\E_\infty$ ring $E$ with an action of a finite group $G$ such that at each geometric point of $\pi_0 E$, the action on the pullback of the Quillen formal group over $E$ is faithful; see \cite[\textsection6.2]{akhilandlennart} for more examples and discussion.
\end{remark}

Let $\M_\Ell^\ori$ denote the \emph{moduli of oriented elliptic curves} of \cite[Df.7.2.9]{ec2}.

\begin{example}\label{periodicisaffine}
	By construction, there is a map of stacks $\M_\Ell^\ori \to \M_\FG^\ori$ induced by the natural transformation $\widehat{(-)} \colon \Ell^\ori(\A) \to \FG^\ori(\A)$ from the category of oriented elliptic curves over an $\E_\infty$ ring $\A$ to the category of oriented formal groups over $R$, sending an elliptic curve to its underlying formal group; see \cite[Pr.7.1.2]{ec2}. In particular, the stack $\M_\Ell^\ori$ is complex-periodic. Moreover, by construction, it is nonconnective spectral Deligne--Mumford; see \cite[Pr.7.2.10]{ec2}. It also has affine diagonal, as this is true classically for $\M_\Ell^\heartsuit$. By \Cref{eqdefofLEFT,equivdefniitionsofrelativelyquasiaffine}, we see that $\M_\Ell^\ori$ is both Landweber exact and chromatically affine, as this is true of the map of classical underlying stacks $\M_\Ell^\heartsuit \to \M_\FG^\heartsuit$.
\end{example}

Let $\overline{\M}_\Ell^\ori$ be the \emph{compactification of the moduli stack of oriented elliptic curves}, defined as the pushout of nonconnective spectral Deligne--Mumford stacks
\[\begin{tikzcd}
    {\Spec \KU\llpar q \rrpar/C_2}\ar[r]\ar[d, "{\T}"]  &   {\Spec \KU\llbracket q \rrbracket/C_2}\ar[d]    \\
    {\M_\Ell^\ori}\ar[r]                                &   {\overline{\M}_\Ell^\ori,}
\end{tikzcd}\]
where $\T$ classifies the Tate curve over $\KU\llpar q \rrpar$ of \cite[Th.A]{globaltate}. By definition, $\Ga(\O_{\overline{\M}_\Ell^\ori}) = \Tmf$ is the $\E_\infty$ ring of \emph{projective topological modular forms}.

\begin{example}\label{projectiveisquasiaffine}
	By construction, the stack $\overline{\M}_\Ell^\ori$ is nonconnective spectral Deligne--Mumford; see \cite[Th.B]{globaltate}. Moreover, the above expression as a pushout immediately shows it is complex-periodic, as all other stacks in the defining span are. As this span also covers $\overline{\M}_\Ell^\ori$, and the other stacks in this span are all Landweber exact by \Cref{modulioftoriwithtateyyyy} and \Cref{periodicisaffine}, we see that $\overline{\M}_\Ell^\ori$ admits a Landweber exact cover, hence is also Landweber exact. Finally, we appeal to \cite[Th.7.2(1)]{akhilandlennart}, which states that the classical map $\overline{\M}_\Ell^\heartsuit \to \M_\FG^\heartsuit$ is quasi-affine. By \Cref{equivdefniitionsofrelativelyquasiaffine}, this also shows that $\overline{\M}_\Ell^\ori$ is chromatically quasi-affine. 
\end{example}

\begin{warn}
    Keep in mind that $\overline{\M}_\Ell^\ori$ is \textbf{not} chromatically affine. Indeed, this follows from the classical computation that the pullback of $\overline{\M}_\Ell^\heartsuit \to \M_\FG^\heartsuit$ against the cover $\Spec \pi_0\MP \to \M_\FG^\heartsuit$ is the weighted projective stack $\mathbf{P}(4,6)$. From another perspective, if $\overline{\M}_\Ell^\ori$ were chromatically affine, then by \cite[Th.C]{osyn}, the descent spectral sequence for $\overline{\M}_\Ell^\ori$ would agree with the Adams--Novikov spectral sequence for $\Tmf$, which is decidably false; see \cite[\textsection1.2.1]{smfcomputation}.
\end{warn}

\begin{remark}
    One can also use the Goerss--Hopkins--Miller construction of $\overline{\M}_\Ell^\ori$ as the pair $(\overline{\M}_\Ell^\heartsuit, \O^\top)$ to see that this stack is Landweber exact. Indeed, it is Landweber exact as using $\O^\top$, any affine \'{e}tale cover of $\overline{\M}_\Ell^\heartsuit$ can be upgraded to an affine \'{e}tale cover of $(\overline{\M}_\Ell^\heartsuit, \O^\top)$ by even periodic Landweber exact $\E_\infty$ rings. In fact, by \cite[Th.A]{uniqueotop}, this is one of the defining features of $\O^\top$.
\end{remark}

\subsection{Complex-periodic stacks of bounded height}\label{ssec:complextheoriesofboundedheight}
Just as we studied the moduli stack of oriented formal groups $\M_\FG^\ori$ in \Cref{ssec:MFGOR}, here we recall the substack $\M_\FGvech^\ori$ of oriented formal groups of height bounded by some fixed height function $\vec{h}$.

\subsubsection{Moduli stack of formal groups of bounded height}
Recall that the classical moduli stack of formal groups $\M_\FG^\heartsuit$ comes equipped with a height filtration. In particular, for each prime $p$, the localized stack $\M_\FG^\heartsuit\times \Spec \Z_{(p)}$ has a filtration by open substacks $\M_\FGh^\heartsuit$ which classify formal groups of height $\leq h$. Our spectral definition will reflect this classical situation, except that we allow for more than one single prime at once.

\begin{mydef}
    A \emph{height function} is a function $\vec{h}\colon \P \to \N\cup \{\infty\}$ from the set of all prime numbers to the set of nonnegative integers including infinity. For a height functor $\vec{h}$, a complex-periodic $\E_\infty$ ring $\A$ has \emph{height $\leq h$} if for each prime $p$, the localization at $p$ of the underlying map of classical stacks $\Spec \pi_0 \A \to \M_\FG^\heartsuit$ factors through $\M_\FGhp^\heartsuit$. A complex-periodic stack $\X$ has \emph{height $\leq \vec{h}$} if it admits a presentation $\X \simeq \colim \Spec \A$ where each complex-periodic $\A$ has height $\leq \vec{h}$.
\end{mydef}

\begin{mydef}
    Fix a height function $\vec{h}$. Define the presheaf $\M_\FGvech^\ori$ as the subpresheaf of $\M_\FG^\ori$ spanned by those oriented formal groups over $\E_\infty$ rings, such that for each prime $p$, when localized at $p$ the underlying classical formal group has height $\leq \vec{h}(p)$. Call this the \emph{moduli stack of oriented formal groups of height $\leq \vec{h}$}. This subpresheaf of the stack $\M_\FG^\ori$ is easily seen to be a stack as well.
\end{mydef}

\begin{remark}
    When localized at a prime $p$, the stack $\M_\FGvech^\ori$ is precisely the definition of $\M_\FGhp^\ori$ appearing in \cite[Df.3.2.1]{rokfiltration}.
\end{remark}

By construction, we obtain the following computation of the points of $\M_\FGvech^\ori$. Say that a complex-periodic $\E_\infty$ ring $\A$ has height $\leq \vec{h}$ if for each prime $p$, the classical Quillen formal group over $\pi_0 \A_{(p)}$ has height $\leq \vec{h}(p)$.

\begin{cor}\label{pointsofmfgvech}
    Let $\A$ be an $\E_\infty$ ring and $\vec{h}$ a height function. Then we have an equivalence of spaces
    \[
\sfD_{E(\vec{h})}(\sfA)\simeq
    \begin{cases}
        \ast    &       \text{if $\A$ is complex-periodic of height $\leq \vec{h}$,} \\
        \varnothing &   \text{if $\A$ is not complex-periodic of height $\leq \vec{h}$.}
    \end{cases}
    \]
\end{cor}

\begin{mydef}
    Let $\vec{h}$ be a fixed height function. A complex-periodic stack $\sfX$ has \emph{height $\leq \vec{h}$} if the (necessarily unique) map $\sfX \to \M_\FG^\ori$ factors through $\M_\FGvech^\ori$. 
\end{mydef}

\subsubsection{Identification with a descent stack and local descendability}
Mirroring \Cref{identicationofMFGOR}, we now want to find a class of $\E_\infty$ rings $E$ such that $\D_E \simeq \M_\FGvech^\ori$.

\begin{prop}\label{descentstackofEVECH}
Fix a height function $\vec{h}$.
\begin{enumerate}
\item Let $f \colon \Spec R \to \M_{\FGvech}^{\heartsuit} \subseteq \M_{\FG}^\heartsuit$ be a flat map. Then $f$ is faithfully flat if and only if $R/I_{p,\vec{h}(p)} \neq 0$ for all primes $p$. Here, $I_{p,h} = (p,v_1,\ldots,v_{h-1})$ is the $h$th Landweber ideal at the prime $p$.
\item There exist Landweber exact even-periodic $\E_\infty$ rings $E(\vec{h})$ of height $\leq \vec{h}$ for which the associated map $\Spec(\pi_0 E(\vec{h})) \to \M_{\FGvech}^\heartsuit$ is faithfully flat.
\item For any choice of $\E_\infty$ ring as above and any other $\E_\infty$ ring $\sfA$, we have equivalences
\[
\sfD_{E(\vec{h})}(\sfA)\simeq
    \begin{cases}
        \ast    &       \text{if $\A$ is complex-periodic of height $\leq \vec{h}$,} \\
        \varnothing &   \text{if $\A$ is not complex-periodic of height $\leq \vec{h}$.}
    \end{cases}
    \]
\end{enumerate}
\end{prop}
\begin{proof}
(1)~~We must prove that if $g\colon \Spec T \to \M_{\FGvech}^\heartsuit$ is any map, then the projection
\[
\Spec R \times_{\M_{\FGvech}^\heartsuit} \Spec T \to \Spec T
\]
is faithfully flat. This is flat by assumption, so we must prove that it is a surjection. To that end, it suffices to prove that if $T = \kappa$ is an algebraically closed field then this map admits a section. As $R/I_{p,\vec{h}(p)} \neq 0$, this follows as formal group laws over $\kappa$ are completely classified by their height.

(2)~~For each prime $p$, let $E_{p,\vec{h}(p)}$ denote a Lubin--Tate spectrum of height $\vec{h}(p)$ at the prime $p$. Then
\[
E(\vec{h}) = \prod_p E_{p,\vec{h}(p)}
\]
satisfies Landweber's exactness criterion as the same is true for $E_{p,\vec{h}(p)}$ and
\[
\pi_\ast E(\vec{h})/I_{p,n} = \pi_\ast E_{p,n}/I_{p,n}
\]
for $n\geq 1$. By the same reasoning (1) implies that $\Spec(\pi_0 E(\vec{h})) \to \M_{\FGhp}^\heartsuit$ is faithfully flat.

(3)~~As $E(\vec{h})$ is Landweber exact and even-periodic, $E(\vec{h}) \to E(\vec{h})\otimes E(\vec{h})$ is faithfully flat, and thus \cref{lem:descentstackbasicproperties}(\ref{item:descentstackimagecriterion}) implies
\[
\sfD_{E(\vec{h})}(\sfA)\simeq
\begin{cases}
\ast & \sfA \to \sfA \otimes E(\vec{h})\text{ faithfully flat} \\
\emptyset & \text{otherwise}.
\end{cases}
\]
It therefore suffices to prove that $\sfA$ is complex periodic of height $\leq \vec{h}$ if and only if $\sfA \to \sfA\otimes E(\vec{h})$ is faithfully flat. If $\sfA \to \sfA \otimes E(\vec{h})$ is faithfully flat, then as $\sfA \otimes E(\vec{h})$ is complex periodic of height $ \leq \vec{h}$, the same is true of $\sfA$. Conversely, suppose $\sfA$ is complex periodic of height $ \leq \vec{h}$, then as $E(\vec{h})$ is Landweber exact we may identify
\[
\Spec(\pi_0(\sfA \otimes E(\vec{h}))) \simeq \Spec(\pi_0 \sfA)\times_{\M_{\FGhp}^{\heartsuit}} \Spec(\pi_0 E(\vec{h})).
\]
This is faithfully flat over $\Spec(\pi_0 \sfA)$ as $\Spec(\pi_0 E(\vec{h})) \to \M_{\FGhp}^{\heartsuit}$ is faithfully flat by assumption.
\end{proof}

\begin{cor}\label{pointsofleqvech}
	Fix a height functor $\vec{h}$ and an $\E_\infty$ ring $E(\vec{h})$ as in part 1 of \Cref{descentstackofEVECH}. Then there is a unique equivalence of $(-1)$-truncated stacks
    \[\D_{E(\vec{h})} \simeq \M_\FGvech^\ori.\]
\end{cor}

\begin{proof}
	By definition $\D_{E(\vec{h})}$ is $(-1)$-truncated and $\M_\FGvech^\ori$ is $(-1)$-truncated as a substack of $\M_\FG^\ori$ which is $(-1)$-truncated by \Cref{affinemappingin}. As both stacks classifying the same subcategory of $\CAlg$ by part 2 of \Cref{descentstackofEVECH} and the definition $\M_\FGvech^\ori$, they must be uniquely equivalent.
\end{proof}

\begin{mydef}\label{def:dvech}
    As \Cref{pointsofleqvech} shows that $\D_{E(\vec{h})}$ only depends on $\vec{h}$, we simply write
    \[\D_{E(\vec{h})} = \D_{\vec{h}}.\]
\end{mydef}

We then immediately obtain the following as a consequence of \Cref{lem:descentstackbasicproperties}.

\begin{cor}
    With $\vec{h}$ and $E(\vec{h})$ as in \Cref{descentstackofEVECH}, then the unique map $\Spec E(\vec{h}) \to \M_\FGvech^\ori$ is an affine faithfully flat cover.
\end{cor}

For the stack $\M_\FGvech^\ori$ to be well-behaved, we need a further restriction on $\vec{h}$.

\begin{mydef}
	A height function $\vec{h}$ is \emph{bounded} if there is a finite nonnegative integer $N$ such that $\vec{h}(p)\leq N$ for all primes $p$.
\end{mydef}

The following is then a version of \Cref{coverisaffine} for $\D_{\vec{h}}$; similar $p$-local versions appear throughout \cite{rokfiltration}.

\begin{theorem}\label{thm:qcohofDvech}
Let $\vec{h}$ be a bounded height function. Then the inclusion
\[
\M_\FGvech^\ori \rightarrowtail \Spec(\spherespectrum)
\]
is locally descendable, and induces an equivalence
\[
\QCoh(\M_\FGvech^\ori) \simeq L_{\vec{h}}\Sp.
\]
\end{theorem}
\begin{proof}
Let $E(\vec{h})$ be an $\E_\infty$ ring as in part 2 of \cref{descentstackofEVECH}. 
By \cref{thm:locallydescendable}, it suffices to prove that $E(\vec{h}) \in \Sp$ is locally descendable. This is a reformulation of the uniform horizontal vanishing lines present in the $E_n$-based Adams spectral sequence, see \cite[Lm.A.4.6]{balderrama2024total}.
\end{proof}

\subsection{Affineness in complex-periodic geometry}\label{ssec:affinenesssection}
If $f\colon \Y \to \X$ is locally descendable, then by \Cref{thm:locallydescendable}, we see that $\D_f \to \X$ is universally $0$-affine. In this situation, if $\sfZ \to \D_f$ is also universally $0$-affine, then the composite $\sfZ \to \X$ is $0$-affine by \Cref{prop:univ0affinefacts}. This kind of argument yields many examples of $0$-affine stacks in complex-periodic geometry. First, a quick definition.

\begin{mydef}
    A complex-periodic stack $\X$ has \emph{bounded height} if the unique map $\X \to \M_\FG^\ori$ factors through $\M_\FGvech^\ori$ for some bounded height function $\vec{h}$.
\end{mydef}

\begin{theorem}\label{thm:affinenessbounded}
    A chromatically universally $0$-affine stack of bounded height $\X$ is universally $0$-affine. In particular, the global sections functor
    \[\Ga\colon \QCoh(\X) \to \Mod_{\Ga(\X)}\]
    is an equivalence of symmetric monoidal categories.
\end{theorem}

\begin{proof}
    There exists a bounded height function $\vec{h}$ for which $f$ factors as $\X \to \D_{\vec{h}} \rightarrowtail \M_\FG^\ori$ by assumption, and as $\D_{\vec{h}} \rightarrowtail \M_\FG^\ori$ is monic and $\X$ is chromatically universally $0$-affine, it follows that $\X \to \D_{\vec{h}}$ is universally $0$-affine. As $\D_{\vec{h}}$ is universally $0$-affine by \cref{thm:qcohofDvech}, \cref{prop:univ0affinefacts}(3) then implies that $\X$ is universally $0$-affine.
\end{proof}

The following emphasizes that \Cref{thm:affinenessbounded} also captures the main theorem of \cite{akhilandlennart}.

\begin{cor}\label{cor:mathewmeierintext}
    Let $\X^\heartsuit$ be a separated and Noetherian classical Deligne--Mumford stack equipped with a quasi-affine and flat map $\X^\heartsuit \to \M_\FG^\heartsuit$ and an even-periodic refinement $\X = (\X^\heartsuit, \O^\top_\X)$. Then $\X$ is $0$-affine.
\end{cor}

\begin{proof}
    By \Cref{ex:quasiaffineareU0A,thm:affinenessbounded}, it suffices to show that $\X$ is a chromatically quasi-affine stack of bounded height $\X$. By \Cref{directcomparisontoMATHEWMEIER}, we see $\X$ is Landweber exact and \Cref{equivdefniitionsofrelativelyquasiaffine} shows that $\X$ is chromatically quasi-affine. To see that $\X$ has bounded height, we argue as in \cite[\textsection4.2]{akhilandlennart}. Recall that the collection of open substacks $\M_\FGvech^\heartsuit \subseteq \M_\FG^\heartsuit$ as the height function $\vec{h}$ varies has a complementary collection of closed substacks $\M_\FGvechgeq^\heartsuit$. As $\X^\heartsuit$ is Noetherian, then $\X^\heartsuit \times_{\M_\FG^\heartsuit} \M_\FGvechgeq^\heartsuit$ must be empty for some bounded height function $\vec{h}$, meaning $\X$ has height $\leq \vec{h}$.
\end{proof}

Finally, we end with some remarks concerning \emph{$0$-semiaffineness}, so when a map of complex-periodic stacks $f\colon \Y\to \X$ induces a cocontinuous pushforward $f_\ast$.

\begin{cor}\label{cor:chromaticallysemiaffine}
    Let $\X$ be a chromatically $0$-semiaffine stack of bounded height. Then $\X$ is $0$-semiaffine, ie, the global sections functor
    \[\Ga\colon \QCoh(\X) \to \Mod_{\Ga(\X)}\]
    preserves colimits.
\end{cor}

\begin{proof}
    By assumption, the map $\X \to \M_\FGvech^\ori$ is semiaffine for some bounded height function $\vec{h}$. Since $\M_\FGvech^\ori \to \Spec \Sph$ is universally $0$-affine by \Cref{thm:qcohofDvech} hence $0$-affine by part 2 of \Cref{prop:univ0affinefacts}, we see that the composite $\X \to \Spec \Sph$ is $0$-semiaffine.
\end{proof}

There are more practical versions of the corollary above given by Mathew--Meier; see \cite[Th.4.14]{akhilandlennart}, for example. Nonetheless, there is some utility to the above statement.

\begin{example}\label{ex:chromsemiafinne}
    Let $\X$ be a complex-periodic stack such that the unique map $\X \to \M_\FG^\ori$ is a relative \emph{qc scheme}, meaning that for each $\Spec \A \to \M_\FG^\ori$ the pullback $\X \times_{\M_\FG^\ori} \Spec \A$ represents a quasi-compact scheme. Then $\X$ is chromatically $0$-semiaffine. Indeed, we can test this against the universal map $\Spec E \to \M_\FGvech^\ori$ for an appropriate $\vec{h}$, in which case the pullback $\Y = \Spec E \times \X$ is represented by a qc scheme by assumption. The descent spectral sequence for the global sections of $\Y$ has a horizontal vanishing line as the scheme underlying $\Y$ has finite cohomological dimension, from which it follows that $\Y \to \Spec E$ induces a cocontinuous pushforward as argued in \cite[Cor.3.25]{akhilandlennart}.
\end{example}

\begin{example}
    An example for such a phenomenon is the \emph{universal oriented elliptic curve} $\mathsf{E}^\ori$ living over $\M_\Ell^\ori$. By definition, the map $\mathsf{E}^\ori \to \M_\Ell^\ori$ is itself represented by a qc scheme, an elliptic curve in particular. As $\M_\Ell^\ori \to \M_\FG^\ori$ is affine, see \Cref{periodicisaffine}, we see that $\mathsf{E}^\ori$ itself is chromatically $0$-semiaffine by \Cref{ex:chromsemiafinne}, and hence $0$-semiaffine by \Cref{cor:chromaticallysemiaffine}.
\end{example}

\subsection{Complex-periodification}\label{sssec:maori}
As a consequence of \Cref{cor:complexperiodification}, the \emph{complex-periodification} of a stack $\X$ can be computed as $\X \times\M_\FG^\ori$.

We will now study such stacks in more detail, focusing on the affine case.

\subsubsection{The complex-periodification of an $\E_\infty$ ring}\label{sssec:complexperiodicifcofeinftyrings}

\begin{mydef}\label{def:complexperiodification}
    For an $\E_\infty$ ring $\A$, write $\M_\A^\ori$ for the product stack
    \[\M_\A^\ori = \D_\MUP \times \Spec \A,\]
    the \emph{complex-periodification} of $\Spec \A$ or simply of $\A$. We write $\M_\A^\heartsuit = (\M_\A^\ori)^\heartsuit$ for the underlying classical stack.
\end{mydef}

\begin{remark}\label{rmk:complexperiodicareblindtoperiodification}
By definition, a complex-periodic $\E_\infty$ ring $\B$ cannot see the difference between $\Spec \A$ and $\M_\A^\ori$:
\[\Map_{\Stk}(\Spec \B, \M_\A^\ori) \simeq \Map_{\Stk}(\Spec \B, \Spec \A) \times \Map_{\Stk}(\Spec \B, \M_\FG^\ori) \simeq \Map_{\CAlg}(\A, \B).\]
\end{remark}

This complex-periodification construction was already considered in \cite{rokeven}, although in a slightly different set-up. In particular, Gregoric's results concerning the complex-periodification of the sphere is immediate from the definition above, and the cases of the $L_h$-local sphere, $\KO$, and $\TMF$, so the equivalences
\[\M_{L_h\Sph}^\ori \simeq \M_\FGh^\ori, \qquad \M_{\KO}^\ori \simeq \Spec \KU / C_2, \qquad \M^\ori_{\TMF} \simeq \M_\Ell^\ori\]
all follow from \Cref{maintheorem:reconstruction}; in fact, \Cref{pointsofleqvech} also implies that for any Landweber exact $\E_\infty$ ring $E$ of height $\leq \vec{h}$ some height function $\vec{h}$, then the complex-periodification of $E$-local sphere $L_E \Sph$ is precisely $\M_\FGvech^\ori$.

Some properties of this construction come right from the definition.

\begin{prop}\label{pr:obviouspropertiesofMAORi}
    Let $\A$ be an $\E_\infty$ ring. 
    \begin{enumerate}
        \item\label{item:MAOridescentstack} The projection $\M_\A^\ori \to \Spec\A$ is monic and realizes $\M_\A^\ori$ as the descent stack for the projection $\Spec\MUP \times\Spec\A\to\Spec\A$. It is an equivalence if and only if $\A$ is complex-periodic.
        \item\label{item:MAOrirelspec} The canonical map $\M_\A^\ori \to \M_\FG^\ori$ is affine.
        \item The natural map of $\E_\infty$ rings $\A \to \Ga(\O_{\M_\A^\ori})$ is $\MUP$-nilpotent completion, ie, it induces an equivalence $\A^\wedge_\MUP \simeq \Ga(\O_{\M_\A^\ori})$.
        \item\label{item:MAOrihopkins} If $\A_\ast\MP$ is even, then $(\A_0\MP,\A_0\MP^{\otimes 2})$ forms a Hopf algebraoid whose associated stack may be identified with $\M_\A^\heartsuit$.
    \end{enumerate}
\end{prop}

\begin{proof}
(1)~~This is a special case of \cref{lem:descentstackbasicproperties}(\ref{item:restrictdescentstack}).

(2)~~By definition, the canonical map $\M_\A^\ori \to \M_\FG^\ori$ is the product of $\Spec \A \to \Spec \S$ with $\M_\FG^\ori$, so this holds as affine morphisms are stable under base change.

(3)~~This follows from \cref{prop:qcohdescentstack}(\ref{item:descentstacknilpotentcompletion}).

(4)~~As $\A_0\MP$ is even, the map $\A\otimes\MP \to \A \otimes \MP \otimes \MP \simeq (\A\otimes\MP)\otimes_\A(\A\otimes\MP)$ is flat, and this is sufficient to endow the pair
\[
(\A_0\MP,\A_0\MP^{\otimes 2}) = (\pi_0(\A \otimes\MP),\pi_0((\A\otimes\MP)\otimes_\A(\A\otimes\MP)))
\]
with the structure of a Hopf algebroid. By definition, the stack $\M_\A^\heartsuit$ is given by
\[
\M_\A^\heartsuit \simeq \colim \Spec(\pi_0 (\MP^{\otimes\bullet+1}\otimes \A)) \simeq \colim \Spec(\pi_0((\MP\otimes\A)^{\otimes_A\bullet+1})),
\]
and flatness again guarantees that this is equivalent to the canonical presentation of the stack associated to the Hopf algebroid $(\pi_0(\A \otimes\MP),\pi_0((\A\otimes\MP)\otimes_\A(\A\otimes\MP)))$.
\end{proof}

\begin{remark}
    One can also identify the \emph{descent spectral sequence} DSS for $\M_\A^\ori$ with the ANSS for $\A$ using the cover of $\M_\A^\ori$ given by $\Spec \A \otimes \MUP^{\otimes(\bullet+1)}$. A discussion of such generalized DSSs is not within the scope of this article. The second-named author will explore this further in forthcoming work.
\end{remark}

\begin{remark}
    There are many references to \emph{Hopkins' stack} associated to a homotopy commutative ring spectrum $R$ in the literature; see \cite[Dfs.3.1.11-12]{peterson_chromaticbook}, \cite[\textsection1.6]{handbooktmf}, and \cite[\textsection5]{carrickdefect}, for example. It follows from \Cref{pr:obviouspropertiesofMAORi}(\ref{item:MAOrihopkins}) that for an $\E_\infty$ ring $\A$ such that $\A_\ast\MUP$ is even, the classicla stack $\M_\A^\heartsuit$ is precisely \emph{Hopkins' stack} associated with $\A$. In particular, $\M_\A^\ori$ is a natural spectral refinement of Hopkins' stack construction.
\end{remark}

The assignment sending $\A$ to $\M_\A^\ori$ is in fact adjoint to the global sections functor. This is a consequence of the following general fact about descent stacks.

\begin{prop}\label{pr:globalsectionshasleftadjoint}
    Let $\Y$ be a stack and $\D_\Y$ be the descent stack associated to the unique map $\Y \to \Spec \Sph$. Then the global sections functor $\Stk_{/\D_\Y}^\op \to \CAlg$ admits a left adjoint sending $\A$ to $\D_\Y \times \Spec \A$.
\end{prop}

\begin{proof}
    This is purely formal: let $\sfZ$ be a stack over $\D_\Y$ which admits a colimit presentation $\sfZ \simeq \colim \Spec \B$. Recall from \Cref{item:fullyfaithfulness} that we have natural equivalences
    \[\Map_{\Stk}(\sfZ, \D_\Y \times \Spec \A) \simeq \lim \Map_{\Stk}(\Spec \B, \D_\Y) \times \Map_{\Stk}(\Spec \B, \Spec \A)\]
    \[\simeq \lim \Map_{\CAlg}(\A, \B) \simeq \Map_{\CAlg}(\A, \Ga(\O_\sfZ)),\]
    where $\Map_{\Stk}(\Spec \B, \D_\Y)$ is contractible as $\Spec \B$ admits a map to $\D_\Y$, and the latter is a subobject of $\Spec \Sph$.
\end{proof}

\begin{cor}\label{cor:globalsectionshasmoiraasleftadjoint}
    The global sections functor $\Ga\colon \Stk^\op_{/\M_\FG^\ori} \to \CAlg$ has left adjoint given by sending $\A$ to $\M_\A^\ori$.
\end{cor}

This left adjoint is clearly not essentially surjective, as all complex-periodifications $\M_\A^\ori$ are affine over $\M_\FG^\ori$. It is, however, monoidal.

\begin{cor}\label{generalizationoffiber}
    Let $\A,\B$ be $\E_\infty$ rings. Then the natural map of stacks
    \[\M_{\A\otimes \B}^\ori \xrightarrow{\simeq} \M_\A^\ori \times_{\M_\FG^\ori} \M_\B^\ori \simeq \M_\A^\ori\times \M_\B^\ori\]
    is an equivalence.
\end{cor}

\begin{proof}
By construction, we have
\begin{align*}
\M_{\A\otimes\B}^{\ori}\simeq\Spec(\A\otimes\B)\times\M_\FG^\ori&\simeq\Spec\A\times\Spec\B\times\M_\FG^\ori\\
&\simeq\Spec\A\times\M_\FG^\ori\times\Spec\B \times\M_\FG^\ori\simeq\M_\A^\ori\times\M_\B^\ori,
\end{align*}
as $\M_\FG^\ori\simeq\M_\FG^\ori\times\M_\FG^\ori$.
\end{proof}

\subsubsection{Examples of complex-periodifications}\label{sssec:maoriexamples}
There are two particular examples of $\M_\A^\ori$ which we want to highlight; these examples produce complex-periodic stacks $\M_\A^\ori$ of \emph{infinite height}, as in both cases $\F_2 \otimes \A$ does not vanish, and hence do not fit into the collection of examples found in \Cref{maintheorem:boundedaffineness}.

Many more interesting examples of complex-periodifications also appear shortly; see \Cref{ssec:reconstructionsection}.

\begin{example}[The stack associated to $\ko$]\label{ex:koperiodification}
    Write $\ko = \tau_{\geq 0}\KO$ for the $\E_\infty$ ring of \emph{connective real $K$-theory}. The truncation map $\ko \to \Z$ shows that $\M_{\ko}^\ori$ does not have bounded height. The underlying classical stack of $\M_{\ko}^\ori$ can be identified with the \emph{moduli stack of quadratic curves}; see \cite[\textsection3.6]{tmfhomology} in the first arXiv version. In particular, $\M_{\ko}^\ori$ is a spectral lift of the moduli stack of quadratic curves.
\end{example}

\begin{example}[The stack associated to $\tmf$]\label{ex:tmfperiodification}
    Write $\tmf = \tau_{\geq 0} \Tmf$ for the $\E_\infty$ ring of \emph{connective topological modular forms}, where $\Tmf$ is the $\E_\infty$ ring of \emph{projective topological modular forms}; see \Cref{projectiveisquasiaffine} for a reminder. Again, there is a truncation map $\tmf \to \Z$ which shows that $\M_{\tmf}^\ori$ cannot have bounded height. One can identify the underlying classical stack of $\M_{\tmf}^\ori$ with the \emph{moduli stack of cubic curves}; see \cite[\textsection5.1]{tmfhomology}.
\end{example}

It is well-known that both $\ko$ and $\tmf$ are $\MU$-nilpotent, and we will see in \Cref{thm:0affinenessforMUPNILPOTENZ} that this implies that the stacks $\M_{\ko}^\ori$ and $\M_{\tmf}^\ori$ are both $0$-affine. They are both also (more or less tautologically) reconstructible; see \Cref{cor:easyreconstruction}.

By construction, the stacks $\M_\qua^\ori$ and $\M_\cub^\ori$ are affine over $\M_\FG^\ori$. However, they are \textbf{not} flat over $\M_\FG^\ori$, nor are either of them Deligne--Mumford stacks; see \cite[Rmk.4.9]{osyn}. In particular, the technology of Mathew--Meier does not apply to these stacks. Nevertheless, our general $0$-affineness results still apply to this situation.

\begin{remark}
    Let $E_h$ be the Lubin--Tate theory of height $h$ corresponding to the Honda formal group law over $\F_{p^h}$. Let us write $EO_h$ for the fixed-points of $E_h$ with respect to the maximal finite subgroup of the Morava stabilizer group. In \cite{gorbounovmahowald}, Gorbounov--Mahowald essentially study a Hopf algebroid model for $\M_{EO_{p-1}}^\heartsuit$ generalizing the Legendre family of elliptic curves, following a suggestion of Hopkins. Currently, it is not clear how to construct ``connective models'' for $eo_{p-1}$, or more generally $eo_h$, for $p\geq 5$; we have $eo_1 = \ko_2^\wedge$ and $eo_2 \simeq \tmf_3^\wedge$. In \cite{hill_eo4}, Hill assumes that such a connective spectrum $eo_4$ exists, and studies the associated Hopf algebroid presentation for $\M_{eo_4}^\heartsuit$. By \Cref{pr:obviouspropertiesofMAORi}(\ref{item:MAOrihopkins}), we see that when such spectra $eo_{p-1}$ exist as $\E_\infty$ rings, the corresponding stacks $\M_{eo_{p-1}}^\ori$ would be spectral lifts of a connective version of the Gorbounov--Mahowald stacks of \cite{gorbounovmahowald}.
\end{remark}

\subsubsection{0-affineness for complex-periodifications}\label{sssec:0affinenessforMAORI}
In general, quasi-coherent sheaves over $\M_\A^\ori$ are the $\QCoh(\M_\FG^\ori)$-linearisation of $\Mod_\A$, 
\begin{equation}
    \QCoh(\M_\FG^\ori) \otimes_{\Sp} \Mod_\A \xrightarrow{\simeq} \QCoh(\M_\A^\ori).
\end{equation}
This follows from \Cref{prop:adjointablesquare0affine,prop:univ0affinefacts}. In particular, $\QCoh(\M_\A^\ori)$ can be just as complicated as $\QCoh(\M_\FG^\ori)$. Additional assumptions on $\A$ simplify the situation.

\begin{mydef}
    A spectrum $X$ is \emph{$\MUP$-nilpotent} if $X$ lies in the thick tensor ideal of $\Sp$ generated by $\MUP$.\footnote{As $\MUP$ is a sum of shifts of $\MU$, the class of $\MUP$-nilpotent spectra agrees with the class of $\MU$-nilpotent spectra.}
\end{mydef}

\begin{example}
    Complex-orientable spectra $X$ are $\MUP$-nilpotent as they are homotopy $\MUP$-modules. By the thick subcategory theorem of \cite[Th.7]{hopkinssmith}, if there exists a type $0$ finite spectrum $F$ such that $X \otimes F$ is complex-orientable, then $X$ is $\MUP$-nilpotent. This includes the Wood-type spectra of \cite{carrickdefect} (see \cite[Pr.2.32]{carrickdefect}), and thus includes many nonperiodic examples such as $\ko$ and $\tmf$. Spectra such as the connective image-of-$J$ spectrum are $\MUP$-nilpotent but not of Wood-type; see \cite[Cor.5.19]{carrickdefect}.
\end{example}

\begin{remark}
    A spectrum $X$ being $\MUP$-nilpotent implies that the ANSS for $X$ has a horizontal vanishing line on some finite page; see \cite[Cor.4.4]{akhilgalois} and \cite[Pr.A.1.6]{balderrama2024total}.
\end{remark}

\begin{theorem}\label{thm:0affinenessforMUPNILPOTENZ}
    Let $\A$ be an $\MUP$-nilpotent $\E_\infty$ ring. Then $\M_\A^\ori \to \Spec \A$ is descendable. In particular, $\M_\A^\ori$ is $0$-affine with $\QCoh(\M_\A^\ori)\simeq\Mod_\A$.
\end{theorem}

\begin{proof}
As $\A$ is $\MUP$-nilpotent, $\A \to \A\otimes\MUP$ is descendable by \cref{lem:ldclosure}(4). As $\M_\A^\ori$ is equivalent to the descent stack for $\A \to \A\otimes\MUP$ by \Cref{pr:obviouspropertiesofMAORi}(\ref{item:MAOridescentstack}), \cref{thm:locallydescendable} implies that $\M_\A^\ori \to \Spec(\A)$ is descendable.
\end{proof}

In particular, we see that the stacks $\M_{\ko}^{\ori}$ and $\M_{\tmf}^{\ori}$ of \Cref{ex:koperiodification,ex:tmfperiodification} are $0$-affine.

\subsection{Reconstruction}\label{ssec:reconstructionsection}
Both \Cref{thm:0affinenessforMUPNILPOTENZ,thm:affinenessbounded} give us a wealth of $0$-affine complex-periodic stacks. In this section, we study a subcategory of such complex-periodic stacks which can be recovered from their global sections: we call such a phenomenon \emph{reconstruction}.

First, we define this phenomenon, before we give a collection of examples. We will often restrict to chromatically affine stacks $\X$. By \Cref{cocontinuousandaffineis0affine}, the stack $\X$ is $0$-affine if and only if it is $0$-semiaffine.

\subsubsection{Definition and basic properties}\label{sssec:basicreconstructionproeprtioes}
Recall the adjunction of \Cref{cor:globalsectionshasmoiraasleftadjoint}:
\[\Ga\colon \Stk_{/\M_\FG^\ori} \rightleftarrows \CAlg^\op \colon \M_{(-)}^\ori\]

\begin{mydef}
    A complex-periodic stack $\X$ is \emph{reconstructible} if the canonical unit map
    \[\X \to \M_{\Ga(\O_\X)}^\ori\]
    is an equivalence of stacks.
\end{mydef}

Our first result follows formally from the definition.

\begin{prop}\label{prop:ffforreconstruction}
Let $\X$ be a reconstructible complex-periodic stack. Then for any complex-periodic stack $\Y$, the map
\[
\Map_{\Stk}(\Y,\X) \to \Map_{\CAlg}(\Gamma(\O_\X),\Gamma(\O_\Y))
\]
is an equivalence.
\end{prop}
\begin{proof}
By the adjunction and equivalence $\M_{\Ga(\O_\X)}^\ori\simeq\X$, we have
\[
\Map_{\CAlg}(\Gamma(\O_\X),\Gamma(\O_\Y))\simeq \Map_{\Stk_{/\M_\FG^\ori}}(\Y,\M_{\Gamma(\O_\X)}^\ori)\simeq \Map_{\Stk_{/\M_\FG^\ori}}(\Y,\X)\simeq\Map_{\Stk}(\Y,\X)
\]
as $\Stk_{/\M_\FG^\ori}\subset\Stk$ is fully faithful.
\end{proof}

A basic supply of examples of reconstructible stacks comes from the following.

\begin{prop}\label{cor:easyreconstruction}
    Let $\A$ be an $\MUP$-nilpotent complete $\E_\infty$ ring. Then $\M_\A^\ori$ is reconstructible.
\end{prop}

For example, the $\MUP$-nilpotent spectra $\ko$ and $\tmf$ of \Cref{ex:koperiodification,ex:tmfperiodification} are $\MUP$-nilpotent complete.

\begin{proof}
    The global sections of $\M_{\A}^\ori$ are precisely the $\MUP$-nilpotent completion of $\A$ by part 3 of \Cref{pr:obviouspropertiesofMAORi}. If $\A \simeq \A^\wedge_{\MUP}$, then the counit map $\M_{\A^\wedge_{\MUP}}^\ori \to \M_\A^\ori$ is an equivalence, and we are done.
\end{proof}

This also shows that there are complex-periodic stacks which are reconstructible but not $0$-affine.

\begin{example}
    The stack $\M_\FG^\ori = \M_{\Sph}^\ori$ is reconstructible as the complex-periodification of the sphere spectrum $\Sph$, as $\Sph$ is $\MUP$-nilpotent complete. On the other hand, $\M_\FG^\ori$ is \textbf{not} $0$-affine, in fact, it is not even $0$-semiaffine. More generally, if $\A$ is an augmented $\E_\infty$ ring, meaning it admits a map of $\E_\infty$ rings $\A \to \Sph$, then $\M_\A^\ori$ is not $0$-affine. Indeed, the augmentation induces a map on affine stacks, which by base change along $\M_\FG^\ori$ induces an affine map $f\colon \M_\FG^\ori \simeq \M_\Sph^\ori \to \M_\A^\ori$. By \Cref{prop:univ0affinefacts}, the pushforward $f_\ast$ preserves colimits, so we write the global sections on $\QCoh(\M_\FG^\ori)$ as the composition
    \[\Ga\colon \QCoh(\M_\FG^\ori) \xrightarrow{f_\ast} \QCoh(\M_\A^\ori) \to \Sp.\]
    If $\M_\A^\ori$ was $0$-semiaffine, then $\M_\FG^\ori$ which be $0$-semiaffine, which we know is false. By \Cref{thm:0affinenessforMUPNILPOTENZ}, this implies that no augmented $\E_\infty$ rings can be $\MUP$-nilpotent.
\end{example}

\subsubsection{Examples of reconstructible stacks}\label{sssec:exampleofreconstruction}
The stacks $\M_\A^\ori$ are all chromatically affine, as the map $\M_\A^\ori \to \M_\FG^\ori$ is base changed from the affine map $\Spec \A \to \Spec \Sph$. Conversely, a stack $\X$ which is also $0$-semiaffine, is reconstructible as soon as it is also chromatically affine.

\begin{theorem}[{\Cref{maintheorem:reconstruction}}]\label{thm:reconstruction}
    Chromatically affine stacks $\X$ which are also $0$-semiaffine are reconstructible.
\end{theorem}

By \Cref{pr:obviouspropertiesofMAORi}(\ref{item:MAOridescentstack}) and \Cref{thm:0affinenessforMUPNILPOTENZ}, such $\X$ include $\M_\A^\ori$ for $\E_\infty$ rings $\A$ which are $\MUP$-nilpotent. More examples come from \Cref{thm:affinenessbounded} and \Cref{cor:chromaticallysemiaffine}. We list a few important instances here.

\begin{example}
    For any height function $\vec{h}$, the map of stacks $\M_\FGvech^\ori \to \M_\FG^\ori$ is clearly affine; this follows straight from the computations of their points from \Cref{affinemappingin,pointsofmfgvech}. If we further assume that $\vec{h}$ is bounded, then we also see that $\M_\FGvech^\ori$ is $0$-affine. By \Cref{thm:reconstruction}, we see obtain the first equivalence
    \[\M_\FGvech^\ori \simeq \M_{\Ga(\O_{\M_\FGvech^\ori})}^\ori \simeq \M_{L_{\vec{h}}\Sph}^\ori,\]
    the second equivalence comes from \Cref{thm:qcohofDvech}.
\end{example}

\begin{example}
    Writing $\M_\Tori^\ori = \Spec \KU / C_2$ for the moduli stack of oriented tori of \Cref{modulioftori}, then \emph{ibid} explicitly shows that the unique map $\M_\Tori^\ori \to \M_\FG^\ori$ is affine. Moreover, since tori always have associated formal group of height $1$, $\M_\Tori^\ori$ has bounded height and by \Cref{thm:affinenessbounded} we see this stack is also $0$-affine. By \Cref{thm:reconstruction}, we have a natural equivalence of stacks
    \[\M_\Tori^\ori \simeq \M_{\KO}^\ori.\]
    Similarly, we have equivalences
    \[\Spec \KU\llbracket q \rrbracket / C_2 \simeq \M_{\KO\llbracket q \rrbracket}^\ori, \qquad \Spec \KU\llpar q \rrpar / C_2 \simeq \M_{\KO\llpar q \rrpar}^\ori\]
    courtesy of \Cref{modulioftoriwithtateyyyy}.
\end{example}

\begin{example}
    By \Cref{periodicisaffine}, we see that the moduli stack of oriented elliptic curves $\M_\Ell^\ori$ is affine over $\M_\FG^\ori$. Moreover, as elliptic curves have associated formal groups of height $\leq 2$, we see that $\M_\Ell^\ori$ is a chromatically affine stack of bounded height. In particular, by \Cref{thm:affinenessbounded}, we see it is $0$-affine. By \Cref{thm:reconstruction}, we obtain the identification
    \[\M_\Ell^\ori \simeq \M_{\TMF}^\ori.\]
\end{example}

\begin{warn}
    Notice that the ``chromatically affine'' assumption in \Cref{thm:reconstruction} and subsequent corollaries cannot be easily removed. For example, the two complex-periodic stacks $\M_{\Tmf}^\ori$ and $\overline{\M}_\Ell^\ori$ both have global sections $\Tmf$, yet the former is chromatically affine and the latter is not. In particular, these two stacks are not equivalent, yet they admit equivalent global sections, violating any potential reconstruction statement. More generally, if the descent spectral sequence for $\X$ does not agree with the Adams--Novikov spectral sequence for $\Ga(\O_\X)$, then reconstruction cannot hold. This is a major point of \cite{smfcomputation,osyn}.
\end{warn}

\begin{example}\label{ex:cooperationalgebras}
    Reconstruction also combines with \Cref{generalizationoffiber} to produce stacks associated to various cooperation algebras and homologies. For example, the $\E_\infty$ ring of $\KO$-cooperations $\KO \otimes \KO$ is the global sections of $\M_\Tori^\ori \times \M_\Tori^\ori$, and $\TMF\otimes \TMF$ is the global sections of $\M_\Ell^\ori \times \M_\Ell^\ori$, and that these fiber products can also be taken over $\M_\FG^\ori$ as desired. Similarly, the $\TMF$-homology of $\KO$ is the global setions of $\M_\Ell^\ori \times \M_\Tori^\ori$. Using this description of cooperations, one can realize the maps
    \[\KO \otimes \KO[\tfrac{1}{n}] \to \prod_{k\geq 1} \KO[\tfrac{1}{n}], \qquad \TMF \otimes \TMF[\tfrac{1}{n}] \to \prod_{k\geq 1} \TMF[\tfrac{1}{n}]\]
    of \cite[Rmk.3.8 \& \textsection4.1]{boss}, given by linearising a product of stable Adams operations $\psi^{n^k}$, geometrically. Let us give the details in the $\TMF$-case. From the above arguments, it suffices to construct this as a map of complex-periodic stacks
    \[\coprod_{k\geq 1} \M_\Ell^\ori \to \M_\Ell^\ori \times_{\M_\FG^\ori} \M_\Ell^\ori,\]
    where we have implicitly inverted $n$. To do that, we produce the commutative diagram of complex-periodic stacks
    \[\begin{tikzcd}
        {\coprod_{k\geq 1} \M_\Ell^\ori}\ar[r, "f"]\ar[d, "g"]  &   {\M_\Ell^\ori}\ar[d]    \\
        {\M_\Ell^\ori}\ar[r]    &   {\M_\FG^\ori,}
    \end{tikzcd}\]
    where $f$ is the coproduct of the identity map, and $g$ is the coproduct of the map defined by the universal oriented elliptic curve equipped with the orientation given by applying the $n^k$-fold multiplication map to the universal orientation. These linearised Adams operations are not immediately clear to be equivalent to those of \cite{heckeontmf,adamsontmf}, although they are all built from the same $n^k$-fold multiplication map on underlying formal groups.
\end{example}

Our proof of \Cref{thm:reconstruction} is quite straight forward, by writing the complex-periodification $\M_{\Ga(\O_\X)}^\ori$ as the colimit of the diagram $\Spec \MUP^{\otimes(\bullet+1)} \otimes \Ga(\O_\X)$, and then passing $\Ga(-)$ over this tensor product to obtain an expression for $\X$. First, we will use a lemma.

\begin{lemma}\label{lm:globalsectionsarenilpotentcomplete}
    Chromatically affine stacks $\X$ which are also $0$-semiaffine are reconstructible and have $\MUP$-nilpotent complete global sections.
\end{lemma}

\begin{proof}
    For some complex-periodic $\E_\infty$ ring $\A$, consider the Cartesian diagram of stacks
    \begin{equation}\label{eq:pullbackinreconstruction}\begin{tikzcd}
        {\X_\A}\ar[d, "{f'}"]\ar[r, "{g'}"] &   {\X}\ar[d, "f"]   \\
        {\Spec \A}\ar[r, "g"] &   {\D_\MUP,}
    \end{tikzcd}\end{equation}
    where $\X_\A$ is affine as $f$ is affine. Setting $\A=\MUP^{\otimes(\bullet+1)}$ and $\X_\A = \X_\bullet \simeq \Spec \B^\bullet$, we have natural equivalences of $\E_\infty$ rings
    \[\lim \Ga(\O_\X) \otimes \MUP^{\otimes(\bullet+1)} \simeq \lim \Ga(\O_\X \otimes \MUP^{\otimes(\bullet+1)}) \simeq \lim \Ga( g^\ast f_\ast \MUP^{\otimes (\bullet+1)})\]
    \[\simeq \lim \Ga( f'_\ast \O_{\X_\bullet})\simeq \lim \B^\bullet \simeq \Ga(\O_\X),\]
    the first as $\X$ is $0$-affine, second from \Cref{relativespectrumfiberofaffine}, third from the Beck--Chevalley condition for affine maps of \Cref{prop:univ0affinefacts}, forth as $\X_\bullet \simeq \Spec \B^\bullet$ is affine, and fifth as $g'$ is faithfully flat and hence an atlas for $\X$.
\end{proof}

\begin{proof}[Proof of \Cref{thm:reconstruction}]
    Consider the Cartesian square of stacks
    \[\begin{tikzcd}
        {\X_\bullet}\ar[r, "g'"]\ar[d, "f'"]  &   {\X}\ar[d, "f"] \\
        {\Spec \MUP^{\otimes(\bullet+1)}}\ar[r, "g"]                 &   {\M_\FG^\ori}
    \end{tikzcd}\]
    for our given $\X$. As $f$ is affine by assumption, we see that $f'$, and hence also $\X_\bullet$, are both affine by base change. We then venture to compute $\X_\bullet = \Spec \Ga(\O_{\X_\bullet})$. The map $g$ is affine as $\M_\FG^\ori$ has affine diagonal by \Cref{fiberofaffine}, hence so is $g'$, so the square above is adjointable by \Cref{prop:univ0affinefacts}. Writing $h\colon \X \to \Spec \Sph$, we can then compute $\Ga(\O_{\X_\bullet})$ as
    \[\Ga(\O_{\X_\bullet}) \simeq h_\ast g'_\ast \O_{\X_\bullet} \simeq h_\ast g'_\ast f'^\ast \O_{\Spec \MUP^{\otimes(\bullet+1)}} \simeq h_\ast f^\ast g_\ast \O_{\Spec \MUP^{\otimes(\bullet+1)}}\]
    \[\simeq h_\ast (\O_\X \otimes \MUP^{\otimes(\bullet+1)}) \simeq h_\ast(\O_\X) \otimes \MUP^{\otimes(\bullet+1)} \simeq \Ga(\O_\X) \otimes \MUP^{\otimes(\bullet+1)},\]
    the first equivalence coming by factoring $\X_\bullet \to \Spec \Sph$, the second as $f'^\ast$ is strong symmetric monoidal, the third as the square above is adjointable, the fourth by \Cref{relativespectrumfiberofaffine}, the fifth as $h_\ast$ preserves colimits as $\X$ is $0$-semiaffine, and the sixth by definition. As $g$ is an atlas,  $g'$ is an atlas, and we obtain the desired equivalence:
    \[\X \simeq \colim \X_\bullet \simeq \colim \Spec \Ga(\O_\X) \otimes \MUP^{\otimes(\bullet+1)} \simeq \colim \Spec (\Ga(\O_\X) \otimes \MUP)^{\otimes_{\Ga(\O_\X)}(\bullet+1)} \simeq \M_{\Ga(\O_\X)}^\ori,\]
    the last equivalences coming from the fact that $\M_{\Ga(\O_\X)}^\ori$ is the descent stack associated to $\Ga(\O_\X) \to \Ga(\O_\X) \otimes \MUP$ by \Cref{pr:obviouspropertiesofMAORi}(\ref{item:MAOridescentstack}) and the simplicial presentation of these descent stacks of \Cref{lem:descentstackbasicproperties}.
\end{proof}

One interpretation of \Cref{thm:reconstruction} is that the global sections functor $\Ga\colon \Stk^\op \to \CAlg$ is fully faithful when restricted to chromatically affine stacks which are $0$-affine. If we restrict to those chromatically affine stacks of a fixed height, then we can easily identify the essential image of $\Ga$ too.

\begin{cor}[{\Cref{maintheorem:reconstruction2}}]\label{thm:reconstructionbounded}
    Fix a bounded height function $\vec{h}$. Then the global sections functor
    \[\Ga\colon \Stk^{\aff}_{/\M_\FGvech^\ori} \to \CAlg(\Sp_{\vec{h}})^\op\]
    is an equivalence.
\end{cor}

\begin{proof}
    The global sections functor $\Ga$ splits as the composition
    \[\Stk^\aff_{/\M_\FGvech^\ori} \xrightarrow{\Phi} \CAlg(\QCoh(\M_\FGvech^\ori))^\op \xrightarrow{\Ga^\op} \CAlg(\Sp_{\vec{h}})^\op,\]
    where the first functor is the equivalence of \Cref{relativespectrumtheorem}, and the second is $\CAlg(-)$ applied to the equivalence of \Cref{thm:qcohofDvech}.
\end{proof}

\begin{cor}\label{cor:reconstructionounded}
    The global sections functor $\Ga\colon \Stk^\op \to \CAlg$ is fully faithful when restricted to chromatically affine stacks of bounded height. The essential image are those $\E_\infty$ rings which are $L_{\vec{h}}$-local for some bounded height function $\vec{h}$.
\end{cor}

\begin{proof}
    For each pair of chromatically affine stacks of bounded height $\X,\Y$, then by taking maxima we obtain a common bounded height function $\vec{h}$ such that both $\X$ and $\Y$ have height $\leq \vec{h}$. It then follows from \Cref{thm:reconstructionbounded} that the map of spaces
    \[\Map_{\Stk}(\Y,\X) \xrightarrow{\Ga} \Map_{\CAlg}(\Ga(\O_\X), \Ga(\O_\Y))\]
    is an equivalence. Essential surjectivity follows by \Cref{thm:reconstructionbounded}.
\end{proof}

As a more explicit example of how one might use reconstruction, we have the following corollary. This was a key part of the construction and uniqueness of the spectral Tate curve constructed by the second- and third-named authors in \cite{globaltate}.

\begin{cor}\label{cor:orientedeillipeiticcurves}
    For any complex-periodic\footnote{If $\X$ is not complex-periodic, then $\Map_{\Stk}(\X, \M_\Ell^\ori)$ is empty, however, the space $\Map_{\CAlg}(\TMF, \Ga(\O_\X))$ may not be: consider $\X = \Spec \KO\llpar q\rrpar$, for example, where there exists a \emph{smooth $q$-expansion map} $\TMF \to \KO\llpar q \rrpar$; see \cite[Rmk.6.6]{globaltate}. In fact, this particular map is explicitly constructed as the global sections of a map $\T\colon \Spec \KU\llpar q\rrpar /C_2 \to \M_\Ell^\ori$, defined by the \emph{spectral Tate curve}, rather than from $\Spec \KO\llpar q \rrpar$.} stack $\X$, the global sections functor induces an equivalence of spaces
    \[\Ell^\ori(\X)^\simeq = \Map_{\Stk}(\X, \M_\Ell^\ori) \xrightarrow{\simeq} \Map_{\CAlg}(\TMF,\Ga(\O_\X).\]
    In particular, an oriented elliptic curve over $\X$ is equivalent data to a map of $\E_\infty$ rings $\TMF \to \Ga(\O_\X)$.
\end{cor}

These sorts of statements were part of the motivation for this article.

\begin{proof}
As $\M_\Ell^\ori$ is reconstructible by \Cref{periodicisaffine} and \Cref{thm:reconstruction}, this follows from \cref{prop:ffforreconstruction}.
\end{proof}

Once a stack satisfies some form of reconstruction, taking pullbacks over that stack is also a purely algebraic process.

\begin{cor}\label{genrealizationofpullbacks}
    Let $\X$ be a chromatically affine stack of bounded height. Then the pullback of a span $\Spec \A \to \X \gets \Spec \B$ is equivalent to $\Spec \A\otimes_{\Ga(\X)} \B$.
\end{cor}

\begin{proof}
    Just as in the proof of \Cref{fiberofaffine}, this follows from the universal property of a pullback, as each $\E_\infty$ ring $\R$ admitting a map to $\X$ satisfies the hypotheses of \Cref{thm:reconstructionbounded}, leading to the equivalences
    \[\Map_{\CAlg}(\A \otimes_{\Ga(\X)} \B, \R) \simeq \Map_{\CAlg}(\A, \R) \underset{\Map_{\Stk}(\Spec \R, \X)}{\times} \Map_{\CAlg}(\B, \R).\qedhere\]
\end{proof}

For example, this statement allows one to compute pullbacks over $\M_\Ell^\ori$:
\begin{equation}\lim (\Spec \A \to \M_\Ell^\ori \gets \Spec \B) \simeq \Spec (\A \otimes_{\TMF} \B)   \end{equation}


\addcontentsline{toc}{section}{References}
\bibliography{reference}
\bibliographystyle{alpha}
\end{document}